\documentclass[a4paper,11pt]{article}
\usepackage[top=3.5cm,bottom=2.5cm,left=2.5cm,right=2.5cm,a4paper]{geometry}
\usepackage{amsfonts,amsmath,amssymb,amsthm}
\usepackage[colorlinks,citecolor=green!60!black]{hyperref} 
\usepackage{float,xcolor,fancyhdr} 
\usepackage{algorithmic,algorithm2e,setspace,float}
\usepackage{ifthen}
\usepackage{subfigure}
\usepackage{natbib} 
\usepackage{enumitem}
\usepackage{mdframed}
\usepackage{soul}

\usepackage{tikz,pgfplots,todonotes}
\usepgfplotslibrary{groupplots}
\usetikzlibrary{positioning,arrows,fit}
\pgfplotsset{compat=newest}
\pgfplotsset{select coords between index/.style 2 args={
    x filter/.code={
        \ifnum\coordindex<#1\fi
        \ifnum\coordindex>#2\fi
    }
}}

\newcommand{\num}[1]{\text{#1}}
\newcommand{\si} [1]{\text{#1}}
\newcommand{\Hz}    {\text{Hz}}
\newcommand{\meter} {\text{m}}
\newcommand{\km}    {\text{km}}

\newcommand{\decibel}{\text{dB}}
\newcommand{\persecond}       {$\,$\text{s}$^{-1}$}

\newcommand{\second}{\text{s}}
\newcommand{\giga}  {\text{G}}

\newcommand{\update}[1]{#1}

\newcommand*\samethanks[1][\value{footnote}]{\footnotemark[#1]}
\newcommand{\myblue} {blue!60!white}

\newtheorem {remark}     {Remark}

\newcommand{\forward}     {\mathcal{F}}
\newcommand{\R}           {\mathbb{R}}
\newcommand{\laplacian}   {\Delta}
\newcommand{\divergence}  {\nabla \cdot}
\newcommand{\ii}          {\mathrm{i}}
\newcommand{\dd}          {\mathrm{d}}
\newcommand{\bx}          {\boldsymbol{x}}

\newcommand{\misfit}      {\mathcal{J}}
\newcommand{\regul}       {\mathcal{I}}
\newcommand{\pressure}    {p}

\newcommand{\data}        {\boldsymbol{y}} 
\newcommand{\pde}       {\mathcal{A}}
\newcommand{\err}       {\mathcal{E}}
\newcommand{\ngrad}[1]  {g_{{#1}}(m)} 
\newcommand{\ngradx}[1] {g_{{#1}}(m,\bx)} 
\newcommand{\basis}     {{\boldsymbol{\psi}}} 
\newcommand{\basisloc}  {\basis_{\text{loc}}}
\newcommand{\decomp}    {{\mathfrak  {m}}}    

\newlength{\modelwidth}     \newlength{\modelheight}
\newlength{\plotwidth}      \newlength{\plotheight}
\newlength{\plotlinewidth}  
\newlength{\plotmarkwidth}  
\newcommand{\modelfile} {}
\newcommand{\datafolder}{}
\newcommand{\datafile}  {}
\newcommand{\dataA}{} \newcommand{\dataB}{} 
\newcommand{\dataC}{} \newcommand{\dataD}{} \newcommand{\dataE}{} 
\newcommand{\modelfileA} {}
\newcommand{\modelfileB} {}
\newcommand{\modelfileC} {}
\newcommand{\modelfileD} {}
\newcommand{\modelfileE} {}
\newcommand{\modelfileF} {}
\newcommand{\modelfileG} {}
\newcommand{\modelfileH} {}
\newcommand{\modelfileI} {}
\newcommand{\dataformula}{}
\newcommand{\tablerrseam}[6]{ #1 & #2 & #3 & #4 & #5 & #6}

\title{Eigenvector Models for Solving the Seismic Inverse Problem for 
       the Helmholtz Equation: Extended Materials}

\author{  \fontsize{10pt}{10pt}
Florian Faucher\thanks{Inria Project-Team Magique 3D, E2S UPPA, CNRS, Pau, France.} 
$^,$\thanks{Faculty of Mathematics, University of Vienna,Oskar-Morgenstern-Platz 1,
            A-1090 Vienna, Austria.} 
$^,$\thanks{\href{mailto:florian.faucher@univie.ac.at}
                                  {\texttt{florian.faucher@univie.ac.at}}.}
~, Otmar Scherzer\samethanks[2]
               $^,$\thanks{Johann Radon Institute for Computational
                           and Applied Mathematics (RICAM), 
                           Linz, Austria.}
~and H\'el\`ene Barucq\samethanks[1]
}
\date{}

\begin{document}
\maketitle

\begin{abstract}
  We study the seismic inverse problem 
  for the recovery of subsurface properties 
  in acoustic media.
  In order to reduce the ill-posedness of the 
  problem, the heterogeneous wave speed parameter
  is represented using a limited 
  number of coefficients associated 
  with a basis of eigenvectors of a diffusion equation, 
  following the \emph{regularization by discretization} 
  approach. 
  We compare several choices for the diffusion coefficient
  in the partial differential equations, 
  which are extracted from the field of image processing.
  We first investigate their efficiency for image 
  decomposition (accuracy of the representation 
  with respect to the number of variables and denoising).
  Next, we implement the method in the quantitative 
  reconstruction procedure for seismic imaging, 
  following the Full Waveform Inversion 
  method, where the difficulty resides in that the  
  basis is defined from an initial model 
  where none of the actual structures is known.
  In particular, we demonstrate that the method is 
  efficient for the challenging reconstruction of 
  media with salt-domes.
  We employ the method in two and three-dimensional experiments, 
  and show that the eigenvector representation 
  compensates for the lack of low-frequency information,
  it eventually serves us to extract guidelines for the 
  implementation of the method.
\end{abstract}

\renewcommand{\thefootnote}{\arabic{footnote}}
\setcounter{footnote}{0}
\section{Introduction}

We consider the inverse problem associated with the 
propagation of time-harmonic waves which occurs, for 
example, in seismic applications, where the mechanical 
waves are used to probe the Earth. 
Following the non-intrusive geophysical setup
for exploration, we work with measured seismograms 
that record the waves at the surface (i.e., 
partial boundary measurements) and one-side 
illumination (back-scattered/reflection data).
In the last decades, this problem has encountered a 
growing interest with the increase in numerical 
capability and the use of supercomputers. 
However, the accurate recovery of the deep subsurface structures 
remains a challenge, due to the nonlinearity and 
ill-posedeness of the problem, the availability of 
partial reflection data only, and the large scale 
domains of investigation.

In the context of seismic, the quantitative reconstruction 
of physical properties using an iterative minimization 
of a cost function originally follows the work of
\cite{Bamberger1979}, \cite{Lailly1983} and \cite{Tarantola1984,Tarantola1987a}
in the time-domain, Pratt et al. (\citeyear{Pratt1990,Pratt1996,Pratt1998})
for the frequency approach. 
The method is commonly referred to as 
\emph{Full Waveform Inversion} (FWI),
which takes the complete observed seismograms for data.
One key of FWI is that the gradient of the misfit functional
is computed using the adjoint-state method~\citep{Lions1971,Chavent1974},
to avoid the formation of the (large) Jacobian matrix;
we refer to \cite{Plessix2006} for a review in geophysical applications. 
Then, Newton-type algorithms represent the traditional
framework to perform the iterative minimization. 
Due to the large computational scale of 
the domain investigated, seismic experiments may have
difficulties to incorporate second order (Hessian)
information in the algorithm, and 
alternative techniques have been proposed, for example, 
in the work of \cite{Pratt1998,Akcelik2002,Choi2008,Metivier2013,Jun2015}.
The quantitative (as opposed to qualitative) reconstruction 
methods based upon iterative minimization are naturally 
not restricted to seismic and we refer, among others, 
to \cite{Ammari2015visco,Barucq2018} and the references 
therein for additional applications using similar techniques.

The main difficulty of FWI 
(in both exploration and global seismology) lies in the 
high nonlinearity of the problem and the presence of 
local minima in the misfit functional, which are due to 
the time shifts and cycle-skipping effect \citep{Bunks1995}, in particular 
when the background velocity (the low-frequency profile) 
is not correctly anticipated 
\citep{Gauthier1986,LuoSchuster1991,Fichtner2008,Faucher2019IP}.
For this reason, the phase information is included in the 
traveltime inversion by \cite{LuoSchuster1991} by using
a cross-correlation function between measurements and 
simulations, where the relative phase shift is given 
by the maximum of the correlation. 
The method is further generalized by \cite{GeeJordan1992}, 
while \cite{VanLeeuwen2010} propose to select the 
phase shift using a weighted norm.
The choice of misfit has further encountered a growing 
interest in the past decade: in application to 
global-scale seismology, \cite{Fichtner2008} compare the 
phase, correlation-based, and envelope misfit functionals, 
the latter being also studied by \cite{Bozdag2011}.
In exploration seismic, comparisons of phase and amplitude 
inversion are performed by \cite{Shin2007b,Shin2007c}.
The $L1$ norm is studied by \cite{Brossier2010} while
approaches based upon optimal transport are considered 
by \cite{Metivier2016,Yang2018}. 
In the context where different fields are measured, 
\cite{Alessandrini2019,Faucher2019FRgWI} advocate 
for a reciprocity-based functional, which further
connects to the correlation-based formulas \citep{Faucher2019FRgWI}.
In the case of accurate knowledge of the background 
velocity, the inverse problem is close to linear or 
quasi-linear as the Born approximation holds and then, 
alternative methods of linear inverse problem can be 
applied, such as the Backus--Gilbert method \citep{BackusGilbert1967,BackusGilbert1968}. 
The difficulty to recover the background velocity variation
has also motivated alternative parametrization of the inverse 
problem: for instance the MBTT (Background/Data-Space Reflectivity) 
reformulation of FWI \citep{Clement2001,Faucher2019IP}. 

In order to diminish the ill-posedness of the inverse 
problem, a \emph{regularization} criterion can be 
incorporated.
It introduces an additional constraint
(in addition to the fidelity between observations 
and simulations), which, however, may be complicated 
to select a priori and problem dependent (with `tuning' parameters).
For instance, we refer to the body 
of work of \cite{Kirsch1996,Isakov2006,Kern2016,Kaltenbacher2018}
and the references therein. 
In the \emph{regularization by discretization} approach,
the model representation plays the role 
of regularizing the functional, by controlling 
(and limiting) the number of unknowns in the problem,
and possibly defining (i.e. constraining) the shape 
of the unknown (e.g., to force smoothness). 
Controlling the number of unknowns influences 
the resolution of the outcome, but also the 
stability and convergence of the procedure.
The use of piecewise constant coefficients appears natural 
for numerical applications, and is also motivated by 
stability results \citep{Alessandrini2005,Beretta2016}. 
However, such a decomposition can lead to an artificial `block' 
representation (cf.~\cite{Beretta2016,Faucher2017}) which 
would not be appropriate in terms of resolution. 
For this reason, a piecewise linear model representation is 
explored by \cite{Alessandrini2018,Alessandrini2019}, still 
motivated by the stability properties. 
We also mention the wavelet-based model reductions, 
that offer a flexible framework and are used for the purpose 
of regularization in seismic tomography 
by \cite{Loris2007,Loris2010}.
In the work of Yuan et al. (\citeyear{Yuan2014,Yuan2015}), FWI is carried out in the
time-domain with a model represented from a 
wavelet-based decomposition. 

In our work, we will use a model decomposition based 
upon the eigenvectors of a chosen diffusion operator,
as introduced by \cite{DeBuhan2013,Grote2017,Grote2019}. 
Note that this decomposition is shown  
(with the right choice of operator) to be related with the
more standard Total Variation (TV) or Tikhonov regularizations. 
The main difference in our work is that we study
several alternatives for the choice of the operator
following image processing techniques, which 
traditionally also relies on such diffusion PDEs 
(e.g.~\cite{Weickert1998}).
We first investigate the performance of the 
decomposition depending on the choice of PDE, 
and, next, the performance of such a model 
decomposition as parametrization of the 
reconstruction procedure in seismic FWI. 
It shows that the efficient choice of PDE should 
change depending on the situation. In addition, 
we provide a series of experiment to extract the 
robust guidelines for the implementation of the 
method in seismic.

We specifically target the reconstruction of 
subsurface \emph{salt domes} (i.e. media 
with high contrasting objects), which is 
particularly challenging, because (in addition 
to the usual restrictive data) of the change of 
the kinematics involved and the lack of low 
frequency data 
\citep{Farmer1996,Chironi2006,Virieux2009,Faucher2019RR}. 
In such cases, the use of the Total-Variation 
regularization \citep{Rudin1992} with FWI is 
becoming more and more prominent, and consists 
in incorporating an additional constraint on 
the model in the minimization problem.
Its efficiency is shown in the context of acoustic 
media with salt-dome contrasts by, e.g., 
\cite{Brandsberg2017,Esser2018,Kalita2019,Aghamiry2019}.
In our work, we study several alternatives and 
demonstrate that the model representation with 
the criterion extracted from \cite{Geman1992} 
appears the most appropriate for the eigenvector 
decomposition method in the presence of 
salt-domes.
We also show the limitations of the method,
in particular, it appears that the decomposition 
fails to represent models which are composed of 
several thin structures.

In Section~\ref{section:inverse_problem}, we 
define the inverse problem associated with the 
Helmholtz equation and introduce the iterative 
method for the reconstruction of the wave speed. 
In Section~\ref{section:model_decomposition}, 
we review several possibilities for the model 
decomposition using the eigenvectors of the 
diffusion operators. 
The process of model (image) decomposition is 
illustrated in Section~\ref{section:numeric_decomposition}. 
Then, in Section~\ref{section:numeric_reconstruction}, 
we carry out the iterative reconstruction with FWI 
experiments in two and three dimensions. 
Here, the model decomposition is based upon the 
initial model, which does not contain a priori 
information on the target, hence increasing the 
differences of performance depending on the 
selection of the basis. 
It allows us to identify the best candidate 
for the recovery of salt dome, and to extract 
some guidelines for applications in quantitative 
reconstruction.

\section{Inverse time-harmonic wave problem}
\label{section:inverse_problem}

\subsection{Forward problem}

We consider a domain $\Omega$ in two or three 
dimensions, with $\Omega \subset \R^2$ or $\Omega \subset \R^3$. 
We focus on acoustic media where, for simplicity, 
the density is taken as a constant, leading us to the 
identification of a single heterogeneous parameter: the wave speed.
The propagation of waves in an acoustic medium with
constant density is given by the scalar pressure 
field $\pressure$, solution to the Helmholtz equation
\begin{equation} \label{eq:helmholtz_equation_without_BC}
  \Big(-\laplacian -\omega^2 c^{-2}(\bx) \Big) 
          \pressure(\bx) = f(\bx), \quad \quad \text{in $\Omega$},
\end{equation}
where $c$ is the wave speed, $f$ the
source, and $\omega$ the angular frequency. 
We now have to specify the boundary conditions
to formulate the appropriate problem.

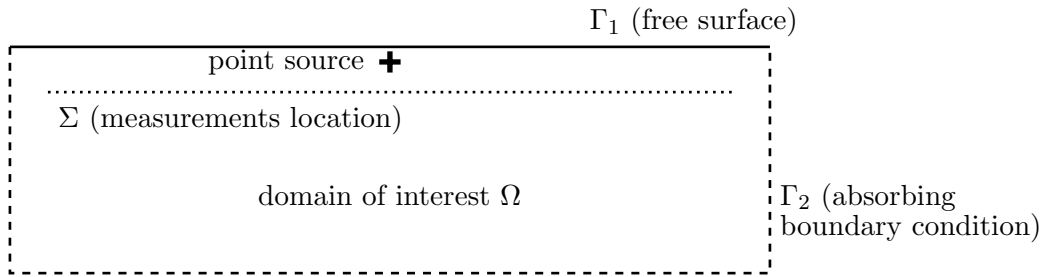
\begin{figure} \centering
\begin{tikzpicture}[]
  \coordinate (d1) at ( 0.,0.); \coordinate (d2) at (10.,0.);
  \coordinate (d3) at (10.,3.); \coordinate (d4) at ( 0.,3.);
  \coordinate (d5) at ( 0.,1.8);\coordinate (d6) at (10.,1.8);
  \coordinate (src)at ( 5.,2.8);
  \coordinate (r1) at ( .5,2.4); \coordinate (r2) at (9.5,2.4);

  \draw[line width=1,color=black,dashed] (d4) to (d1) to (d2) to (d3);
  \draw[line width=1,color=black]        (d3) []to[] (d4);
  \draw[mark=+,mark size=4,mark options={color=black,line width=2}] plot[] coordinates{(src)};
  \draw[line width=1,color=black,dotted]  (r1) []to[] (r2);

  \draw (d3) node[left=1. of d3,black ,anchor=south]{$\Gamma_1$ (free surface)};
  \draw (r1) node[below=0.7of r1 ,black ,anchor=south west]{$\Sigma$ (measurements location)};
  \draw (d3) node[below=2 of d3,anchor=west]{$\Gamma_2$ (absorbing};
  \draw (d3) node[below=2 of d3,anchor=west,yshift=-4mm]{boundary condition)};
  \draw (src)node[below=1.5 of src,anchor=north]{domain of interest $\Omega$};
  \draw (src)node[left =0.2 of src,anchor=east,black] {point source};

\end{tikzpicture}
 \caption{Illustration of the two-dimensional computational
          domain. The sign `plus' indicates the (point) source
          and the dotted line is the location of the 
          discrete measure points. The complete acquisition 
          consists in several sources, all are located at
          the same depth, only the lateral position
          varies; the measure points are fixed.
          }
\label{fig:sketch_domain}
\end{figure}

Following a seismic setup, the boundary of $\Omega$, $\Gamma$ 
is separated into two. The upper (top) boundary, $\Gamma_1$, 
represents the interface between the Earth and the air, and 
is subsequently represented via a \emph{free surface} boundary 
condition, where the pressure is null. 
On the other hand, other part of the boundary 
corresponds to the numerical need for restricting the area 
of interest. Here, conditions must ensure that entering waves 
are not reflected back to the domain 
(i.e., because the area of interest is only a part of the Earth),
see Figure~\ref{fig:sketch_domain}. 
The two most popular formulations to handle such numerical 
boundary are either the Perfectly Matched Layers 
(PML, \cite{Berenger1994}), or outgoing artificial boundary 
conditions.
In our case, we use Absorbing Boundary Conditions 
(ABC, \cite{Engquist1977}) so that the complete 
problem writes as
\begin{equation} \label{eq:helmholtz_equation}
    \left\lbrace \begin{aligned}
    \big(-\laplacian -\omega^2 c^{-2}(\bx) \big) \pressure(\bx) & = f(\bx), 
                                      \quad \quad & \text{in $\Omega$}, \\
    \pressure(\bx) &= 0 , \quad \quad & \text{on $\Gamma_1$}, \\
    \partial_\nu \pressure(\bx) - \ii \omega c^{-1}(\bx) \pressure(\bx)
    &= 0 , \quad \quad & \text{on $\Gamma_2$}, \\
    \end{aligned} \right.
\end{equation}
where $\partial_\nu$ is the normal derivative.
We recall that $\pressure$ denotes the pressure, 
$\omega$ denotes the frequency of the induced 
source and $c$ denotes the wave speed.

The inverse problem aims the recovery 
of the wave speed $c$ in~\eqref{eq:helmholtz_equation},
from a discrete set of measurements 
(i.e. partial data), which corresponds to 
observation of the wave propagations. 
More precisely, our data consist in measurements of 
the pressure field solution to 
\eqref{eq:helmholtz_equation}, at the (discrete) 
device locations.
We refer to $\Sigma$ for this set of positions, where
the receivers are located, and define the forward 
map at frequency $\omega$, $\forward_\omega$ 
(which links the model to the data), such that,
\begin{equation} \label{eq:forward_map}
 \forward_\omega ~:~  m \, \, \rightarrow \, \, \pressure(\bx)\mid_\Sigma.
\end{equation}
We have introduced $m(\bx) :=c^{-2}(\bx)$, which is 
also our choice of parameter for the reconstruction, 
see Remark~\ref{rk:parametrization}.
In seismic, the data are further generated from 
several point sources (excited one by one) and 
all devices (sources and receivers) remain near 
the surface.
All devices (sources and receivers) 
remain near the surface ($\Gamma_1$), as 
illustrated with $\Sigma$ in Figure~\ref{fig:sketch_domain}.

\subsection{Quantitative reconstruction using iterative minimization}

The inverse problem aims the reconstruction of the 
unknown medium squared slowness $m$ (i.e. the wave 
speed) from data $\data_\omega$ 
that connects to the forward map $\forward_\omega(m_\dagger)$
for a reference (target) model $m_\dagger$ with
\begin{equation}
  \data_\omega = \forward_\omega(m_\dagger) + \mathfrak{E}_\omega,
\end{equation}
where $\mathfrak{E}_\omega$ represents the noise in the 
measurements (from the inaccuracy of the devices, model error, etc), 
possibly frequency dependent.

For the reconstruction, we follow the 
Full Waveform Inversion (FWI) method
\citep{Tarantola1984,Pratt1998}, which
relies on an iterative minimization of a misfit 
functional defined as the difference between the 
observed data and the computational simulations:
\begin{equation} \label{eq:minimization}
 \min_m \misfit(m) = \dfrac{1}{2} \sum_\omega 
                     \Vert \forward_\omega(m) - \data_\omega \Vert^2,
\end{equation}
where we use the standard least-squares minimization 
but alternatives have also been studied, as indicated 
in the introduction. 

\begin{remark}[Multi-frequency algorithm]\label{rk:frequency_progression}
 For the choice of frequency in Problem~\ref{eq:minimization}, 
 applications commonly use a sequence of increasing frequencies 
 during the iterative process
 \citep{Bunks1995,Pratt1990,Sirgue2004,Brossier2009,Barucq2018,Faucher2017}.
 Namely, one starts with a low frequency and minimize the functional 
 for the fixed frequency content in the data. 
 Once the misfit functional stagnates, or after
 a prescribed number of iterations, the frequency 
 is updated (increased) and the iterations continue, 
 cf. Algorithm~\ref{algo:fwi}.
 Moreover, the use of sequential frequency 
 (instead of band of frequencies)
 is advocated by \cite{Faucher2019RR}, because
 it enlarges the size of the basins of attraction.
\end{remark}

\begin{remark}[Parametrization of the unknown]\label{rk:parametrization}
 For the reconstruction, we invert the squared slowness
 $m = c^{-2}$ instead of the velocity. The choice
 of this parameter is first motivated by the 
 Helmholtz equation~\eqref{eq:helmholtz_equation}.
 However, it (i.e., velocity, slowness or squared slowness inversion)
 can lead to an important difference in the efficiency 
 of the reconstruction procedure. 
 It is discussed, for example, by \cite{Tarantola1986,Brossier2011,Kohn2012};
 in particular, we motivate our choice from the comparison of 
 reconstructions provided in the context of seismic 
 by~\cite[Section~5.4]{Faucher2017}.
\end{remark}

Then, an iterative minimization algorithm 
is used for the resolution of Problem~\ref{eq:minimization} 
in the framework of the Newton methods. 
Starting with an initial guess $m^{(0)}$, 
the model is updated at each iteration $k$, using 
a search direction $s^{(k)}$, such that
\begin{equation} \label{eq:model_update_general}
  m^{(k+1)} = m^{(k)} + \mu s^{(k)}, 
  \qquad k \, > \, 0.
\end{equation}
Several possibilities exist for the search direction 
(e.g., Newton, Gauss-Newton, BFGS, gradient descent, 
etc.) and we refer to \cite{Nocedal2006} for an extensive
review of the methods. The scalar coefficient  $\mu$ 
is approximated using a line search algorithm \citep{Nocedal2006}. 
In our implementation, we rely on a gradient-based optimization,
with the nonlinear conjugate gradient method, and a backtracking 
line search \citep{Nocedal2006}. 
A review of the performance of first order-based 
minimization algorithms and the influence 
of line search step selection is further 
investigated by~\cite{Barucq2018} in the context 
of inverse scattering.
The computation of the gradient of the misfit
functional is carried out using the adjoint-state method
\citep{Chavent1974,Plessix2006}, which specific steps 
for complex-valued fields can be found in 
\cite[Appendix~A]{Faucher2019IP}

\section{Regularization by discretization: model decomposition}
\label{section:model_decomposition}

In this section, we introduce a representation of 
the unknown, i.e., the \emph{model decomposition}, based
upon the eigenvectors of a diffusion equation.
The objective is to reduce the dimension of the 
unknown to mitigate the ill-posedness of the inverse
problem. 
We provide several possibilities for the choice of
the eigenvectors, following the literature in image 
processing.

\subsection{Regularization and diffusion operators}

The resolution of the inverse problem using a 
quantitative method introduces an optimization 
problem~\eqref{eq:minimization} where the misfit functional 
$\misfit$ accounts for a \emph{fidelity} term.
The resolution of such a problem is 
also common in the context of image processing 
(e.g., for denoising or edge enhancement)
where the fidelity term corresponds to the 
matching between the original and processed images.
It is relatively common (for both the quantitative 
reconstruction methods and in image processing) 
to incorporate an additional term in the minimization, 
for the purpose of regularization. 
The primary function  of this additional term is to reduce 
the ill-posedness of the problem, by adding a constraint.
It has been the topic of several studies, we refer to, 
e.g., \cite{Kirsch1996,Isakov2006,Charbonnier1994,Robert1996,
      Rudin1992,Vogel1996,Lobel1997,Kern2016,Qiu2016,Kaltenbacher2018}.
The \emph{regularized} minimization problem writes as 
\begin{equation} \label{eq:minimization_with_regularization}
 \min_m \misfit_r(m) = \dfrac{1}{2} \sum_\omega 
                     \Vert \forward_\omega(m) - \data_\omega \Vert^2
                   + \regul(m),
\end{equation}
where $\regul$ stands for the regularization term.

In many applications such as image processing, 
$\regul$ is usually defined to only depend on the 
gradient of the variable (image), such that 
\begin{equation} \label{eq:regularization_term}
  \regul(m) = 
  \int_\Omega \phi(\vert \nabla m \vert) \, \, \dd \Omega \, ,
\end{equation}
where $\phi \in L^2(\Omega)$.
In particular, the minimum of $\regul$ with respect to
$m$ verifies the Euler--Lagrange equations \citep{Evans2010,Dubrovin1992}.
In one dimension, it is given by \cite[Theorem~31.1.2]{Dubrovin1992}
and is extended for higher dimensions with~\cite[Theorem~37.1.2]{Dubrovin1992}
(further simplified in our case because $\regul$ only depends on the gradient).
It states that the minimizer of $\regul$ is the solution of the diffusion equation:
\begin{equation} \label{eq:diffusion_main}
  \divergence \bigg( \dfrac{\phi'(\vert \nabla m \vert)}
                          {\vert \nabla m \vert} \nabla m \bigg) = 0 
                    \quad \quad \text{in $\Omega$.} 
\end{equation}
For the sake of clarity, we introduce the following notation:
\begin{equation} \label{eq:diffusion_operator}
  \pde(\mathrm{y}, \, \eta) := -\divergence \big( \eta(\mathrm{y}) \nabla \big),
  \quad \quad \quad \text{with} \quad \quad 
  \eta(\mathrm{y}) = \dfrac{\phi'(\vert \nabla \mathrm{y} \vert)}{\vert \nabla \mathrm{y} \vert}.
\end{equation}
In the following, we present several choices for
the diffusion PDE coefficient $\eta$, following
image processing theory.

\begin{remark}
  The minimization of $\misfit_r$ 
  in Problem~\eqref{eq:minimization_with_regularization}
  can be performed using traditional gradient 
  descent or Newton type algorithms. 
  Another alternative, in particular when rewriting 
  with the Euler--Lagrange formulation in the context of 
  image processing, is to recast the problem as a time 
  dependent evolution one, see, e.g., the work of
  \cite{Weickert1998,Catte1992,Rudin1992,Alvarez1992}.
\end{remark}

\begin{remark}
  The diffusion equation~\eqref{eq:diffusion_main} is 
  obtained using the fact that $\phi = \phi(\vert \nabla m \vert)$ 
  only depends on $\nabla m$. In case of dependency of
  the function with $m$, or higher order derivatives, 
  the Euler--Lagrange formulation must be adapted.
\end{remark}

\subsection{Diffusion coefficients from image processing}

There exist several possibilities for the choice
of diffusion coefficient $\eta$ (also referred to
as the \emph{weighting function}) in~\eqref{eq:diffusion_operator},
inherited from image processing theory and applications. 
In the following, we investigate the most common formulations, 
see Table~\ref{table:list_of_formula}, for which we have mainly
followed the ones that are reviewed by 
\cite[Table~1]{BlancFeraud1995} and~\cite{Robert1996}.
Furthermore, we incorporate a scaling 
coefficient $\beta > 0$ for the diffusion coefficient,  
which impacts on the magnitude. 
For consistency in the different models, 
the norms that we employed are 
scaled with the maximal values so that they 
remain between $0$ and $1$. We define,
\begin{equation} 
\begin{aligned}
&
\ngradx{1} = \dfrac{\vert \nabla m(\bx) \vert}  {\gamma_1}, \qquad 
\ngradx{2} = \dfrac{\vert \nabla m(\bx) \vert^2}{\gamma_2}, \qquad 
\text{with} \\
&
\gamma_1 = \max \big( \vert \nabla m(\bx) \vert \big)    \, , \quad
\gamma_2 = \max \big( \vert \nabla m(\bx) \vert^2 \big)  \, , \quad
\vert \nabla m(\bx) \vert
   = \sqrt{\sum_{k=1}^\mathfrak{d} \bigg( \dfrac{\partial m}{\partial x_k} \bigg)^2},
\end{aligned}
\end{equation}
where $\mathfrak{d}$ is the space dimension
($\mathfrak{d}=2$ or $\mathfrak{d}=3$ in our 
experiments) and $\bx$ the space coordinates: 
$\bx = \{x_k\}_{k=1}^\mathfrak{d}$.
In order to simplify the formulas, we will omit 
the space dependency in the following.
Note that in the numerical experiments, we calculate 
the eigenvalues and eigenvectors from the \emph{linear} 
differential operator $\mathcal{A}(m,\eta)$ defined 
in~\eqref{eq:diffusion_operator}, where the 
diffusion coefficient $\eta$ is taken from the
\emph{nonlinear} PDE model.

\begin{remark}
  We can make the following comments regarding 
  the nine diffusion coefficients that are 
  introduced in Table~\ref{table:list_of_formula}.
  \begin{itemize}
    \item The PDE~\eqref{eq:diffusion_main} using the Tikhonov 
          diffusion coefficient $\eta_9$ coincides with the 
          Laplace equation.
    \item For the formulation of $\eta_4$ and $\eta_8$, we 
          have to impose a threshold as the coefficient is 
          not defined for the points where the gradient is
          zero. In the computations, we impose that 
          $\eta_4=\eta_8=1$ for the points $\bx_i$ where 
          $\nabla m(\bx_i) < \num{$10^{-12}$}$.
    \item The first Perona--Malik formula $\eta_1$, is very
          similar to the Lorentzian approach, $\eta_6$: 
          only the position of $\beta$ differs. 
          Namely, the Perona--Malik formula would rather use 
          small $\beta$ while the Lorentzian formula would use
          large $\beta$.
    \item The second Perona--Malik formula, $\eta_2$, is very
          similar to the Gaussian criterion $\eta_7$, which only
          includes an additional dependency on $\beta$.
    \item The formulation of $\eta_8$ corresponds
          to the Total Variation (TV) regularization \citep{Vogel1996}.
  \end{itemize}
\end{remark}

\begin{table} \centering
\caption{List of formula for the coefficient 
         of the diffusion operator which is used 
         to decompose the image.}
\label{table:list_of_formula}
\begin{tabular}{@{}l|lcc}
 \bf reference (name)          & \bf definition    & 
 \bf $\beta \rightarrow 0$     &
 \bf $\beta \rightarrow \infty$ \\ \hline \hline
 \cite{Perona1988,Perona1990}
 &
 $\eta_{1}(m,\beta) = \dfrac{\beta}{\beta + \ngrad{2}}$ &  0  & 1 \\ \hline
 \cite{Perona1988,Perona1990}
 &
 $\eta_{2}(m,\beta) = \exp\bigg(-\dfrac{\ngrad{2}}{\beta}\bigg)$ 
 &  0  & 1 \\ \hline
 \cite{Geman1992}
 &
 $\eta_{3}(m,\beta) = \dfrac{2\beta}{\big(\beta + \ngrad{2}\big)^2}$ 
 &  0  & 0 \\ \hline
 \cite{Green1990Image}
 &
 $\eta_{4}(m,\beta) = \tanh\bigg( \dfrac{\ngrad{1}}{\beta}  \bigg)
                        \bigg( \dfrac{1}{\beta \ngrad{1}}     \bigg)$
 &  $+\infty$  & 0 \\ \hline
\cite{Charbonnier1994}
 &
 $\eta_{5}(m,\beta) = \dfrac{1}{\beta}
                   \bigg( \dfrac{\beta + \ngrad{2}}{\beta} \bigg)^{-1/2}$
 & $+\infty$ & 0 \\ \hline  
 \cite{Grote2019}  \\ (Lorentzian) 
 &
 $\eta_{6}(m,\beta) = \dfrac{\beta}{\big( 1 + \beta \ngrad{2} \big)^2} $
 & 0 & 0 \\ \hline
 \cite{Grote2019} \\ (Gaussian) 
 &
 $\eta_{7}(m,\beta) 
 = \bigg( \beta \exp\bigg( \dfrac{\ngrad{2}}{\beta} \bigg) \bigg)^{-1}$
 & 0 & 0 \\ \hline  
 \cite{Rudin1992} \\ (Total Variation, TV) 
 &
 $\eta_{8}(m) = \dfrac{1}{\ngrad{1}}$ 
 &  n/a  & n/a \\ \hline
 Tikhonov                                  &
 $\eta_{9} = 1$ \phantom{$\dfrac{1}{2}$}
 &  n/a  & n/a \\ \hline
\end{tabular}
\end{table}


\subsection{Eigenvector model decomposition in FWI}

In our work, we employ the \emph{regularization 
by discretization} approach: instead of adding the
regularization term $\regul$ in the 
minimization problem, we remain with
Problem~\eqref{eq:minimization}, and use a
specific representation for the model (unknown). 
We follow the work of 
\cite{DeBuhan2010,DeBuhan2013,Grote2017,Grote2019}
with the ``Adaptive Inversion'' or 
``Adaptive Eigenspace Inversion'' method. 
Namely, the unknown is represented via a 
decomposition into the basis of eigenvectors 
computed from a diffusion PDE. 
The purpose is to control the number of 
unknowns in the representation, and consequently 
reduce the ill-posedness of the inverse problem. 
The decomposition uses the steps (given in \cite{Grote2017}) 
depicted in Algorithm~\ref{algo:decomposition_ev}.

\begin{algorithm}[ht!] \centering 
  \setstretch{.5}
\noindent\fbox{%
\parbox{.9\textwidth}{
    \textbf{Eigenvector decomposition:} 
        given an initial model $m(\bx)$, a 
        selected integer value $N > 0$, and the 
        selected diffusion coefficient $\eta$.
    \begin{enumerate}
      \item Compute $\mathfrak{m}_0$, the solution of 
            the linear PDE
            \begin{equation} \label{eq:methodEV_m0}
            \left\lbrace \begin{aligned}
               \pde(m,\eta) \mathfrak{m}_0 & = 0,   \quad \quad &\text{in $\Omega$,} \\
                            \mathfrak{m}_0 & = m,   \quad \quad &\text{on $\Gamma$.}
            \end{aligned} \right.
            \end{equation}
      \item Compute the subset of $N$ eigenfunctions $\{\psi_k\}_{k=1,\ldots,N}$
            which are associated to the $N$ \emph{smallest} eigenvalues
             $\{\lambda_k\}_{k=1,\ldots,N}$ such that, for all $k$,
            \begin{equation} \label{eq:methodEV_eigenvector}
            \left\lbrace \begin{aligned}
               \pde(m,\eta) \psi_k & = \lambda_k \psi_k, \quad \quad &\text{in $\Omega$,} \\
                            \psi_k & = 0,                \quad \quad &\text{on $\Gamma$.}
            \end{aligned} \right.
            \end{equation}
      \item Compute the model decomposition using $N$ eigenvectors:
            \begin{equation} \label{eq:methodEV_decomposition}
              \mathfrak{m} = \mathfrak{m}_0 + \sum_{k=1}^N \alpha_k \psi_k,
            \end{equation}
            where $\alpha_k$ is a scalar and $\psi_k$ 
            a vector. Here, the set of $\alpha$ is chosen 
            to minimize $\Vert \mathfrak{m} - m \Vert^2$; 
            and the 
            $\psi_k$, $k=1,\ldots,N$ are the 
            eigenvectors associated with the $N$ smallest
            eigenvalues $\lambda_k$, computed in Step (ii).
    \end{enumerate}
}}
 \BlankLine 
 \caption{The model decomposition using the eigenvectors 
          of the diffusion equation associated with the 
          smallest eigenvalues. 
          We refer to the model decomposition as 
          $\decomp(m,N,\basis)$ where $\basis(m,\eta,N)$ 
          is the set of eigenvectors, 
          see~\eqref{eq:def:basis_decomposition}
          and~\eqref{eq:def:model_decomposition}.
          }
 \label{algo:decomposition_ev}
\end{algorithm}


Following Algorithm~\ref{algo:decomposition_ev}, 
we introduce the notation, 
\begin{equation} \label{eq:def:basis_decomposition}
  \begin{aligned}
      \basis(m,\eta,N) = \{\psi_k\}_{k=1,\ldots,N} \quad
      &\text{the set of $N$ eigenvectors associated with the
             model $m$ and} \\
      &\text{th diffusion coefficient $\eta$,
             computed from~\eqref{eq:methodEV_m0}
             and \eqref{eq:methodEV_eigenvector};}
  \end{aligned}
\end{equation}
and 
\begin{equation} \label{eq:def:model_decomposition}
  \decomp(m,N,\basis) 
      \qquad \text{is the decomposition of the model $m$ using $N$ 
                  vectors from $\basis$,~\eqref{eq:methodEV_decomposition}.}
\end{equation}

Therefore, the model is represented via $N$ 
coefficients $\alpha$ in~\eqref{eq:methodEV_decomposition}
in the basis given by the diffusion operator. 
The reconstruction procedure follows an 
iterative minimization of~\eqref{eq:minimization}, and 
performs successive update of the coefficients $\alpha$.
The key is that $N$ is much smaller than the dimension 
of the original representation of $m$, but allows 
an accurate resolution, as we illustrate in
Sections~\ref{section:numeric_decomposition} 
and~\ref{section:numeric_reconstruction}. 
Algorithm~\ref{algo:fwi} details the procedure.

\begin{algorithm}[ht!]
\setstretch{1}
\textbf{Inputs:}
  the measurements $\data_\omega$; 
  the initial model $m^{(0)}$;
  the selected number of iterations $n_{\text{iter}}$; 
  the list of frequencies $\omega_i$, $i=1$,$\ldots$,$n_\omega$;
  the decomposition dimension associated with frequency:
  $N_i$, for $i=1,\ldots,n_\omega$.
\vspace*{.2\baselineskip}
\noindent\fbox{%
\parbox{.9\textwidth}{
  Using Algorithm~\ref{algo:decomposition_ev},
  \begin{enumerate}
    \item compute the eigenvector basis for $m^{(0)}$
          using the selected $\eta$ and the highest 
          integer $N_{\max} = \max\big(\{N_i\}_{i=1}^{n_\omega}\big)$,
          \begin{equation}
            \basis =\basis(m,\eta,N_{\max}).
          \end{equation}
    \item Decompose the initial model using the initial 
          decomposition dimension:
          \begin{equation}
          \decomp^{(0)} = \decomp(m^{(0)},N_1,\basis)
                = \mathfrak{m}_0 + \sum_{l=1}^{N_1} \alpha_l^{(0)} \psi_l.
          \end{equation}
  \end{enumerate}
}}
\textbf{Frequency loop} \For{$i  \in \{1, \ldots, n_\omega \}$}{

  \vspace*{.25\baselineskip}
  Set $N=N_i$.
  \vspace*{.25\baselineskip}
    
  \textbf{Optimization loop} \For{$j \in \{1, \ldots, n_{\text{iter}} \}$}{
    \begin{description}
      {\setstretch{1.1}
      \item Set $k := (i-1)n_{\text{iter}} + j - 1$.
      \item Solve the Helmholtz equation~\eqref{eq:helmholtz_equation} 
            at frequency $\omega_i$ with model $\mathfrak{m}^{(k)}$.
      \item Compute the misfit functional $\misfit$ in~\eqref{eq:minimization}.
      \item Compute the gradient of the misfit functional $\nabla_\alpha\misfit$
            using the adjoint-state method.
      \item Compute the search direction $s^{(k)}$, see Remark~\ref{rk:minimization_tech}.
      \item Compute the descent step $\mu$ using the line search algorithm,
            see Remark~\ref{rk:minimization_tech}.
      \item Update the coefficient $\alpha$ with 
            $\alpha^{(k+1)} = \alpha^{(k)} - \mu s^{(k)}$.
      \item Update the model:
            \begin{equation*} 
              \decomp^{(k+1)} = 
              \decomp_0 + \sum_{l=1}^N \alpha_l^{(k+1)} \psi_l.
            \end{equation*} 
      }
    \end{description}
  }
} 
 \BlankLine 
\caption{The iterative minimization algorithm (FWI) using the model
         decomposition to control the number of variables. 
         The model is represented in the eigenvector basis, 
         for which the weights are updated along with the 
         iterations.}
 \label{algo:fwi}
\end{algorithm}


\begin{remark}[Minimization algorithm] 
  \label{rk:minimization_tech}
  For the minimization procedure depicted in 
  Algorithm~\ref{algo:fwi}, we use a non-linear 
  conjugate gradient method for the search 
  direction. 
  This method has the advantage that it only 
  necessitates the computation of the gradient 
  of the cost function \citep{Nocedal2006}.
  Then, to control the update step $\mu$ in 
  Algorithm~\ref{algo:fwi}, a line search 
  algorithm is typically employed 
  \citep{Eisenstat1994,Nocedal2006,Chavent2015,Barucq2018}.
  This operation is complex in practice 
  because an accurate estimation would require 
  intensive computational operations (with an 
  additional minimization problem to solve). 
  Here,  we employ a simple backtracking 
  algorithm \citep{Nocedal2006}.
\end{remark}

\begin{remark}[Gradient computation]
  \label{rk:gradient_computation}
  The gradient of the cost function is 
  computed using the first order adjoint-state 
  method \citep{Lions1971,Chavent1974},
  which is standard in seismic application
  \citep{Plessix2006}. It avoids the formation 
  of a dense Jacobian matrix and instead requires 
  the resolution of an additional PDE, 
  which is the adjoint of the
  forward PDE, with right-hand sides defined
  from the difference between the measurements 
  and the simulations, 
  see~\cite{Plessix2006,Faucher2017,Barucq2018,Faucher2019IP}
  for more details.

  In our implementation, the gradient is first computed
  with respect to the original (nodal) representation 
  and we use the chain rule to retrieve the gradient with
  respect to the decomposition coefficients $\alpha$:
  \begin{equation}
    \dfrac{\partial \misfit}{\partial \alpha} = 
        \dfrac{\partial \misfit}{\partial m} \dfrac{\partial m}{\partial \alpha}.
  \end{equation}
  It is straightforward, from~\eqref{eq:methodEV_decomposition},
  that the derivation for a chosen coefficient ${\alpha_l}$ 
  gives $\partial_{\alpha_l} m = \psi_l$. Therefore, it 
  is computationally easy to introduce the formulation
  with respect to the eigenvector decomposition from an 
  already existing `classical' 
  (i.e. when the derivative with respect to the model is 
  performed) formulation: it only
  necessitates one additional step with the eigenvectors.
\end{remark}

\begin{remark}
  In Algorithm~\ref{algo:fwi}, the basis of 
  eigenvectors remains the same for the 
  \emph{complete set of iterations}, and is 
  extracted from the initial model. Only the 
  number of vectors taken for the representation,
  $N_i$, changes. Namely, from $N_1$ to 
  $N_2 > N_1$, the decomposition using $N_2$ 
  still has the same $N_1$ first eigenvectors
  in its representation (with different weights 
  $\alpha$), and additional $(N_2 - N_1)$ eigenvectors.
  As an alternative, we investigate
  the performance of an algorithm where the 
  basis changes at each frequency (i.e. it is 
  recomputed from the current iteration model), 
  see Appendix~\ref{appendix:multi_basis_experiments}.
\end{remark} 

\subsection{Numerical implementation}

Our code is developed in \texttt{Fortran90}, it uses both 
\texttt{mpi} and \texttt{OpenMP} parallelism and run on 
cluster\footnote{The experiments have been performed
                 on the cluster PlaFRIM 
                (Plateforme F\'ed\'erative pour la 
                 Recherche en Informatique et Math\'ematiques,
                 \url{https://www.plafrim.fr/fr})
                 with the following node specification:
                 2 Dodeca-core Haswell Intel Xeon E5--2680 v3 (2.5\si{\giga\Hz});
                 128\si{\giga}o RAM; Infiniband QDR TrueScale: 
                 40\si{\giga}b\si{\persecond}, 
                 Omnipath 100\si{\giga}b\si{\persecond}.}
for efficiency.
The forward wave operator is discretized using a Finite 
Differences scheme, e.g.~\cite{Virieux1984,Operto2009,Wang2011}. 
The discretization of the Helmholtz operator generates a 
large sparse matrix, for which we use the direct solver 
\textsc{Mumps} (\cite{Amestoy2001,Amestoy2006}) for its factorization 
and the resolution of linear system. 
This solver is particularly optimized and designed for 
this type of linear algebra problems, i.e. large, sparse matrices.
Our preference for direct solver instead of iterative
ones is mainly motivated by two reasons:
\begin{enumerate}
 \item seismic acquisition is composed of a large amount of sources, 
       i.e. a large amount of right-hand sides (rhs) to be processed 
       for the linear system. Using direct solver, the resolution time
       is very low once the factorization is performed, hence it is 
       well adapted for the multi-rhs seismic configuration.
 \item For the minimization algorithm, the gradient is computed 
       via the adjoint-state method (see 
       Remark~\ref{rk:gradient_computation}).
       It means that an additional linear system has to be solved,
       which is actually the adjoint of the forward one. 
       Here, the factors obtained from the factorization of the 
       forward operator can be directly reused, and allow a 
       reduced computational cost, see~\cite{Barucq2018} and
       the references therein.
\end{enumerate}

The next step is the computation of 
the eigenvectors associated with the smallest eigenvalues 
for the diffusion operator. We use the package 
\textsc{Arpack}\footnote{\url{www.caam.rice.edu/software/ARPACK/},
  \textsc{Arpack} uses sequential computation, hence, contrary to
  the rest of our code, this part does not use parallelism. Future
  developments include the implementation of the parallel version
  of the package: \textsc{Parpack}.}, 
which is devoted to solve large sparse eigenvalue problems 
using iterative methods. More precisely, it uses implicitly 
restarted Lanczos or Arnoldi methods, respectively for symmetric 
and non-symmetric matrices, \cite{Lehoucq1996}. 
Several options are available in the package, including the 
maximum number of iterations allowed, or a tolerance parameter
for the accuracy of acceptable solution\footnote{
    We have observed important reduction of time cost when allowing 
    some flexibility in the accuracy with this threshold criterion. 
    However, in the computational experiments, we do not use this 
    option, as the numerical  efficiency is not the primary objective 
    of our study.
}.

\begin{remark}[Eigenvectors associated with the lowest eigenvalues]
  The Lanczos and Arnoldi methods are particularly efficient
  to compute the largest eigenvalues and associated eigenvectors 
  of the matrix, and only require matrix vector multiplication. 
  However, we are interested in the lowest eigenvalues for 
  our decomposition.
  The idea is simply to use that the lowest eigenvalues
  of the discretized diffusion matrix, say $A$, 
  are simply the largest eigenvalues of the matrix $A^{-1}$.
  Then, the matrix-vector multiplication, say $Av$
  for a vector $v$, becomes a resolution of a 
  linear system $A^{-1}v$. 
  It may appear computationally expensive but it 
  is not thanks to the use of the direct solver \textsc{Mumps}
  (see above), which, once the factorization is obtained, is 
  very efficient for the resolution procedure. Hence, the computation 
  follows the steps\footnote{\textsc{Arpack} 
          has the possibility to compute the smallest 
          eigenvalues using matrix-vector multiplication,
          however, we have observed a drastic increase of the 
          computational time compared to using the inverse
          matrix and resolution of linear system.}:
  \begin{enumerate} 
    \item compute the (sparse) matrix discretization 
          of the selected diffusion operator: $A$;
    \item compute the factorization of the matrix $A$
          using \textsc{Mumps},
    \item use the package \textsc{Arpack} to compute the 
          largest eigenvalues of $A^{-1}$, by replacing
          the matrix-vector multiplication step in the 
          iterations by the resolution of a linear 
          system using \textsc{Mumps}.
  \end{enumerate}
\end{remark}

\begin{figure}[!ht] 
\centering
\renewcommand{\datafile}{figures/arpack/data/time-decomposition.txt}
\renewcommand{\dataA}   {eigenvalues_perona1}
\renewcommand{\dataB}   {eigenvalues_lorentz}
\renewcommand{\dataC}   {eigenvalues_rudin}
\subfigure[The computational time to retrieve the eigenvectors
           associated with the lowest eigenvalues,~\eqref{eq:methodEV_eigenvector}, 
           using sequential \textsc{Arpack}.]
          {\includegraphics[scale=1]{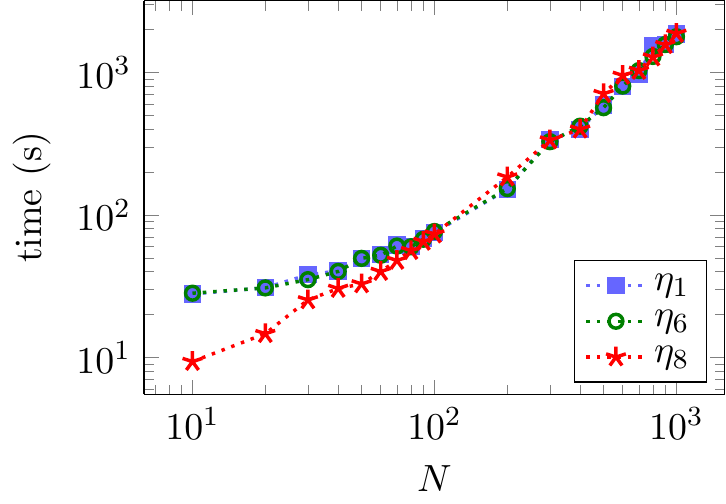}
           \label{plot:arpack:timeplot_ev}} \hspace*{-.20cm}
\renewcommand{\dataA}   {lapack_perona1}
\renewcommand{\dataB}   {lapack_lorentz}
\renewcommand{\dataC}   {lapack_rudin}
\subfigure[The computational time to obtain the coefficients
           $\alpha$ in~\eqref{eq:methodEV_decomposition}
           using least squares method and sequential 
           \textsc{Lapack}.]
          {\includegraphics[scale=1]{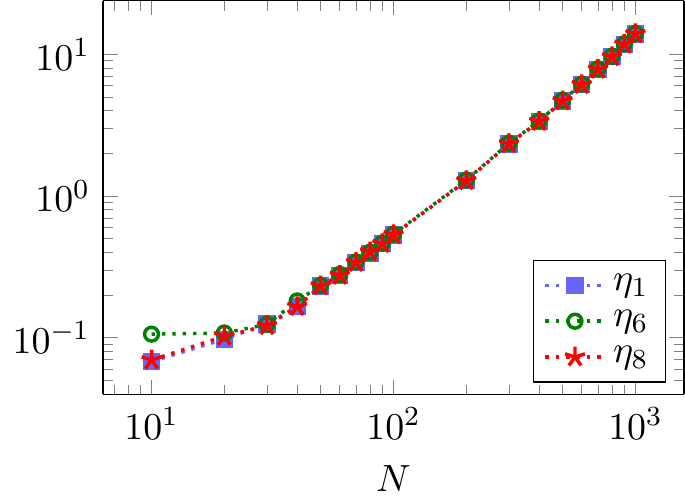}
           \label{plot:arpack:timeplot_lapack}}  
\caption{Comparison of the computational time for the 
         eigenvector decomposition for different $\eta$. 
         The computation of the eigenvectors uses the discretized 
         matrix of a diffusion operator of size 
         $\num{277221} \times \num{277221}$. It corresponds
         with the decomposition of the model 
         Figure~\ref{fig:salt_model} which is further 
         illustrated in Section~\ref{section:numeric_decomposition}.
         Other formulations of $\eta$ from Table~\ref{table:list_of_formula}
         are not shown for clarity, but follow the exact same
         pattern.}
\label{plot:arpack:timeplot}
\end{figure}

Finally, the last step is to retrieve the appropriate 
coefficients $\alpha_k$ in~\eqref{eq:methodEV_decomposition} 
for the decomposition. It basically consists in the 
resolution of a dense linear system (from least squares method).
We use \textsc{Lapack}, \cite{lapack} (contrary to 
\textsc{Mumps}, \textsc{Lapack} is adapted to dense 
linear system). Note that, because we usually
consider a few hundreds of coefficients for the 
decomposition, this operation remains relatively cheap compared 
to the eigenvectors computation. We compare the computational
time for the eigenvectors computation and model decomposition 
in Figure~\ref{plot:arpack:timeplot} for different values 
of $N$ and $\eta$. We first note that the choice of $\eta$ 
does not really modify the computational time. Then, 
we see that the two operations are mostly 
linear in $N$, and that the time to solve the least squares 
problem with \textsc{Lapack} is much smaller than the time
to compute the eigenvectors with \textsc{Arpack}, namely, 
hundred time smaller.  For the largest case: $N=1000$ for 
a squared matrix of size \num{277221}, it takes about 30min
to retrieve the eigenvectors, and $10$ \si{\second} to 
compute the $\alpha$ (in our applications we usually take 
$N < 100$).

\setlength{\plotwidth}  {5.90cm}
\setlength{\plotheight} {4cm}

\section{Illustration of model decomposition}
\label{section:numeric_decomposition}

\begin{figure}[ht!] \centering
  \graphicspath{{figures/marmousi/}}
  \renewcommand{\modelfile}{cp_marmousi_true}
  \subfigure[Marmousi velocity model.]
            {\includegraphics[scale=1]{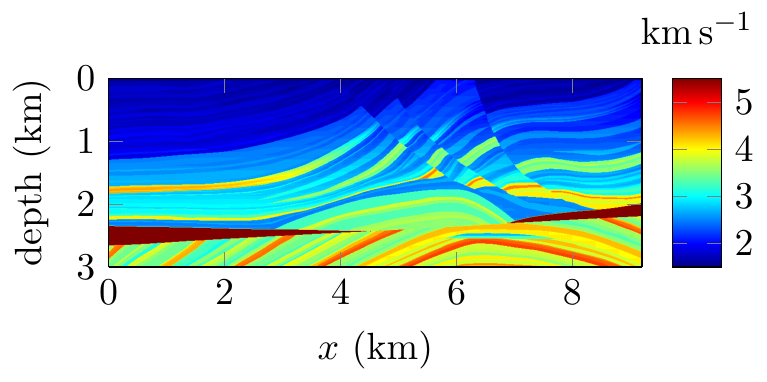} 
             \label{fig:marmousi_model}}  
  \graphicspath{{figures/salt/}}
  \renewcommand{\modelfile}{cp_salt_true}
  \subfigure[Salt velocity model.]
            {\includegraphics[scale=1]{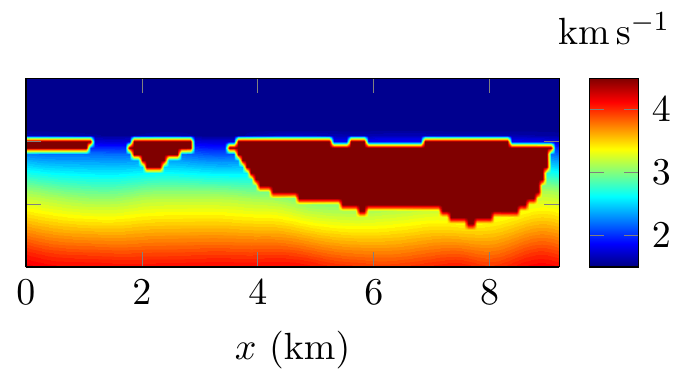}
             \label{fig:salt_model}}
  \caption{Seismic velocity models used to illustrate 
           the eigenvector decomposition. They are both
           of size $9.2\times3$ \si{\km} and represented 
           with $921\times301$ nodal coefficients.}
  \label{fig:models}
\end{figure}

First, we illustrate the eigenvector model 
decomposition with geophysical media in 
two dimensions. The original model is 
represented on a structured grid by 
$n_x \times n_z$ coefficients, and we have 
here $921\times301 = \num{277221}$ coefficients. 
We consider three media of different nature: 
\begin{itemize}
  \item the Marmousi velocity model, which consists in 
        structures and faults, see Figure~\ref{fig:marmousi_model}; 
  \item a model encompassing salt domes: objects of high 
        contrast velocity, see Figure~\ref{fig:salt_model};
  \item eventually, the SEAM Phase I velocity model
        which consists in both salt and layer structures, 
        see Figure~\ref{fig:seam_decomposition}.
\end{itemize}
All three models uses the same 
number of coefficients for their representations,
and the first two are actually of the same size
($9.2\times3$ \si{\km}).

\setlength{\modelwidth} {7.00cm}
\setlength{\modelheight}{3.50cm}

We perform the decomposition 
of the models by application of 
Algorithm~\ref{algo:decomposition_ev}, and 
steps~\eqref{eq:methodEV_m0},~\eqref{eq:methodEV_eigenvector}
and~\eqref{eq:methodEV_decomposition} (and we recall that
we use the linear PDE problem).
We study the main parameters of the decomposition: 
\begin{itemize}
 \item the choice of $\eta$, with the possibilities given in 
       Table~\ref{table:list_of_formula},
 \item the choice of the scaling parameter $\beta$ in the formulation
       of $\eta$ (Table~\ref{table:list_of_formula}),
 \item the number of eigenvectors $N$ employed for the 
       decomposition in~\eqref{eq:methodEV_decomposition}.
\end{itemize}
The accuracy of the decomposition is estimated using
the $L^2$ norm of the relative difference between 
the decomposition and the original representation such that
\begin{equation} \label{eq:error_relative}
  \err = 100 \dfrac{\Vert m - \mathfrak{m} \Vert}{\Vert m \Vert},
  \quad \quad \text{Relative Error (\%)},
\end{equation}
where $m$ is the original model (Figure~\ref{fig:models}) 
and $\mathfrak{m}$ the decomposition using the 
basis of eigenvectors from Algorithm~\ref{algo:decomposition_ev}.

\subsection{Decomposition of noise-free models}
\label{subsection:decomposition_noise_free}

We decompose the salt and Marmousi models using 
the nine possibilities for $\eta$, that are given
in Table~\ref{table:list_of_formula}. For the choice
of scaling coefficient $\beta$ (which does not 
affect $\eta_8$ and $\eta_9$), we roughly cover an 
interval from $10^{-7}$ up to $10^{6}$, namely:
$\{$$10^{-7}$, $10^{-6}$,  $10^{-5}$, $10^{-4}$,  $10^{-3}$, 
$10^{-2}$, \num{5$.10^{-2}$}, $10^{-1}$, \num{5$.10^{-1}$}, \num{1}, 
\num{5},   \num{10},   $10^{2}$,  $10^{3}$,   $10^{4}$,  $10^{5}$, $10^{6}$$\}$.
In Tables~\ref{table:norm_allmethods_marmousi} 
and~\ref{table:norm_allmethods_salt}, 
we show the best relative error (i.e. minimal value)
obtained for the Marmousi model of 
Figure~\ref{fig:marmousi_model}
and the salt model of Figure~\ref{fig:salt_model}.
We test all choices of $\eta$ and values of 
$N$ between $10$ (coarse) and $500$ (refined). 
The corresponding values of the scaling parameter 
$\beta$ which gives the best (i.e. the minimal) error
are also given in parenthesis. 

\begin{table}[ht!]
\caption{Minimal relative error obtained and associated scaling 
         coefficient: $\err$ ($\beta$), for the decomposition 
         of the Marmousi model Figure~\ref{fig:marmousi_model}.
         The definition of $\eta$ is given Table~\ref{table:list_of_formula}.}
\label{table:norm_allmethods_marmousi}
{\normalsize{ \begin{center}
\begin{tabular}{@{}llllllll}
 \renewcommand{\arraystretch}{1.00}
  Coeff.              & $N=10$                                              & $N=20$ & $N=50$  & $N=100$ & $N=250$ & $N=500$ \\ \hline
 $\eta_1$             & {\textcolor{\myblue}{\bf 6\%}}    ($10^{-6}$) &  {\textcolor{\myblue}{\bf 5\%}}    ($10^{-6}$)      &    {\textcolor{\myblue}{\bf 4\%}}    ($10^{-6}$)       &   {\textcolor{\myblue}{\bf 4\%}}    ($10^{-6}$)  &    {\textcolor{\myblue}{\bf 3\%}}    ($10^{-6}$)  &  {\textcolor{\myblue}{\bf 3\%}}    ($10^{-7}$) \\ \hline 
 $\eta_2$             &                           14\%  (5.$10^{-2}$) &                          13\%    (5.$10^{-2}$)      &                            12\%    (5.$10^{-2}$)       &                            9\%      ($10^{-2}$)  &                             7\%      ($10^{-2}$)  &                           5\%      ($10^{-2}$) \\ \hline 
 $\eta_3$             &                            8\%    ($10^{-4}$) &                           7\%      ($10^{-3}$)      &                             6\%      ($10^{-3}$)       &   5\%    ($10^{-3}$)  &                             5\%      ($10^{-3}$)  &                           4\%      ($10^{-2}$) \\ \hline 
 $\eta_4$             &                           14\%  (5.$10^{-1}$) &                          14\%      ($10^{1}$)      &                            13\%      ($10^{1}$)       &                           13\%      ($10^{-1}$)  &                            12\%    (5.$10^{-1}$)  &                          10\%    (5.$10^{1}$) \\ \hline 
 $\eta_5$             &                           13\%    ($10^{-7}$) &                          12\%      ($10^{-6}$)      &                            12\%      ($10^{-5}$)       &                           10\%      ($10^{-5}$)  &                            10\%      ($10^{-5}$)  &                           6\%      ($10^{-7}$) \\ \hline 
 $\eta_6$             &                            8\%  (5.$10^{3}$) &                           7\%      ($10^{3}$)      &                             6\%    (5.$10^{2}$)       &   5\%    ($10^{3}$)  &                             5\%      ($10^{3}$)  &                           4\%      ($10^{2}$) \\ \hline 
 $\eta_7$             &                           14\%  (5.$10^{-2}$) &                          13\%    (5.$10^{-2}$)      &                            12\%    (5.$10^{-2}$)       &                           11\%    (5.$10^{-2}$)  &                             9\%    (5.$10^{-2}$)  &                           7\%    (5.$10^{-2}$) \\ \hline 
 $\eta_8$             &                           15\%      (n/a)     &                          14\%        (n/a)          &                            13\%        (n/a)           &                           12\%        (n/a)      &                            10\%        (n/a)      &                           9\%        (n/a) \\ \hline 
 $\eta_9$             &                           14\%      (n/a)     &                          14\%        (n/a)          &                            14\%        (n/a)           &                           13\%        (n/a)      &                            12\%        (n/a)      &                          11\%        (n/a)          \\ \hline 
\end{tabular} \end{center} }}
\end{table}


\begin{table}[ht!]
\caption{Minimal relative error obtained and associated scaling 
         coefficient: $\err$ ($\beta$), for the decomposition 
         of the salt model Figure~\ref{fig:salt_model}.
         The definition of $\eta$ is given Table~\ref{table:list_of_formula}.}
\label{table:norm_allmethods_salt}
{\normalsize{ \begin{center}
\renewcommand{\arraystretch}{1.00}
\begin{tabular}{@{}llllllll}
  Coeff.              & $N=10$                                              & $N=20$ & $N=50$  & $N=100$ & $N=250$ & $N=500$ \\ \hline
 $\eta_1$             & {\textcolor{\myblue}{\bf \num{4}\%}}  ($10^{-3}$) & {\textcolor{\myblue}{\bf \num{4}\%}}  ($10^{-3}$) & {\textcolor{\myblue}{\bf \num{3}\%}}  ($10^{-2}$)  & {\textcolor{\myblue}{\bf \num{2}\%}}  ($10^{-2}$) & {\textcolor{\myblue}{\bf \num{1}\%}}  ($10^{-2}$) & {\textcolor{\myblue}{\bf \num{1}\%}}   ($10^{-2}$) \\ \hline 
 $\eta_2$             &       \num{9}\%  ($10^{-1}$)                      &                            \num{7}\%  ($10^{-1}$) &                          \num{5}\%    ($10^{-1}$)  &                          \num{5}\%    ($10^{-1}$) &                          \num{4}\%    ($10^{-1}$) &                          \num{3}\%     ($10^{-1}$) \\ \hline 
 $\eta_3$             &       \num{8}\%  (5.$10^{-2}$)                    & {\textcolor{\myblue}{\bf \num{4}\%}}(5.$10^{-2}$) & {\textcolor{\myblue}{\bf \num{3}\%}} (5.$10^{-2}$) & {\textcolor{\myblue}{\bf \num{3}\%}} (5.$10^{-2}$)& {\textcolor{\myblue}{\bf \num{2}\%}} (5.$10^{-2}$)& {\textcolor{\myblue}{\bf \num{1}\%}} (5.$10^{-2}$) \\ \hline 
 $\eta_4$             &      \num{14}\%  ($10^{-5}$)                      &                           \num{12}\%(5.$10^{1}$) &                          \num{9}\%   (5.$10^{1}$) &                          \num{8}\%   (5.$10^{1}$)&                          \num{6}\%    ($10^{1}$) &                          \num{3}\%   (5.$10^{1}$) \\ \hline 
 $\eta_5$             & {\textcolor{\myblue}{\bf \num{3}\%}}  ($10^{-6}$) & {\textcolor{\myblue}{\bf \num{3}\%}}  ($10^{-6}$) & {\textcolor{\myblue}{\bf \num{2}\%}}  ($10^{-5}$)  & {\textcolor{\myblue}{\bf \num{1}\%}}  ($10^{-5}$) & {\textcolor{\myblue}{\bf \num{1}\%}}  ($10^{-5}$) & {\textcolor{\myblue}{\bf \num{1}\%}}   ($10^{-4}$) \\ \hline 
 $\eta_6$             &       \num{6}\%  (\num{50})                       &                            \num{5}\%  ($10^{1}$) & {\textcolor{\myblue}{\bf \num{4}\%}}  ($10^{1}$)  & {\textcolor{\myblue}{\bf \num{3}\%}}  ($10^{1}$) & {\textcolor{\myblue}{\bf \num{2}\%}}  ($10^{1}$) & {\textcolor{\myblue}{\bf \num{1}\%}}   ($10^{1}$) \\ \hline 
 $\eta_7$             &       \num{9}\%  ($10^{-1}$)                      &                            \num{7}\%  ($10^{-1}$) &                          \num{5}\%    ($10^{-1}$)  &                          \num{5}\%    ($10^{-1}$) &                          \num{3}\%    ($10^{-1}$) &                          \num{3}\%     ($10^{-1}$) \\ \hline 
 $\eta_8$             &      \num{59}\%  (n/a)                            &                           \num{22}\%  (n/a)       &                         \num{16}\%    (n/a)        &                         \num{13}\%    (n/a)       &                         \num{11}\%    (n/a)       &                          \num{7}\%     (n/a)       \\ \hline 
 $\eta_9$             &      \num{20}\%  (n/a)                            &                           \num{15}\%  (n/a)       &                         \num{13}\%    (n/a)        &                         \num{10}\%    (n/a)       &                          \num{7}\%    (n/a)       &                          \num{5}\%     (n/a)       \\ \hline 
 \end{tabular} \end{center} }}
\end{table}

\begin{figure}[ht!] \centering
  \pgfmathsetmacro{\ymin}{1.5} 
  \renewcommand{\datafolder}{figures/decomposition_error_allmethods/data/}
  \renewcommand{\dataA}      {Var50}
  \renewcommand{\dataformula}{perona1}
  \pgfmathsetmacro{\ymax}{20}
  \subfigure[Decomposition using $\eta_1$.]{
              \includegraphics[scale=1]{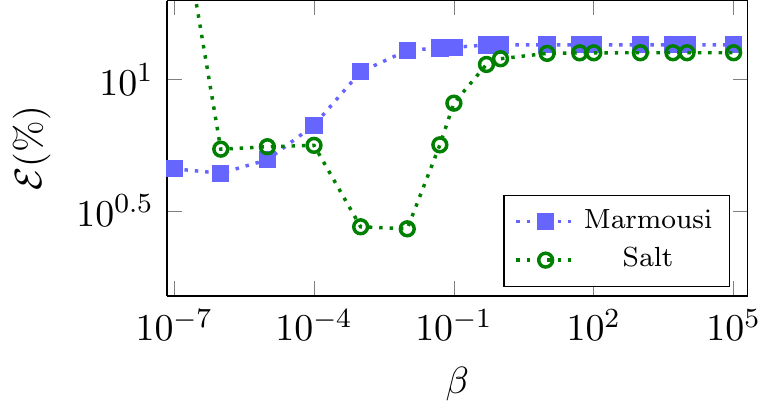}} 
              \hspace*{-.75cm}
  \pgfmathsetmacro{\ymax}{100}
  \renewcommand{\dataformula}{geman} 
  \subfigure[Decomposition using $\eta_3$.]{
              \includegraphics[scale=1]{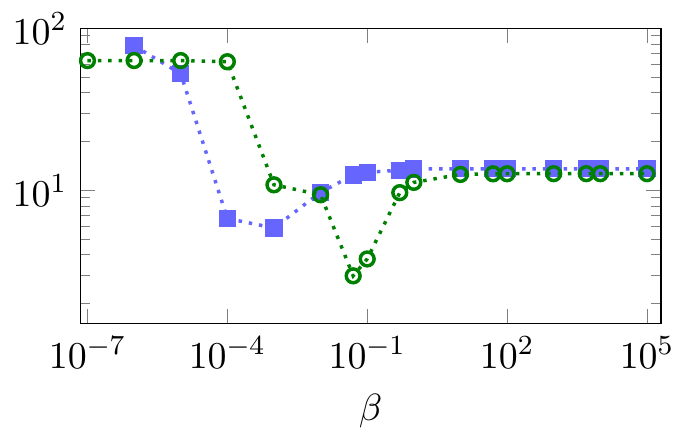}}
  \pgfmathsetmacro{\ymax}{18}
  \renewcommand{\dataformula}{aubert} 
  \subfigure[Decomposition using $\eta_5$.]{
              \includegraphics[scale=1]{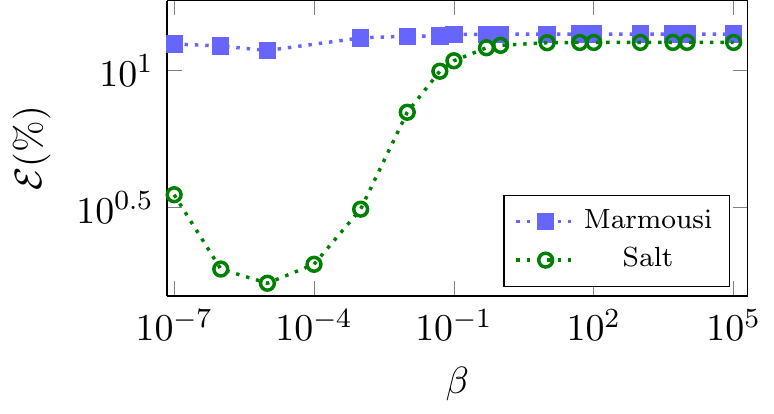}} 
              \hspace*{-.75cm}
  \pgfmathsetmacro{\ymax}{50}
  \renewcommand{\dataformula}{lorentz} 
  \subfigure[Decomposition using $\eta_6$.]{
              \includegraphics[scale=1]{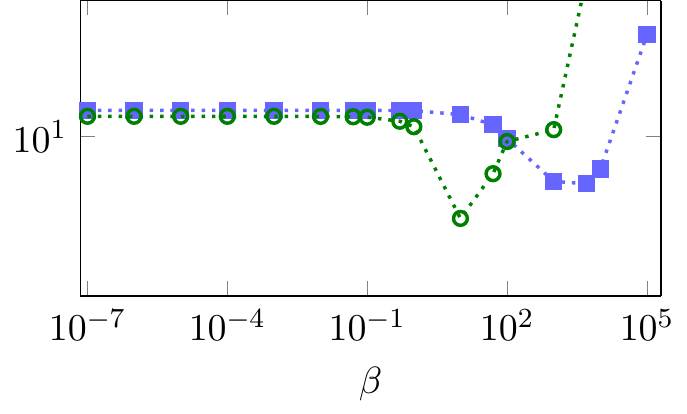}}
  \caption{Decomposition of the Marmousi and salt 
           velocity models of  
           Figures~\ref{fig:marmousi_model}
           and~\ref{fig:salt_model} using $N=50$ 
           eigenvectors and following Algorithm~\ref{algo:decomposition_ev}. 
           The relative error is computed
           from~\eqref{eq:error_relative} for four 
           selected formulation of $\eta$ 
           (see Table~\ref{table:list_of_formula})
           and different scaling parameter $\beta$.}
  \label{fig:decomposition_error_per_formula}
\end{figure}

As expected, we observe that the more eigenvectors 
are chosen (higher $N$), the better will be the 
decomposition. When using $500$ eigenvectors, 
which represents about 2\% of the original number
of coefficients ($921\times301$ in 
Figure~\ref{fig:models}), the error is of a few percent
only. This can be explained by the redundancy of information
provided by the original fine grid where the model
is represented (e.g. the upper part of Figure~\ref{fig:salt_model}
and the three salt bodies are basically constant).
Comparing the methods and models, we see that
\begin{itemize}
 \item the Marmousi model (Table~\ref{table:norm_allmethods_marmousi})
       is harder to retrieve than the salt model 
       (Table~\ref{table:norm_allmethods_salt}) as it gives
       higher errors. In particular for low $N$, the salt 
       model can be acutely decomposed (possibly $3$\% error with $N=10$).
 \item For both models, it appears that four methods stand out: 
       $\eta_1$ (Perona--Malik), $\eta_3$ (\cite{Geman1992}),
       $\eta_5$ (\cite{Charbonnier1994})
       and $\eta_6$ (Lorentzian), with a slight advantage 
       towards $\eta_1$.
 \item The scaling coefficient that minimizes the error is 
       consistent with respect to $N$, with similar amplitude.
       However, changing the model may require the modification
       of $\beta$: between the salt and Marmousi decomposition, 
       the optimal $\beta$ is quite different for $\eta_1$,
       and also for $\eta_6$.
\end{itemize}
To investigate further the last point, we show 
the evolution of relative error $\err$ with 
respect to the scaling coefficient $\beta$ for 
the decomposition of the Marmousi and salt 
models in Figure~\ref{fig:decomposition_error_per_formula},
where we compare four selected formulations 
for $\eta$. 
We observe some flexibility in the choice of 
$\beta$ that gives an accurate decomposition 
for $\eta_1$ and $\eta_5$. 
On the other hand, $\eta_3$ and $\eta_6$ show
sharp functions, which means that the selection 
of $\beta$ has more influence in these cases
(and must be carefully taken). 
In addition, the range of efficient $\beta$ 
changes depending on the model decomposed, 
except for $\eta_5$. It demonstrates that the 
choice of $\beta$ for optimality is not trivial 
in general, and is model-dependent.

\setlength{\plotwidth}  {5.9cm}
\setlength{\plotheight} {3.0cm}

\begin{figure}[ht!] \centering
  \includegraphics[scale=1]{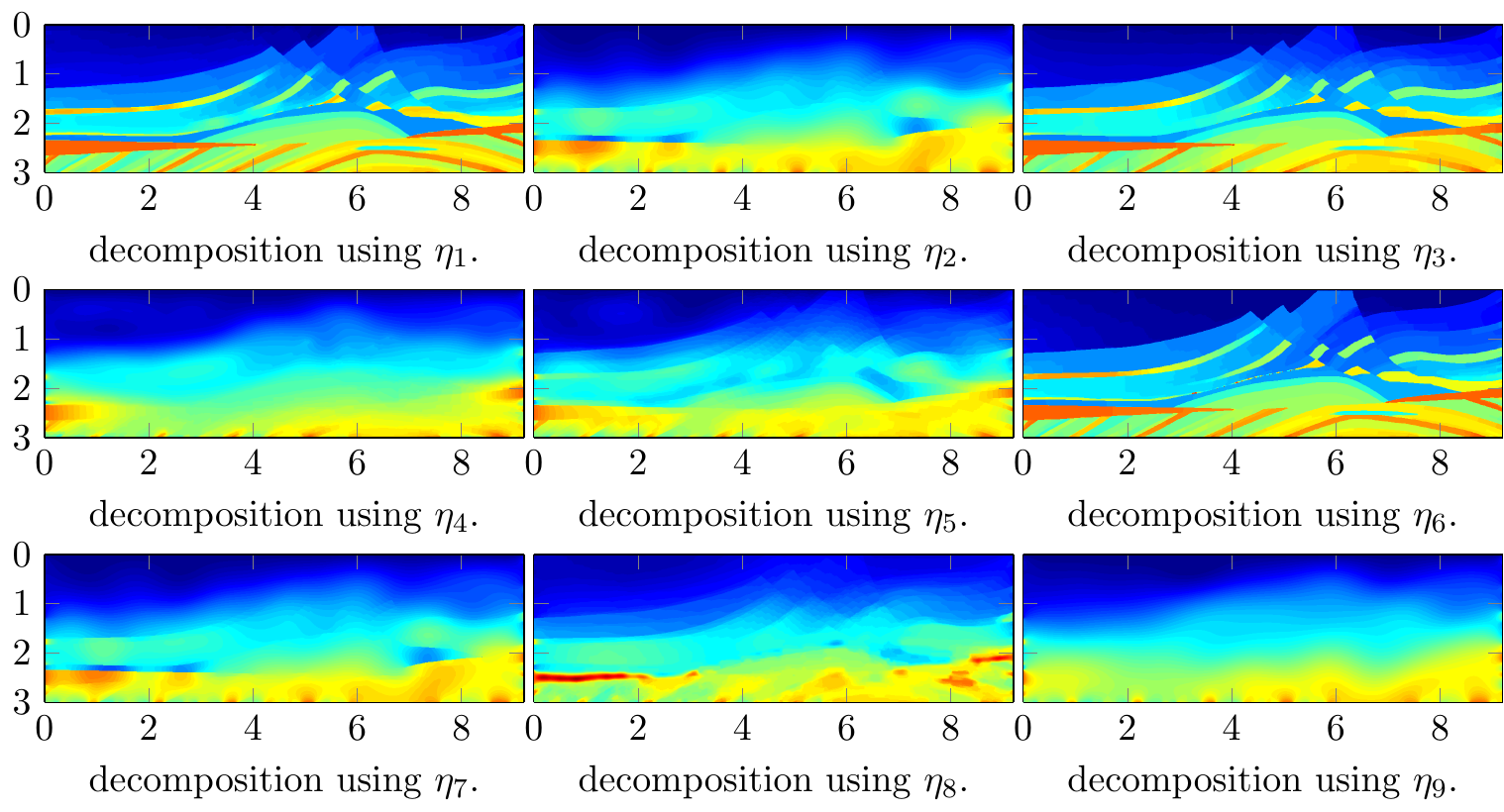}
  \vspace*{-.5\baselineskip}
  \caption{Decomposition of the noise-free Marmousi 
           velocity model of Figure~\ref{fig:marmousi_model} 
           using $N=50$, the formulation for $\eta$ are 
           from Table~\ref{table:list_of_formula} and 
           the respective $\mu$ values extracted 
           from Table~\ref{table:norm_allmethods_marmousi}.
           The color scale varies
           between \num{1500} and 
           \num{5500}\si{\meter\persecond}.}
  \label{fig:models_marmousi_decomposition}
\end{figure}

\begin{figure}[ht!] \centering
  \includegraphics[scale=1]{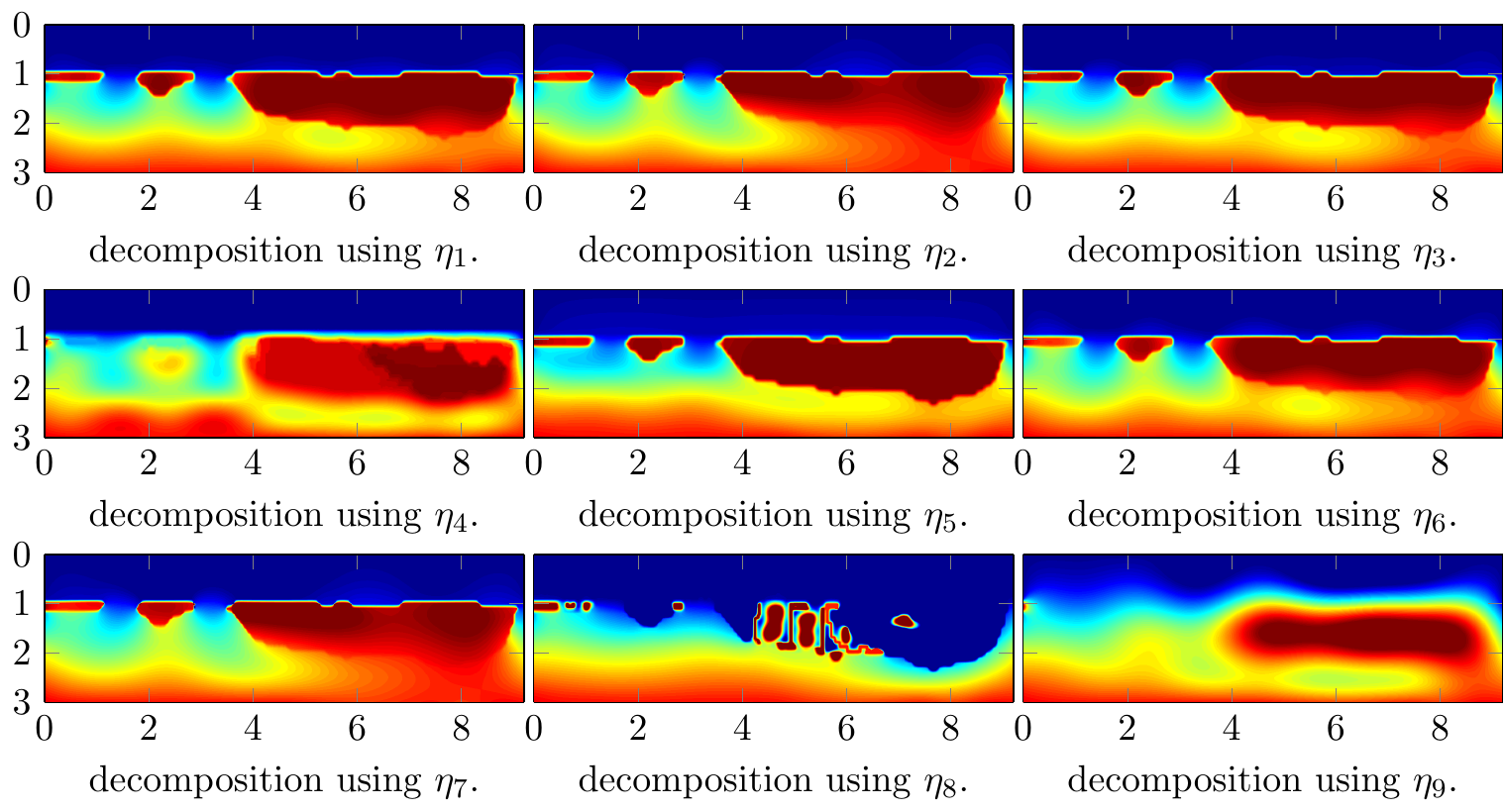}
  \vspace*{-.5\baselineskip}
  \caption{Decomposition of the noise-free salt
           velocity model of Figure~\ref{fig:salt_model} 
           using $N=20$, the formulation for $\eta$ are 
           from Table~\ref{table:list_of_formula} and 
           the respective $\mu$ values extracted 
           from Table~\ref{table:norm_allmethods_salt}.
           The color scale varies
           between \num{1500} and 
           \num{4500}\si{\meter\persecond}.}
  \label{fig:models_salt_decomposition}
\end{figure}

We then picture the resulting images 
obtained after the decomposition of both models,
see Figures~\ref{fig:models_marmousi_decomposition}
and~\ref{fig:models_salt_decomposition}.
The pictures illustrate correctly the observations
of the tables and the differences between the formulation. 
The salt model is usually well recovered with all formulations, 
while the Marmousi model is more hardly discovered, except 
with $\eta_1$, $\eta_3$ and $\eta_6$.
Those three formulations are the only ones able to capture
the structures.


\setlength{\modelwidth} {6.50cm}
\setlength{\modelheight}{3.10cm}
\renewcommand{\modelfile}{marmousi-decomposition_50ev}
\graphicspath{{figures/marmousi_decomposition/}}

\renewcommand{\modelfile}{salt-decomposition_20ev}
\graphicspath{{figures/salt_decomposition/}}

\subsection{Decomposition of noisy models (denoising)}

\begin{figure}[ht!] \centering
  \graphicspath{{figures/marmousi/}}
  \renewcommand{\modelfile}{cp_marmousi_true_noise20}
  \subfigure[Marmousi velocity model with noise.]
            {\includegraphics[scale=1]{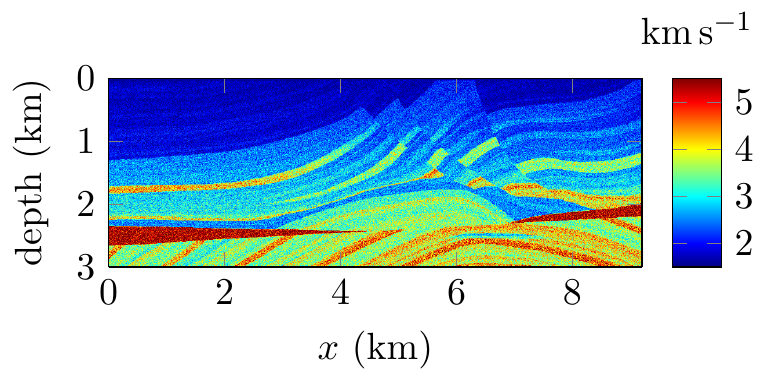}
             \label{fig:marmousi_model_noise20}}  
  \graphicspath{{figures/salt/}}
  \renewcommand{\modelfile}{cp_salt_true_noise20}
  \subfigure[Salt velocity model with noise.]
            {\includegraphics[scale=1]{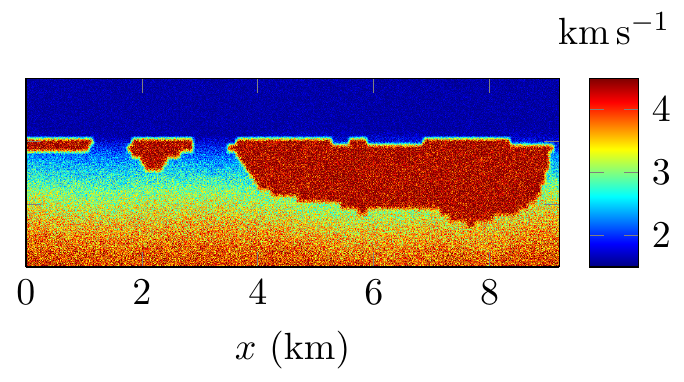}
             \label{fig:salt_model_noise20}}
  \caption{Seismic velocity models of Figure~\ref{fig:models}
           with added noise. For every coefficient, the noisy 
           one is obtained by randomly taking a value between 
           $\pm 20\%$ of the noiseless coefficient.}
  \label{fig:models_noise}
\end{figure}

We incorporate noise in the representation for 
getting closer to the reality of applications 
where few information on the model is available.
Hence, we reproduce the model decomposition, 
this time working with noisy pictures. 
For every nodal velocity (of Figure~\ref{fig:models}), 
we recast the values using an uniform distribution 
that covers $\pm 20\%$ of the noiseless value. 
The resulting media are illustrated in 
Figure~\ref{fig:models_noise}.

\setlength{\modelwidth} {7.00cm}
\setlength{\modelheight}{3.50cm}

We apply the model decomposition using 
the different formulations of $\eta$ 
and choice of scaling coefficient 
$\beta$, following the procedure employed 
for the noiseless model. In 
Tables~\ref{table:norm_allmethods_marmousi_noise}
and~\ref{table:norm_allmethods_salt_noise},
we show the evolution of best relative error
with $N$, for the noisy Marmousi and salt models 
respectively. Here, the relative error is computed
from the difference between the noiseless model
and the decomposition of the noisy one.
The objective of the regularization is 
to preserve the structures while smoothing
out the noise effect.

\begin{table}[ht!]
\caption{Minimal relative error obtained and associated scaling 
         coefficient: $\err$ ($\beta$). The error is computed 
         with respect to the noiseless model 
         Figure~\ref{fig:marmousi_model} but the decomposition 
         uses the noisy model of Figure~\ref{fig:marmousi_model_noise20}.
         The definition of $\eta$ is given Table~\ref{table:list_of_formula}.}
\label{table:norm_allmethods_marmousi_noise}
{\normalsize{ \begin{center}
 \renewcommand{\arraystretch}{1.00}
 \begin{tabular}{@{}llllllll}
  Coeff.              &  $N=10$ & $N=20$ & $N=50$  & $N=100$ & $N=250$  \\ \hline
 $\eta_1$             &  16\%      ($10^{-4}$)  &  15\%      ($10^{-4}$)  &   14\%      ($10^{-7}$)  & 13\%      ($10^{-7}$)  & 11\%      ($10^{-7}$)   \\ \hline 
 $\eta_2$             &  16\%    (5.$10^{-2}$)  &  16\%    (5.$10^{-2}$)  &   15\%    (5.$10^{-2}$)  & 14\%    (5.$10^{-2}$)  & 13\%    (5.$10^{-2}$)   \\ \hline 
 $\eta_3$             &  16\%      ($10^{-2}$)  &  14\%      ($10^{-3}$)  &   13\%      ($10^{-4}$)  & 12\%      ($10^{-4}$)  & 10\%      ($10^{-4}$)   \\ \hline 
 $\eta_4$             &  16\%      ($10^{1}$)  &  16\%      ($10^{1}$)  &   16\%    (5.$10^{1}$)  & 15\%    (5.$10^{1}$)  & 14\%      ($10^{2}$)   \\ \hline 
 $\eta_5$             &  16\%      ($10^{-6}$)  &  16\%      ($10^{-7}$)  &   15\%      ($10^{-7}$)  & 15\%      ($10^{-7}$)  & 14\%      ($10^{-6}$)   \\ \hline 
 $\eta_6$             &  16\%      ($10^{2}$)  &  14\%      ($10^{3}$)  &   13\%      ($10^{4}$)  & 12\%      ($10^{4}$)  & 10\%      ($10^{4}$)   \\ \hline 
 $\eta_7$             &  16\%    (5.$10^{-2}$)  &  16\%    (5.$10^{-2}$)  &   15\%    (5.$10^{-2}$)  & 14\%    (5.$10^{-2}$)  & 13\%    (5.$10^{-2}$)   \\ \hline 
 $\eta_8$             &  17\%        (n/a)      &  16\%        (n/a)      &   16\%        (n/a)      & 15\%        (n/a)      & 14\%        (n/a)       \\ \hline 
 $\eta_9$             &  17\%        (n/a)      &  17\%        (n/a)      &   16\%        (n/a)      & 16\%        (n/a)      & 14\%        (n/a)       \\ \hline 
 \end{tabular} \end{center} }}
\end{table}



\begin{table}[ht!]
\caption{Minimal relative error obtained and associated scaling 
         coefficient: $\err$ ($\beta$). The error is computed 
         with respect to the noiseless model 
         Figure~\ref{fig:salt_model} but the decomposition 
         uses the noisy model of Figure~\ref{fig:salt_model_noise20}.
         The definition of $\eta$ is given Table~\ref{table:list_of_formula}.}
\label{table:norm_allmethods_salt_noise}
{\normalsize{ \begin{center}
 \renewcommand{\arraystretch}{1.00}
 \begin{tabular}{@{}llllllll}
  Coeff.              &  $N=10$ & $N=20$ & $N=50$  & $N=100$ & $N=250$  \\ \hline
 $\eta_1$             &                           14\%      ($10^{-6}$) &                          14\%      ($10^{-6}$) &                          11\%      ($10^{-5}$) &                           9\%      ($10^{-5}$)  &  5\%      ($10^{-3}$)   \\ \hline 
 $\eta_2$             &                           17\%    (5.$10^{-2}$) &                          15\%    (5.$10^{3}$) &                          11\%    (5.$10^{-2}$) &                           8\%    (5.$10^{-2}$)  &  6\%    (5.$10^{-2}$)   \\ \hline 
 $\eta_3$             &  {\textcolor{\myblue}{\bf 10\%}}    ($10^{-4}$) & {\textcolor{\myblue}{\bf 10\%}}    ($10^{-4}$) & {\textcolor{\myblue}{\bf  8\%}}    ($10^{-4}$) & {\textcolor{\myblue}{\bf  6\%}}    ($10^{-4}$)  &  5\%      ($10^{-2}$)   \\ \hline 
 $\eta_4$             &                           18\%      ($10^{-3}$) &                          15\%      ($10^{0}$) &                          12\%      ($10^{-1}$) &                          10\%      ($10^{-3}$)  &  6\%      ($10^{-3}$)   \\ \hline 
 $\eta_5$             &                           17\%      ($10^{-7}$) &                          15\%    (5.$10^{3}$) &                          12\%      ($10^{-4}$) &                           9\%      ($10^{-7}$)  &  6\%      ($10^{-5}$)   \\ \hline 
 $\eta_6$             &  {\textcolor{\myblue}{\bf 10\%}}    ($10^{4}$) & {\textcolor{\myblue}{\bf  9\%}}    ($10^{4}$) & {\textcolor{\myblue}{\bf  8\%}}    ($10^{4}$) & {\textcolor{\myblue}{\bf  6\%}}  (5.$10^{3}$)  &  5\%      ($10^{2}$)   \\ \hline 
 $\eta_7$             &                           17\%    (5.$10^{-2}$) &                          15\%      ($10^{5}$) &                          11\%      ($10^{-2}$) &                           8\%    (5.$10^{-2}$)  &  6\%    (5.$10^{-2}$)   \\ \hline 
 $\eta_8$             &                           18\%        (n/a)     &                          16\%        (n/a)     &                          12\%        (n/a)     &                           9\%        (n/a)      &  6\%        (n/a)       \\ \hline 
 $\eta_9$             &                           21\%        (n/a)     &                          15\%        (n/a)     &                          13\%        (n/a)     &                          11\%        (n/a)      &  8\%        (n/a)       \\ \hline 
 \end{tabular} \end{center} }}
\end{table}


The decomposition of noisy pictures
requires more eignevectors for an 
accurate representation. Then, the salt model, 
with high contrast objects, still behaves 
better than the many structures of the 
Marmousi model.  For the decomposition of 
the noisy Marmousi model, none of the 
formulations really stands out and the 
error never reaches below 10\% using at 
most $N=250$.

\begin{figure}[ht!] \centering
  \includegraphics[scale=1]{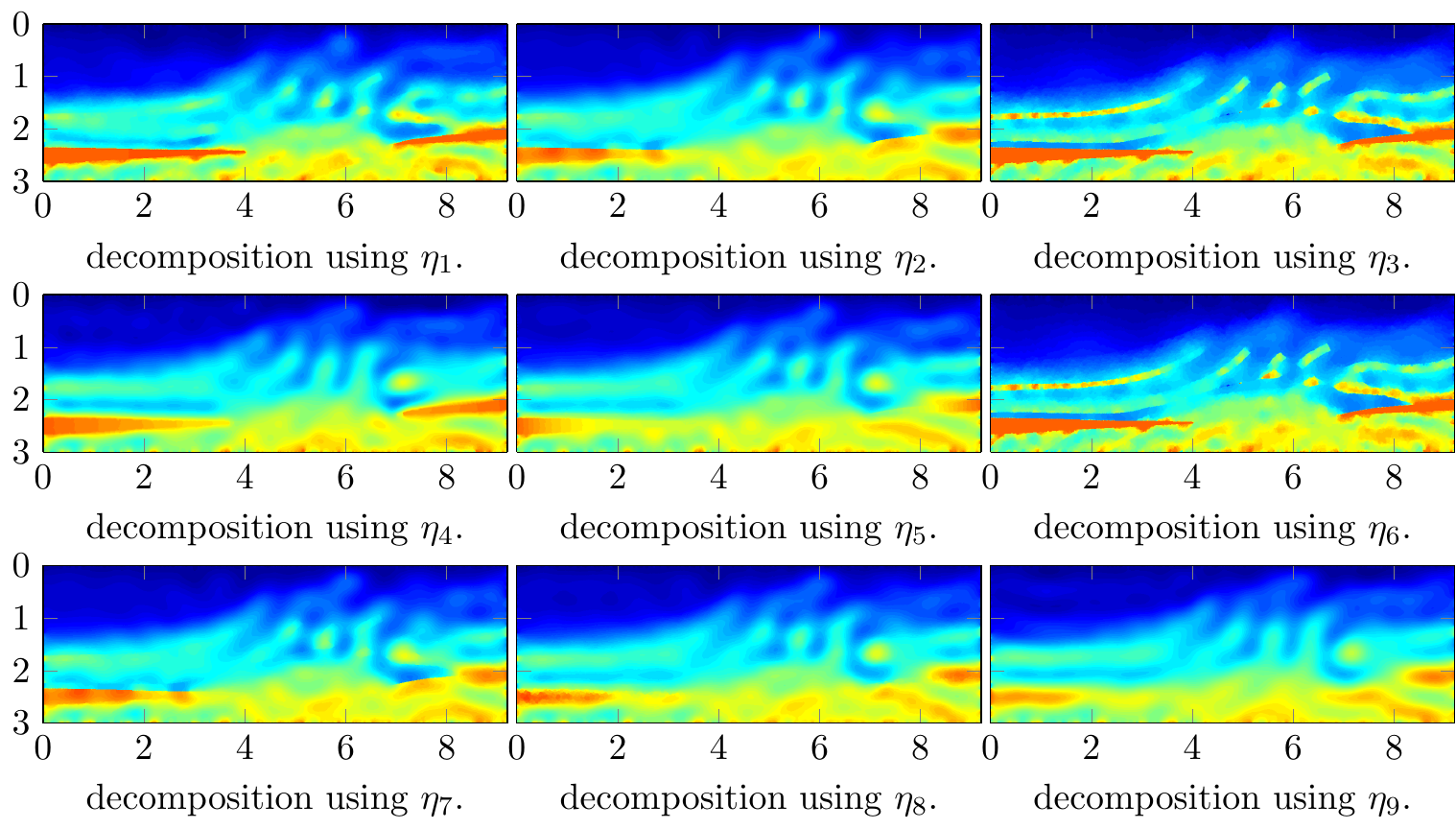}
  \caption{Decomposition of the Marmousi velocity model 
           of Figure~\ref{fig:marmousi_model_noise20} 
           using $N=200$
           and the formulations of $\eta$ from 
           Table~\ref{table:list_of_formula}. The selected
           value of $\beta$ for every method corresponds 
           to the value given in 
           Table~\ref{table:norm_allmethods_marmousi_noise}. 
           The color scale follows the one of 
           Figure~\ref{fig:marmousi_model_noise20}
           with values between \num{1500} and 
           \num{5500}\si{\meter\persecond}.}
  \label{fig:models_marmousi_decomposition_noise}
\end{figure}

\begin{figure}[ht!] \centering
  \includegraphics[scale=1]{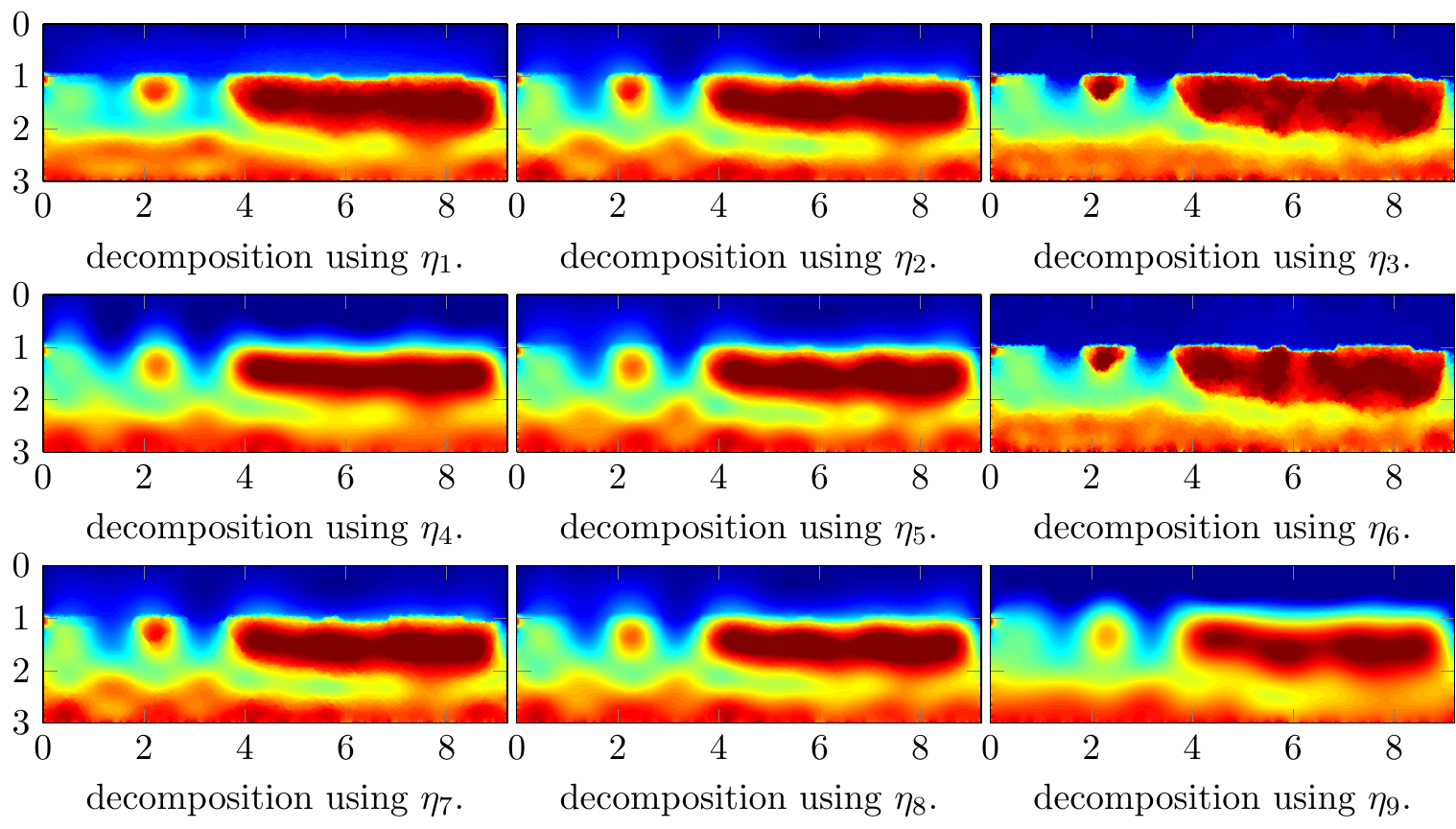}
  \vspace*{-.5\baselineskip}
  \caption{Decomposition of the salt velocity model 
           of Figure~\ref{fig:salt_model_noise20} using $N=50$
           and the formulations of $\eta$ from 
           Table~\ref{table:list_of_formula}. The selected
           value of $\beta$ for every method corresponds 
           to the value given in 
           Table~\ref{table:norm_allmethods_salt_noise}. 
           The color scale follows the one of 
           Figure~\ref{fig:salt_model_noise20}
           with values between \num{1500} and 
           \num{4500}\si{\meter\persecond}.}
  \label{fig:models_salt_decomposition_noise}
\end{figure}

In Figures~\ref{fig:models_marmousi_decomposition_noise}
and~\ref{fig:models_salt_decomposition_noise}, we 
picture the resulting decomposition 
for the two media. For the decomposition of the Marmousi
model, we use $N=200$; and $N=50$ for 
the salt model. It corresponds to 
higher values compared to the pictures shown for the 
noiseless models (Figures~\ref{fig:models_marmousi_decomposition}
and~\ref{fig:models_salt_decomposition}).

The decomposition of the salt model remains acceptable, 
and we easily distinguish the main contrasting object.
The smaller objects also appear, in a smooth representation.
The formulations using $\eta_1$, $\eta_3$ and $\eta_6$
provide sharper boundary for the contrasting 
objects, in particular for the upper interface.
Regarding the decomposition of the noisy Marmousi 
model, it illustrates the limitation of the 
method, where none of the formulations is really able 
to reproduce the structures, and most edges are lost.
In particular, the central part of the model is mostly
missing and the amplitude of the values has been 
reduced. It seems that $\eta_1$, $\eta_3$ and $\eta_6$
are slightly more robust and gives (relatively speaking) 
the best results. To conclude, these three formulations
appear less sensitive (for both media) to noise than 
the other ones.

\setlength{\modelwidth} {6.50cm}
\setlength{\modelheight}{3.25cm}
\renewcommand{\modelfile}{marmousi-decomposition_200ev}
\graphicspath{{figures/marmousi_decomposition_noise20/}}
\setlength{\modelwidth} {6.50cm}
\setlength{\modelheight}{3.25cm}
\renewcommand{\modelfile}{salt-decomposition_50ev}
\graphicspath{{figures/salt_decomposition_noise20/}}

In the context of image decomposition, we have 
can draw the following conclusions for the the 
decomposition using the basis of eigenvectors.
\begin{itemize}
  \item The method is efficient to represent media 
        with contrasting shapes (e.g., salt domes), 
        even when noise is contained in the images. 
        In this case, the choice of $\eta$ does 
        not really affect the representation of the 
        objects, and all methods behave quite well, see
        Figure~\ref{fig:models_salt_decomposition_noise}.
  \item The performance of the decomposition 
        strongly depends on the media, and diminishes 
        with thin structures as in the Marmousi model.
        In this case, an appropriate choice 
        of formulation ($\eta_1$, $\eta_3$, $\eta_6$
        from Table~\ref{table:norm_allmethods_marmousi})
        can provide the accurate representation 
        for noise-free picture but when incorporating noise, the 
        performance deteriorates and the edge contrasts are 
        lost.
\end{itemize}

\subsection{Sub-surface salt with layers: SEAM Phase I model}

\begin{figure}[ht!] \centering
  \graphicspath{{figures/seam/main/}}
  \renewcommand{\modelfile}{cp_true}
  \subfigure[Noise-free velocity model.]
            {\includegraphics[scale=1]{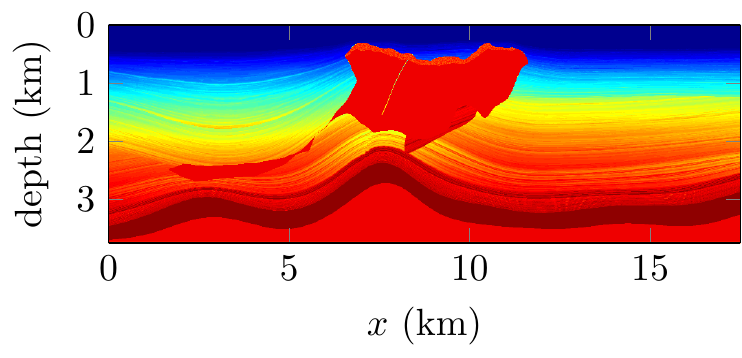}
             \label{fig:seam_true}}  
  \renewcommand{\modelfile}{cp_true_noise20}
  \subfigure[Velocity model with noise.]
            {\includegraphics[scale=1]{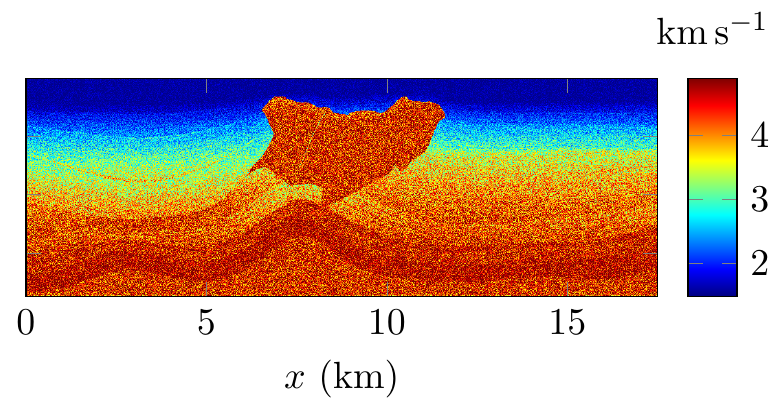}
             \label{fig:seam_true_noise20}}
  \caption{SEAM velocity model of size 
           $17.5 \times 3.75$ \si{\km}, represented
           on a Cartesian grid of size $921 \times 301$.}
  \label{fig:seam_decomposition}
\end{figure}

We have seen that the decomposition behaves well when a 
contrasting object with sharp contrast belongs to the medium,
while structures/layers are hardly represented. 
We pursue our investigation with a common geophysical 
configuration where both salt-domes and layers exist in the 
subsurface. 
We use a velocity model extracted from the SEAM
(SEG Advanced Modeling Program ) Phase I 
benchmark\footnote{see \url{https://wiki.seg.org/wiki/Open_data}.}
a consider a medium of size $17.5 \times 3.75$ \si{\km}.
Per consistency with the previous experiment, it 
is represented using a grid of $921 \times 301$ points,
see Figure~\ref{fig:seam_decomposition}.

\setlength{\modelwidth} {8.00cm}
\setlength{\modelheight}{3.80cm}

Similarly to our previous experiment, 
we investigate the noiseless model and a noisy
one, which incorporates $\pm 20$ \% error. 
The relative error and corresponding scaling 
coefficients for both models are given 
in Tables~\ref{table:norm_allmethods_seam}
and~\ref{table:norm_allmethods_seam_noise20}.
The relative error is of a few percent for 
high $N$, and we observe important differences
between the formulation. Here, $\eta_8$,
$\eta_9$  and $\eta_4$ give the worst results.

\begin{figure}[ht!] \centering
  \graphicspath{{figures/seam/decomposition/}}
  \renewcommand{\modelfile}{cp_decomposition-noiseless_20ev}
  \subfigure[Decomposition of noiseless model using $N=20$.]
            {\includegraphics[scale=1]{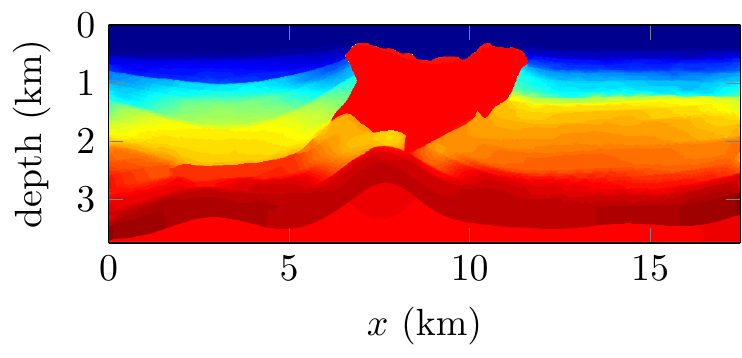}
             \label{fig:seam_decomposition_results_noisefree}}  
  \renewcommand{\modelfile}{cp_decomposition-noise20_100ev}
  \subfigure[Decomposition of noisy model using $N=100$.]
            {\includegraphics[scale=1]{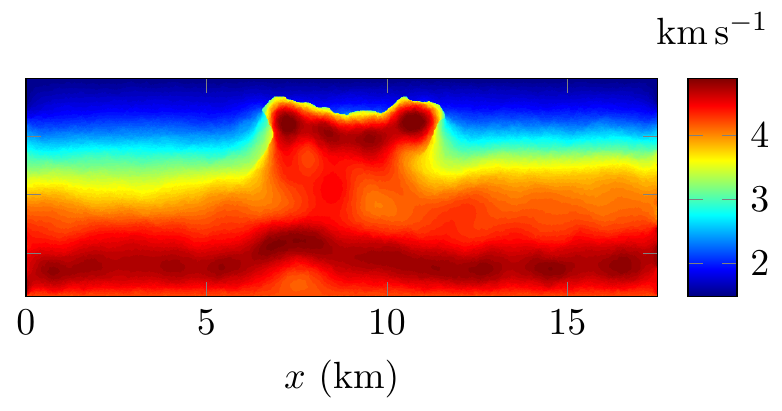}}
  \caption{Decomposition of the SEAM Phase I
           velocity model of Figure~\ref{fig:seam_decomposition} 
           using $\eta_1$. The relative errors for all
           formulations can be found in 
           Tables~\ref{table:norm_allmethods_seam}
           and~\ref{table:norm_allmethods_seam_noise20}.}
  \label{fig:seam_decomposition_results}
\end{figure}

We further illustrate the decomposition
in Figure~\ref{fig:seam_decomposition_results},
using $\eta_1$ only for the sake of clarity.

\begin{table}[ht!]
\caption{Minimal relative error obtained and associated scaling 
         coefficient: $\err$ ($\beta$), for the decomposition 
         of the SEAM model Figure~\ref{fig:seam_true}.
         The definition of $\eta$ is given Table~\ref{table:list_of_formula}.}
\label{table:norm_allmethods_seam}
{\normalsize{ \begin{center}
\begin{tabular}{@{}lllllll}
 \renewcommand{\arraystretch}{1.00}
  Coeff.              & $N=10$ & $N=20$ & $N=50$  & $N=100$ & $N=200$ \\ \hline
 \tablerrseam{$\eta_1$}{4\% ($10^{-6}$)}
                       {3\% ($10^{-6}$)}
                       {3\% ($10^{-6}$)}
                       {2\% ($10^{-6}$)}
                       {2\% ($10^{-7}$)}   \\ \hline 
 \tablerrseam{$\eta_2$}{8\% (5.$10^{-2}$)}
                       {6\% (5.$10^{-2}$)}
                       {5\% (  $10^{-1}$)}
                       {3\% (  $10^{-2}$)}
                       {3\% (  $10^{-2}$)}   \\ \hline 
 \tablerrseam{$\eta_3$}{6\% ($10^{-3}$)}
                       {5\% ($10^{-2}$)}
                       {4\% ($10^{-3}$)}
                       {3\% ($10^{-2}$)}
                       {2\% ($10^{-3}$)}   \\ \hline 
\tablerrseam{$\eta_4$} {10\% (  5)}
                       { 9\% ( 10)}
                       { 7\% ($10^{2}$)}
                       { 6\% ($10^{2}$)}
                       { 5\% (  5)}        \\ \hline 
\tablerrseam{$\eta_5$} {7\% ($10^{-7}$)} 
                       {6\% ($10^{-6}$)} 
                       {4\% ($10^{-6}$)} 
                       {3\% ($10^{-5}$)} 
                       {3\% ($10^{-5}$)}    \\ \hline 
\tablerrseam{$\eta_6$} {6\% ($10^{3}$)}
                       {5\% ($10^{2}$)}
                       {4\% ($10^{3}$)}
                       {3\% ($10^{2}$)}
                       {2\% ($10^{3}$)}   \\ \hline 
\tablerrseam{$\eta_7$} {8\% ($5.10^{-2}$)}
                       {6\% ($5.10^{-2}$)}
                       {5\% ($  10^{-1}$)}
                       {4\% ($5.10^{-2}$)}
                       {3\% ($5.10^{-2}$)}   \\ \hline 
\tablerrseam{$\eta_8$} {18\% (n/a)} 
                       {17\% (n/a)} 
                       {14\% (n/a)} 
                       {12\% (n/a)} 
                       { 9\% (n/a)}    \\ \hline 
\tablerrseam{$\eta_9$} {12\% (n/a)} 
                       {10\% (n/a)} 
                       { 7\% (n/a)} 
                       { 7\% (n/a)} 
                       { 6\% (n/a)}    \\ \hline 
\end{tabular} \end{center} }}
\end{table}

\begin{table}[ht!]
\caption{Minimal relative error obtained and associated scaling 
         coefficient: $\err$ ($\beta$), for the decomposition 
         of the SEAM model Figure~\ref{fig:seam_true_noise20}.
         The definition of $\eta$ is given Table~\ref{table:list_of_formula}.}
\label{table:norm_allmethods_seam_noise20}
{\normalsize{ \begin{center}[ht!]
\begin{tabular}{@{}lllllll}
 \renewcommand{\arraystretch}{1.00}
  Coeff.              & $N=10$ & $N=20$ & $N=50$  & $N=100$ & $N=200$ \\ \hline
 \tablerrseam{$\eta_1$}{11\% ($10^{-3}$)}
                       { 8\% ($10^{-3}$)}
                       { 6\% ($10^{-6}$)}
                       { 5\% ($10^{-7}$)}
                       { 4\% ($10^{-5}$)}  \\ \hline
 \tablerrseam{$\eta_2$}{8\% ($5.10^{-2}$)}
                       {6\% ($5.10^{-2}$)}
                       {6\% ($5.10^{-2}$)}
                       {4\% ($5.10^{-2}$)}
                       {3\% ($5.10^{-2}$)}  \\ \hline
 \tablerrseam{$\eta_3$}{6\% ($10^{-2}$)}
                       {5\% ($10^{-2}$)}
                       {4\% ($10^{-3}$)}
                       {4\% ($10^{-3}$)}
                       {4\% ($10^{-2}$)}  \\ \hline
\tablerrseam{$\eta_4$} {11\% ($10^{-2}$)}
                       {10\% ($0.5$)    }
                       { 7\% ($10^{-4}$)}
                       { 6\% ($10^{-3}$)}
                       { 5\% ($10^{-3}$)}  \\ \hline
\tablerrseam{$\eta_5$} {11\% ($10^{-4}$)}
                       {10\% ($1$)      }
                       { 7\% ($10^{-2}$)}
                       { 6\% ($10^{-4}$)}
                       { 5\% ($10^{-5}$)} \\ \hline
\tablerrseam{$\eta_6$} {8\% ($10^{2}$)}
                       {6\% ($10^{2}$)}
                       {5\% ($10^{2}$)}
                       {4\% ($10^{2}$)}
                       {3\% ($10^{2}$)}   \\ \hline
\tablerrseam{$\eta_7$} {8\% ($5.10^{-2}$)}
                       {6\% ($5.10^{-2}$)}
                       {6\% ($5.10^{-2}$)}
                       {4\% ($5.10^{-2}$)}
                       {3\% ($5.10^{-2}$)} \\ \hline
\tablerrseam{$\eta_8$} {11\% (n/a)}
                       {11\% (n/a)}
                       { 7\% (n/a)}
                       { 6\% (n/a)}
                       { 5\% (n/a)} \\ \hline
\tablerrseam{$\eta_9$} {12\% (n/a)}
                       {10\% (n/a)}
                       { 7\% (n/a)}
                       { 7\% (n/a)}
                       { 6\% (n/a)} \\ \hline
\end{tabular} \end{center} }}
\end{table}

\setlength{\modelwidth} {8.00cm}
\setlength{\modelheight}{3.80cm}

This model, which encompasses salt and layers, is well
recovered with a decomposition using $N=20$ when there is 
no noise. In case of noise, it needs higher $N$ and 
the contrasting shapes are smoothed out. 
Nonetheless, compared to the Marmousi model, the 
decomposition is able to capture the main features.

\section{Experiments of reconstruction with FWI}
\label{section:numeric_reconstruction}

In this section, we perform seismic imaging 
following the FWI Algorithm~\ref{algo:fwi}
for the identification of the subsurface 
physical parameters.
We focus on media encompassing salt domes, 
as we have shown that the decomposition is 
more appropriate, and 
it also gives a more challenging configuration
in seismic applications.
We follow a seismic context, where the forward
problem is given by~\eqref{eq:helmholtz_equation}
and the data are restricted to be acquired near the
surface (see Figure~\ref{fig:sketch_domain}).
The challenges of our experiments is threefold:
\begin{enumerate}
  \item the recovery of salt domes is recognized to be 
        difficult \citep{Faucher2019RR};
  \item we consider an initial guess that has no information 
        on the subsurface structures and where the 
        background velocity amplitude is incorrect. 
  \item we avoid the use of the (unrealistically) low frequencies
        (below $2$ \si{\Hz} in exploration seismic). 
\end{enumerate}

\subsection{Reconstruction of two-dimensional salt model}
\label{subsection:fwi_2Dsalt_reconstruction}

\begin{figure}[ht!]
  \centering
  \renewcommand{\modelfile}{trace_15s_[-1_1]e5_noise20db}
  \pgfmathsetmacro{\rcvmin} {1}  \pgfmathsetmacro{\rcvmax} {183}
  \pgfmathsetmacro{\cmin} {-0.1} \pgfmathsetmacro{\cmax}  {0.1}
  \pgfmathsetmacro{\tmin}  {0.}  \pgfmathsetmacro{\tmax} {15.0}
  \pgfmathsetmacro{\tmaxlim} {10.0}
  \subfigure[Time-domain trace with included noise.]
            {{\raisebox{0mm}{\includegraphics[scale=1]{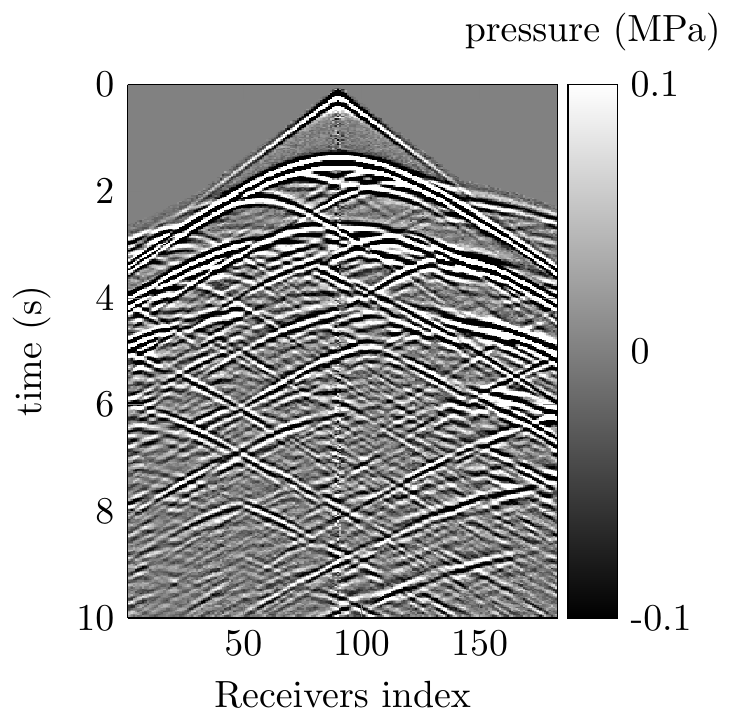}}}} \hfill
  \renewcommand{\datafile}{figures/salt_data/data_haven.txt}
  \renewcommand{\dataA}   {data_true_2hz_real}
  \renewcommand{\dataB}   {data_true_4hz_real}
  \subfigure[Real parts of the discrete Fourier transform.]
            {{\raisebox{7mm}{\includegraphics[scale=1]{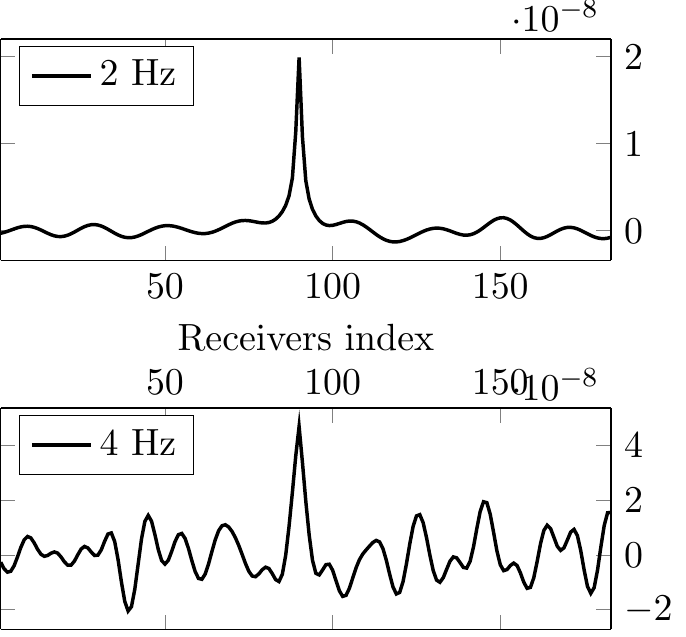}}}}
 \caption{Data associated with a single centered source. 
          The data are first generated
          in the time-domain, then we incorporate white noise
          and proceed to the Fourier transform. In this 
          experiment, the complete seismic acquisition is 
          composed of $91$ independent sources and $183$ 
          receivers for each source.}
  \label{fig:fwi_data_salt}
\end{figure}

\begin{figure}
  \centering
  \graphicspath{{figures/salt/}}
    \renewcommand{\modelfile}{cp_salt_true}
  \subfigure[Target model.]
            {\includegraphics[scale=1]{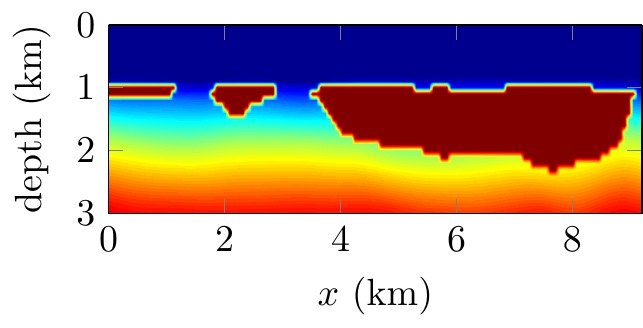}
             \label{fig:fwi_true_salt}} \hspace*{1cm}
  \renewcommand{\modelfile}{cp_salt_start}
  \subfigure[Initial guess for inversion.]
            {\includegraphics[scale=1]{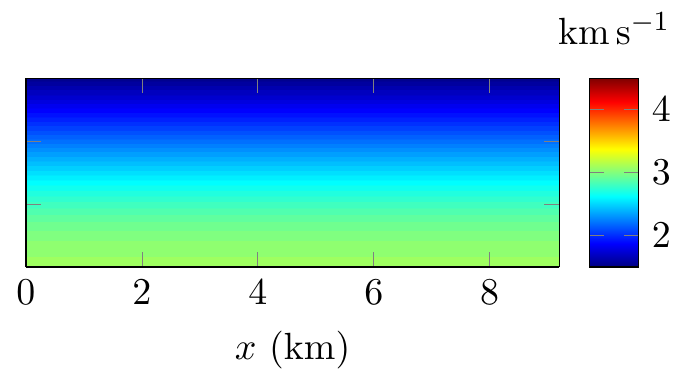}
             \label{fig:fwi_start_salt}}
  \caption{Target model and starting model for 
           FWI. The models are of size $9.2 \si{\km} 
           \times 3 \si{\km}$. The initial model 
           corresponds to a one-dimensional variation
           of low velocity.}
  \label{fig:fwi_salt_models}
\end{figure}

We first consider a two-dimensional salt model 
of size
$9.2\times3$ \si{\km}, which consists in three 
domes, see Figure~\ref{fig:fwi_true_salt}.
We generate the data 
using $91$  sources and $183$ receivers 
(i.e. data points) per source. Both devices are 
located near the surface: the sources are positioned
on a line at $20$ \si{\meter} depth and the receivers 
at $80$ \si{\meter} depth.
In order to establish a realistic situation 
despite having a synthetic experiment, 
the data are generated in the time-domain and 
we incorporate white noise in the measurements.
The level of the signal to noise ratio 
in the seismic trace is of $10$ \si{\decibel},
the noise is generated independently for all 
receivers record associated with every source.
Then we proceed to the discrete Fourier transform
to obtain the signals to be used in the reconstruction
algorithm.
In Figure~\ref{fig:fwi_data_salt}, we show 
the time-domain data with noise and the corresponding
frequency data for one source, located at $20$ \si{\meter}
depth, in the middle of the $x$-axis.

\setlength{\modelwidth} {5.95cm}
\setlength{\modelheight}{7.00cm}
\setlength{\plotwidth}  {6.20cm}
\setlength{\plotheight} {4.50cm}
\graphicspath{{figures/salt_data/}}

For the reconstruction of the salt dome model,
the starting and true model are given in 
Figure~\ref{fig:fwi_salt_models}. We do not 
assume any a priori knowledge of the contrasting 
object in the subsurface, and start with a 
one-dimensional variation, which has a drastically 
lower amplitude (i.e., the background velocity is 
unknown).

\setlength{\modelwidth} {7.00cm}
\setlength{\modelheight}{3.50cm}

\subsubsection{Fixed decomposition, single frequency reconstruction}

\begin{figure}
  \centering
  \graphicspath{{figures/salt_fwi_dataHnoise20db/00_classic/start2hz/}}
  \renewcommand{\modelfile}{cp_2hz_classic}
  \includegraphics[scale=1]{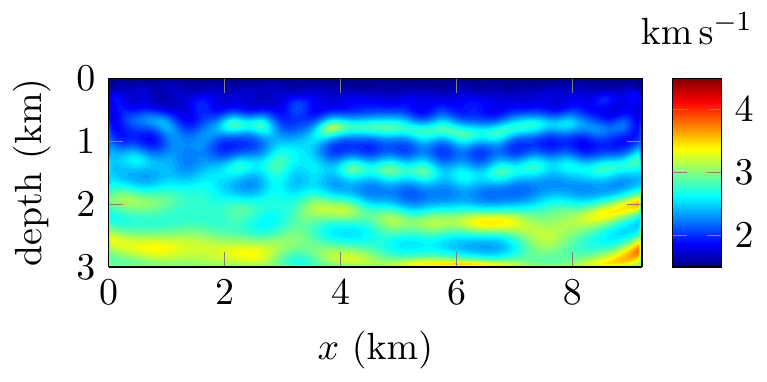}
  \caption{Reconstruction of the salt velocity model 
           from the starting medium Figure~\ref{fig:fwi_start_salt}
           using $2$ \si{\Hz} frequency data.
           The reconstruction \emph{does not} apply eigenvector 
           decomposition. The model is parametrized following
           the domain discretization, using piecewise constant
           representation with one value per node on a 
           $\num{921}\times\num{301}$ grid.}
  \label{fig:fwi_salt_classic}
\end{figure}

\begin{figure}[ht!] \centering
  \includegraphics[scale=1]{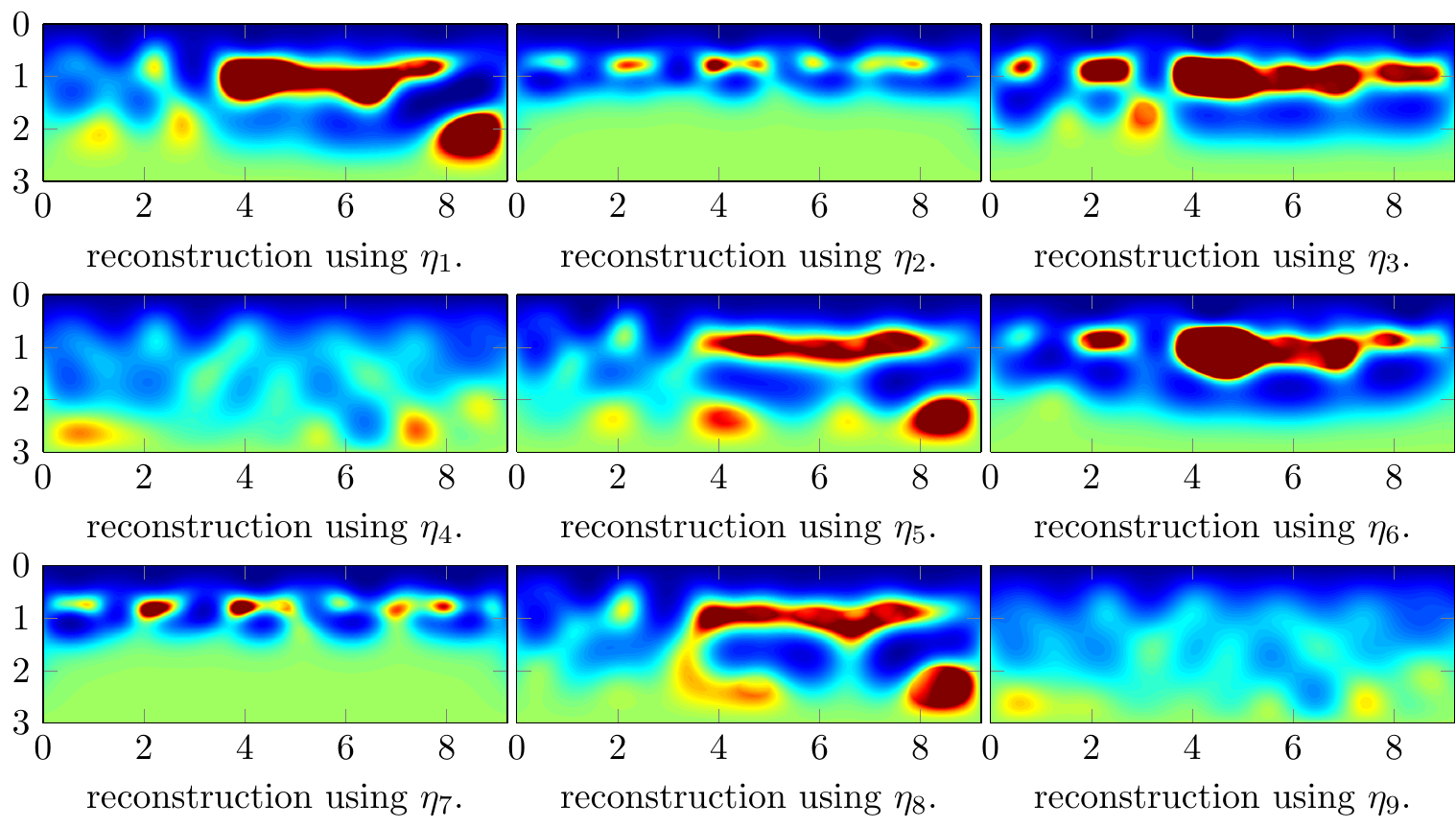}
  \caption{Reconstruction of the salt velocity model 
           from the starting medium Figure~\ref{fig:fwi_start_salt}    
           using $2$ \si{\Hz} frequency data.          
           The eigenvector decomposition employs $N=50$ 
           and the formulations of $\eta$ from 
           Table~\ref{table:list_of_formula}. It uses the same
           color scale as Figure~\ref{fig:fwi_salt_models}.}
  \label{fig:fwi_salt_50ev_2hz}
\end{figure}

\begin{figure}
  \centering
  \renewcommand{\modelfile} {cp_2hz-30iter_100ev_geman}  
  \subfigure[Reconstruction using $\eta_3$.]
            {\includegraphics[scale=1]{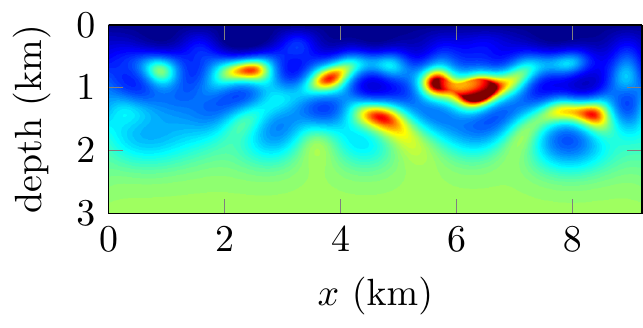}} \hspace*{1cm}
  \renewcommand{\modelfile} {cp_2hz-30iter_100ev_rudin}   
  \subfigure[Reconstruction using $\eta_8$ (Total Variation).]
            {\includegraphics[scale=1]{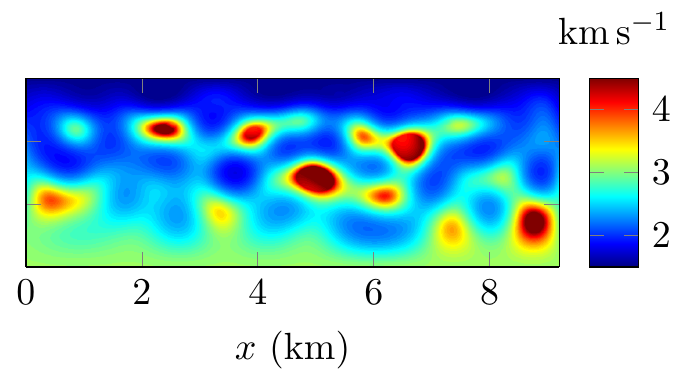}}
  \caption{Reconstruction of the salt velocity model 
           from the starting medium Figure~\ref{fig:fwi_start_salt}         
           using $2$ \si{\Hz} frequency data.          
           The eigenvector decomposition employs $N=100$.
           The comparison of all formulations of $\eta$ 
           (Table~\ref{table:list_of_formula}) is
           pictured in Appendix~\ref{appendix:additional_figures},
           Figure~\ref{fig:fwi_salt_100ev_2hz}.}
  \label{fig:fwi_salt_100ev_2hz_only2}
\end{figure}

We first only use $2$ \si{\Hz} frequency data, 
and perform $180$ iterations for the minimization. 
In Figure~\ref{fig:fwi_salt_classic}, we show 
the reconstruction where the decomposition has 
not been employed, i.e. the model representation 
follows the original piecewise constant decomposition 
of the model (one value per node).
In Figure~\ref{fig:fwi_salt_50ev_2hz},  
we compare the reconstruction using 
Algorithm~\ref{algo:fwi} for the different 
formulations of $\eta$ given in 
Table~\ref{table:list_of_formula}, 
using a fixed $N=50$. 
We use $N=100$ for Figure~\ref{fig:fwi_salt_100ev_2hz_only2}. 
For the sake of clarity, we focus on $\eta_3$
(the most effective formulation), and $\eta_8$
(which relates to the Total Variation 
regularization)
and move the complete pictures with comparison of 
all formulations in 
Appendix~\ref{appendix:additional_figures}, 
Figure~\ref{fig:fwi_salt_100ev_2hz} for $N=100$.

\setlength{\modelwidth} {6.50cm}
\setlength{\modelheight}{3.25cm}
\graphicspath{{figures/salt_fwi_dataHnoise20db/01_fixed-basis_single-freq_fixed-N/start2hz/}}
\renewcommand{\modelfileA} {cp_2hz_50ev_perona1_180iter}  
\renewcommand{\modelfileB} {cp_2hz_50ev_perona2_180iter}  
\renewcommand{\modelfileC} {cp_2hz_50ev_geman_180iter}    
\renewcommand{\modelfileD} {cp_2hz_50ev_green_180iter}    
\renewcommand{\modelfileE} {cp_2hz_50ev_aubert_180iter}   
\renewcommand{\modelfileF} {cp_2hz_50ev_lorentz_180iter}  
\renewcommand{\modelfileG} {cp_2hz_50ev_gaussian_180iter} 
\renewcommand{\modelfileH} {cp_2hz_50ev_rudin_180iter}    
\renewcommand{\modelfileI} {cp_2hz_50ev_tikhonov_180iter} 

\setlength{\modelwidth} {7.00cm}
\setlength{\modelheight}{3.50cm}
\graphicspath{{figures/salt_fwi_dataHnoise20db/01_fixed-basis_single-freq_fixed-N/start2hz/}}

We observe that
\begin{itemize}
 \item the traditional FWI algorithm (without decomposition),
       see Figure~\ref{fig:fwi_salt_classic}, fails to 
       recover any dome. It only shows some thin layers
       of increasing velocity, with amplitudes much lower 
       than the original ones.
 \item The decomposition using $N=50$ is able 
       to discover the largest object with 
       formulation $\eta_1$, $\eta_3$, $\eta_5$, 
       $\eta_6$ and $\eta_8$,
       see Figure~\ref{fig:fwi_salt_50ev_2hz}.
       The best result is given by $\eta_3$ 
       which recovers the three domes; while
       $\eta_1$, $\eta_5$ and $\eta_8$ show 
       artifacts in the lower right corner. 
       The other decompositions fail. 
       We note that, due to the lack of velocity background 
       information, the positions of the domes are slightly 
       above the correct ones to compensate for the
       low travel times.
 \item In this experiment, the method behaves much better 
       with a restrictive number of eigenvectors. 
       With $N=100$ 
       (Figures~\ref{fig:fwi_salt_100ev_2hz_only2} 
        and~\ref{fig:fwi_salt_100ev_2hz}), 
       the iterative reconstruction only results in 
       artifacts.
       The restrictive number of eigenvectors provides 
       a regularization of the problem by reducing
       the number of parameters, which is 
       crucial. For instance, the stability is known
       to deteriorate exponentially with the number of
       parameters in the piecewise constant case 
       \citep{Beretta2016}.
\end{itemize}

Opposite to the decomposition of images 
(Section~\ref{section:numeric_decomposition}), 
the quantitative reconstruction using a model 
represented with a basis of eigenvectors from
a diffusion PDE
shows drastic differences between the formulations, 
where the procedure can fail depending on the choice of $\eta$. 
In addition, the number of 
eigenvectors for the representation has to be carefully selected,
see Subsection~\ref{subsection:implementation}.



\subsubsection{Experiments with increasing $N$ and multiple frequencies}

\begin{figure}[ht!] \centering
  \graphicspath{{figures/salt_fwi_dataHnoise20db/01_fixed-basis_single-freq_fixed-N/start2hz/}}
  \renewcommand{\modelfile} {cp_2hz_50ev_geman_180iter}  
  \subfigure[Reconstruction Experiment 1 using $\eta_3$ and 180 iterations.]
            {\includegraphics[scale=1]{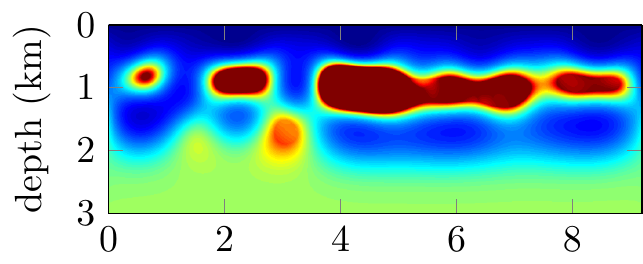}
             \label{fig:fwi_salt_experiments_only2_A}} \hspace*{1cm}
  \renewcommand{\modelfile} {cp_2hz_50ev_rudin_180iter}   
  \subfigure[Reconstruction Experiment 1 using $\eta_8$ (Total Variation) and 180 iterations.]
            {\includegraphics[scale=1]{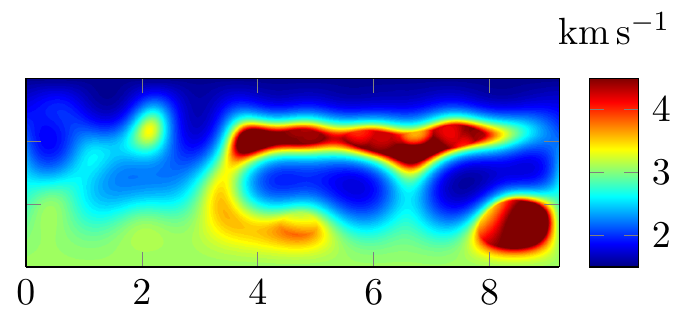}} \\[-2mm]
  \graphicspath{{figures/salt_fwi_dataHnoise20db/02_fixed-basis_single-freq_multi-N/start2hz/}}
  \renewcommand{\modelfile} {cp_2hz_50ev_geman_180iter-100ev}  
  \subfigure[Reconstruction Experiment 2 using $\eta_3$.]
            {\includegraphics[scale=1]{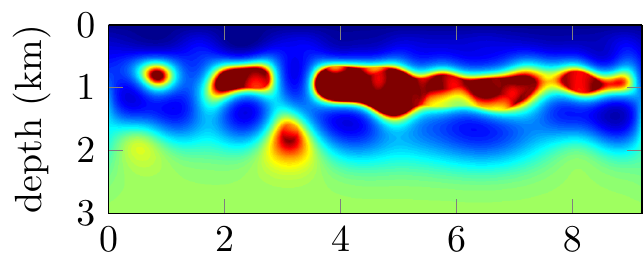}} \hspace*{1cm}
  \renewcommand{\modelfile} {cp_2hz_50ev_rudin_180iter-100ev}   
  \subfigure[Reconstruction Experiment 2 using $\eta_8$ (Total Variation).]
            {\includegraphics[scale=1]{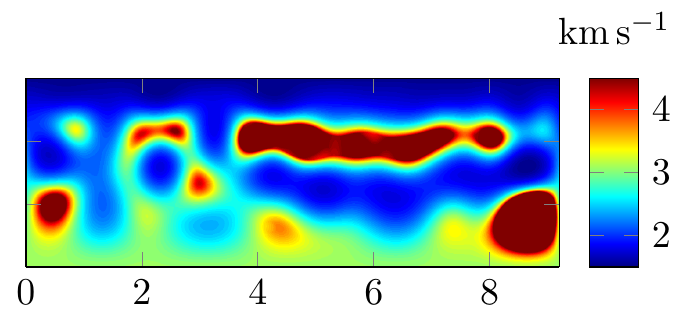}} \\[-2mm]
  \graphicspath{{figures/salt_fwi_dataHnoise20db/03_fixed-basis_multi-freq_fixed-N/start2hz/}}  
  \renewcommand{\modelfile} {cp_5hz_50ev_geman}  
  \subfigure[Reconstruction Experiment 3 using $\eta_3$.]
            {\includegraphics[scale=1]{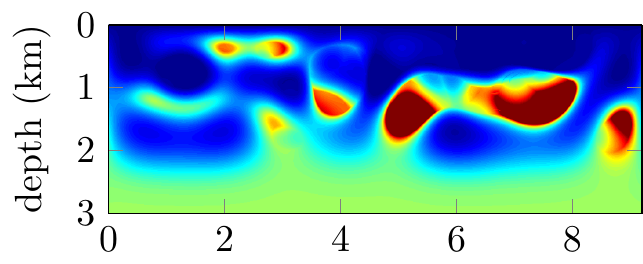}} \hspace*{1cm}
  \renewcommand{\modelfile} {cp_5hz_50ev_rudin}  
  \subfigure[Reconstruction Experiment 3 using $\eta_8$ (Total Variation).]
            {\includegraphics[scale=1]{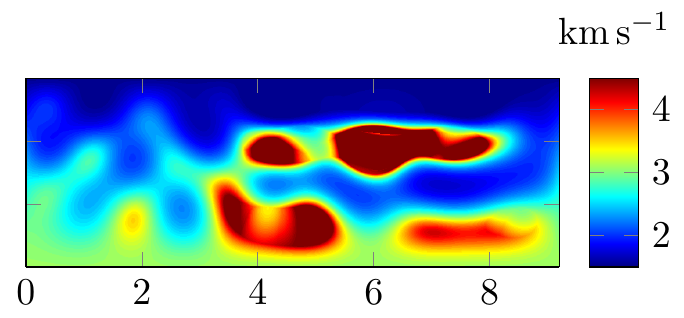}} \\[-2mm]
  \graphicspath{{figures/salt_fwi_dataHnoise20db/04_fixed-basis_multi-freq_multi-N/start2hz/}}
  \renewcommand{\modelfile} {cp_start50ev_5hz-80ev_geman}  
  \subfigure[Reconstruction Experiment 4 using $\eta_3$.]
            {\includegraphics[scale=1]{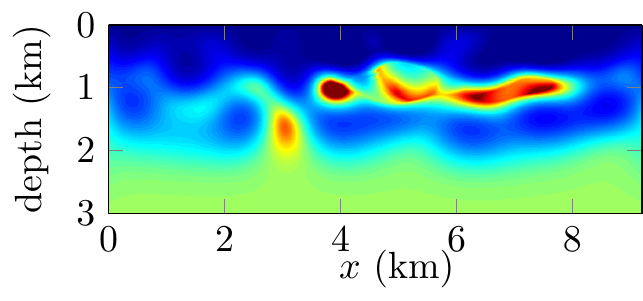}} \hspace*{1cm}
  \renewcommand{\modelfile} {cp_start50ev_5hz-80ev_rudin}   
  \subfigure[Reconstruction Experiment 4 using $\eta_8$ (Total Variation).]
            {\includegraphics[scale=1]{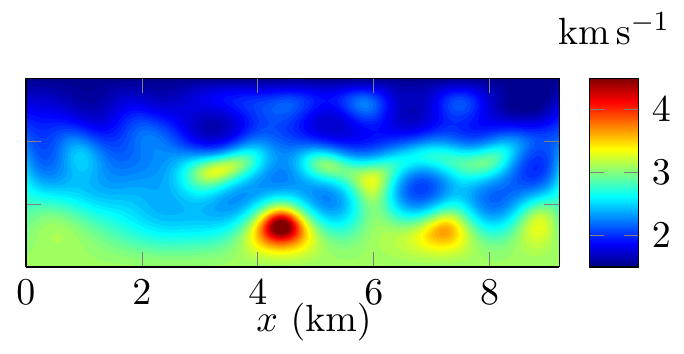}}
  \caption{Results of Experiments~2, 3 and 4 of 
           Table~\ref{table:fwi_experiments_2Dsalt}
           for the reconstruction of the salt velocity 
           model Figure~\ref{fig:fwi_salt_models}.
           The comparison of all formulations of $\eta$ 
           (Table~\ref{table:list_of_formula}) is
           pictured in Appendix~\ref{appendix:additional_figures},
           Figures~\ref{fig:fwi_salt_multiN_single-freq_2hz}, 
           \ref{fig:fwi_salt_50ev_multi-freq_5hz} and 
           \ref{fig:fwi_salt_multiN_multi-freq_5hz}.}
  \label{fig:fwi_salt_experiments_only2}
\end{figure}

We investigate the performance of the eigenvector 
decomposition for multiple frequency data, and with
progressive evolution of the number of eigenvectors
in the representation $N$. We have a total of four
different experiments, which are summarized in 
Table~\ref{table:fwi_experiments_2Dsalt}. The 
reconstructions, for $\eta_3$ and $\eta_8$, 
are shown Figure~\ref{fig:fwi_salt_experiments_only2}.
The results for all $\eta$ of 
Table~\ref{table:list_of_formula} are pictured in
Appendix~\ref{appendix:additional_figures},
Figures~\ref{fig:fwi_salt_multiN_single-freq_2hz}, 
\ref{fig:fwi_salt_50ev_multi-freq_5hz} and 
\ref{fig:fwi_salt_multiN_multi-freq_5hz}.

\begin{table} \begin{center}
\caption{List of experiments for the reconstruction of the 
         two-dimensional salt dome model 
         (Figure~\ref{fig:fwi_start_salt}). For each combination
         of frequency and associated number of eigenvectors 
         $N$ in the decomposition, $30$ iterations are 
         performed ($n_{\text{iter}}$ in Algorithm~\ref{algo:fwi}).}
\label{table:fwi_experiments_2Dsalt}
\renewcommand{\arraystretch}{1.0}
\begin{tabular}{@{}|l|c|c|c|c|}
\hline
              & \bf{frequency list} & \bf{list of} $N$ & \bf{total iterations} 
                                    & \bf{Reference} \\ \hline
Experiment~1  & $2$ \si{\Hz}        & 50                          &  180 
                                    & Figures~\ref{fig:fwi_salt_50ev_2hz}, 
                                              \ref{fig:fwi_salt_100ev_2hz}             \\ \hline
Experiment~2  & $2$ \si{\Hz}        & $\{50,60,70,80,90,100\}$    & 180 
                                    & Figures~\ref{fig:fwi_salt_experiments_only2}, 
                                              \ref{fig:fwi_salt_multiN_single-freq_2hz} \\ \hline
Experiment~3  & $\{2,3,4,5\}$\si{\Hz}
                                    & 50                          & 120 
                                    & Figures~\ref{fig:fwi_salt_experiments_only2}, 
                                              \ref{fig:fwi_salt_50ev_multi-freq_5hz}    \\ \hline
Experiment~4  & $\{2,3,4,5\}$\si{\Hz}       
                                    & $\{50,60,70,80\}$           & 120 
                                    & Figures~\ref{fig:fwi_salt_experiments_only2}, 
                                              \ref{fig:fwi_salt_multiN_multi-freq_5hz}  \\ \hline

\end{tabular} \end{center}
\end{table}


\setlength{\modelwidth} {7.00cm}
\setlength{\modelheight}{3.50cm}

From these experiments using multiple $N$ and/or
frequency contents, we observe the performance of
the method.
The best results are obtained using a 
        single $2$ \si{\Hz} frequency with 
        either progression of increasing $N$
        or constant $N$: Experiments 2 and 1. 
        The progression of $N$, Experiment~2, 
        appears the most robust.

\begin{figure}[ht!] \centering
  \pgfmathsetmacro{\tlim}  {4}
  \pgfmathsetmacro{\cmax}  {0.1}
  \graphicspath{{figures/salt_fwi_dataHnoise20db/data/}}
  \renewcommand{\modelfile}{sismogram_true_15sec_1e5_max2hz}
  \subfigure[\update{Measured seismic trace, from the 
             model of Figure~\ref{fig:fwi_true_salt},
             filtering out the frequency content above $2$ \si{\Hz}.}]
            {\includegraphics[scale=1]{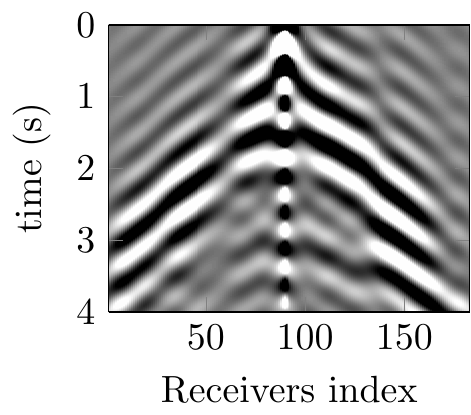}} \hfill
  \renewcommand{\modelfile}{sismogram_reconstruction_15sec_1e5_max2hz}
  \subfigure[\update{Seismic trace from the reconstructed
             velocity model of Figure~\ref{fig:fwi_salt_experiments_only2_A},
             filtering out the frequency content above $2$ \si{\Hz}.}]
            {\includegraphics[scale=1]{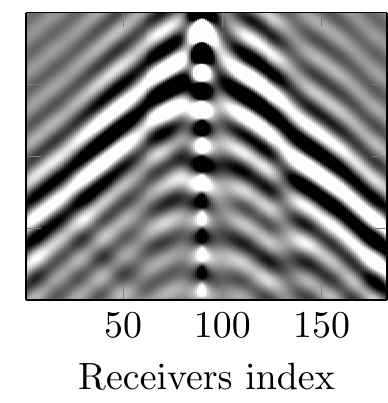}} \hfill
  \renewcommand{\modelfile}{sismogram_start_15sec_1e5_max2hz}
  \subfigure[\update{Seismic trace from the starting
             velocity model, Figure~\ref{fig:fwi_start_salt},
             filtering out the frequency content above $2$ \si{\Hz}.}]
            {\includegraphics[scale=1]{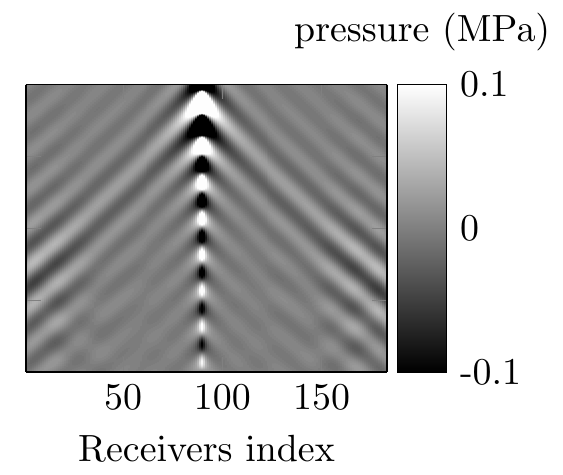}} \\[-2em]

  \pgfmathsetmacro{\cmax}  {1}
  \renewcommand{\modelfile}{sismogram_true_15sec_1e6}
  \subfigure[\update{Measured seismic trace, from the 
             model of Figure~\ref{fig:fwi_true_salt}.}]
            {\includegraphics[scale=1]{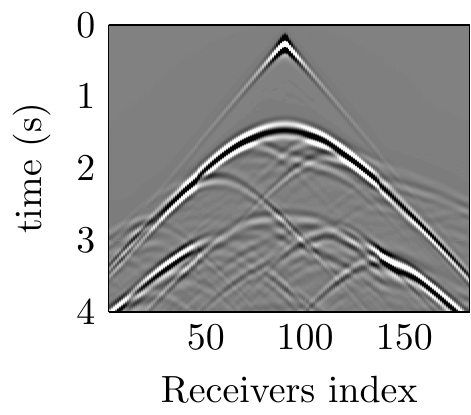}} \hfill
  \renewcommand{\modelfile}{sismogram_reconstruction_15sec_1e6}
  \subfigure[\update{Seismic trace from the reconstructed
             velocity model of Figure~\ref{fig:fwi_salt_experiments_only2_A}.}]
            {\includegraphics[scale=1]{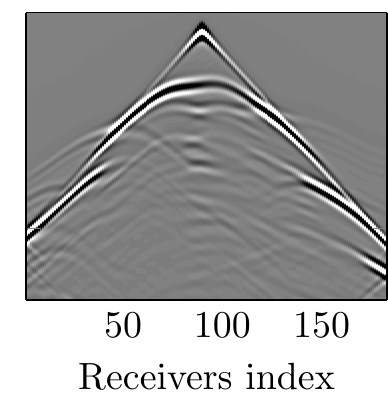}} \hfill
  \renewcommand{\modelfile}{sismogram_start_15sec_1e6}
  \subfigure[\update{Seismic trace from the starting
             velocity model, Figure~\ref{fig:fwi_start_salt}.}]
            {\includegraphics[scale=1]{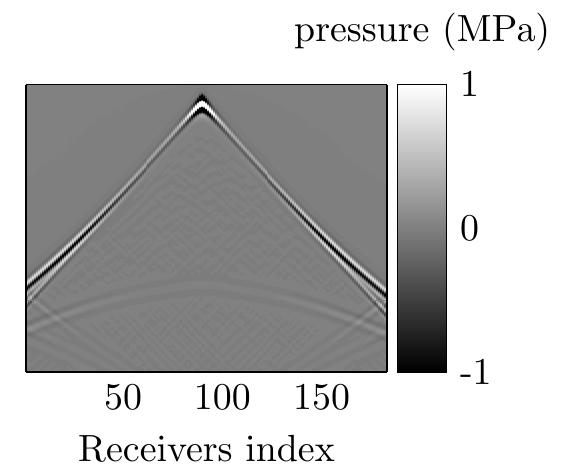}}

  \caption{\update{Comparison of the time-domain seismic traces for 
           a central shot using different velocity models of
           the salt medium reconstruction experiment.}}
  \label{fig:salt_fwi_timetrace}
\end{figure}

In this experiment, using multiple frequencies 
does not improve the results. It is probably due 
to the lack of knowledge of the velocity background 
which prevents us from recovering finer scale
(i.e., kinematic error, \cite{Bunks1995}).
In particular, local minima in the misfit functional
become more and more prominent in the high-frequency 
regime \citep{Bunks1995,Faucher2019IP}.
We illustrate the performance of the reconstruction 
with the results in the data-space: in 
Figure~\ref{fig:salt_fwi_timetrace}, we show the 
time-domain seismograms for the true, starting and 
reconstructed velocity models. 
We observe that when we filter out frequencies above 2 \si{\Hz}
(first line of Figure~\ref{fig:salt_fwi_timetrace}), the trace from 
the reconstructed model is indeed very similar to the measured one.
However, when encompassing all frequency contents (bottom line 
of Figure~\ref{fig:salt_fwi_timetrace}), important differences arise,
in particular, one can see the travel time of the first reflection 
which is earlier with the recovered model. 
This indicates that the location of the salt in 
the reconstructed velocity is above its `true' 
position.

\begin{figure}[ht!] \centering
  \renewcommand{\datafile}{figures/salt_fwi_dataHnoise20db/data/frequency_data.txt}
  \renewcommand{\dataA}   {data_true_2hz_real}
  \renewcommand{\dataB}   {data_start_2hz_real}
  \renewcommand{\dataC}   {data_rec_2hz_real}
  \pgfmathsetmacro{\ymin}{-1.5e-3} \pgfmathsetmacro{\ymax}{3e-3}
  \subfigure[Comparison of the real parts of the pressure fields at $2$ \si{\Hz}.]
            {\includegraphics[scale=1]{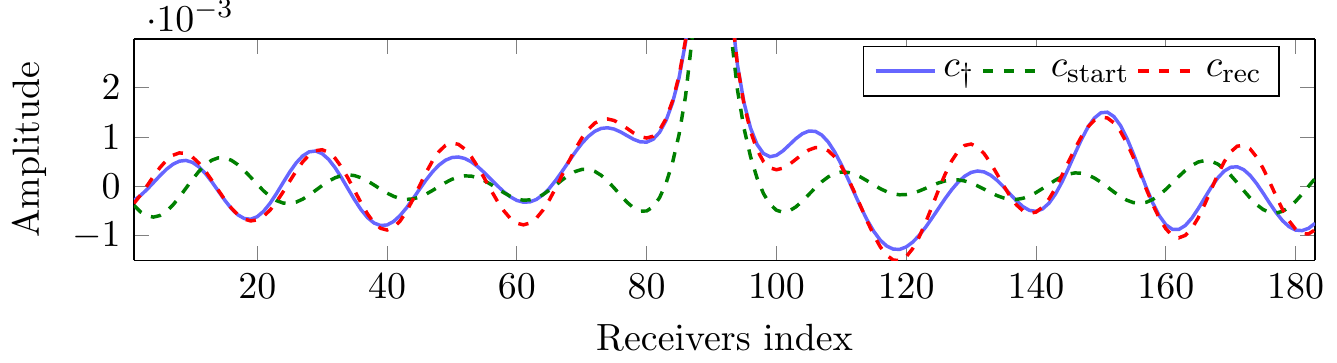}} \\
  \renewcommand{\dataA}   {data_true_4hz_real}
  \renewcommand{\dataB}   {data_start_4hz_real}
  \renewcommand{\dataC}   {data_rec_4hz_real}
  \pgfmathsetmacro{\ymin}{-2e-2} \pgfmathsetmacro{\ymax}{3e-2}
  \subfigure[Comparison of the real parts of the pressure fields at $4$ \si{\Hz}.]
            {\includegraphics[scale=1]{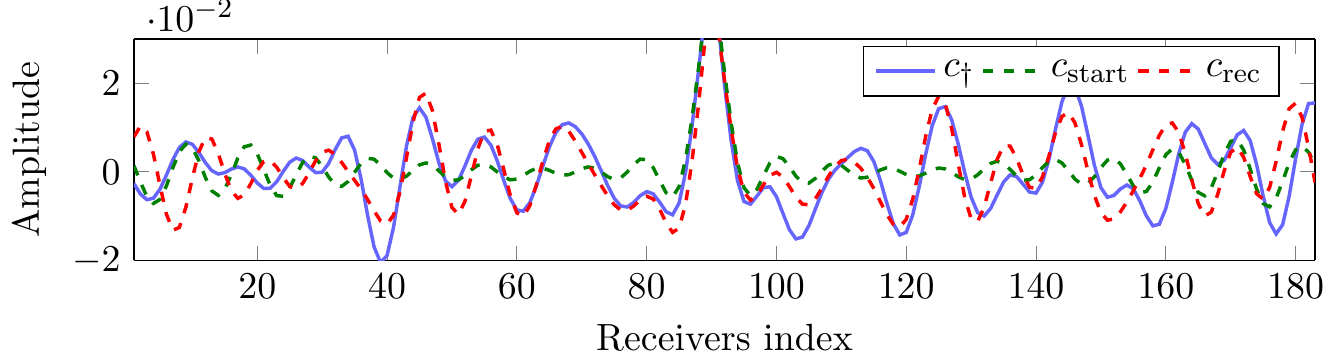}} 
  \caption{\update{Comparisons of the frequency-domain data (pressure field) at
           the receivers location for a central source. We compare the data 
           obtained from the target wave speed with salt of 
           Figure~\ref{fig:fwi_true_salt}, 
           the starting wave speed of Figure~\ref{fig:fwi_start_salt} 
           and the reconstruction of Figure~\ref{fig:fwi_salt_experiments_only2_A}.
           Those three models are respectively denoted $c_{\dagger}$,
           $c_{\text{start}}$ and $c_{\text{rec}}$.}}
  \label{fig:salt_fwi_freq_dataplot}
\end{figure}

Then, while we incorporate the higher frequency in the 
minimization procedure, the FWI is not amenable to improve 
the results (see Figure~\ref{fig:fwi_salt_experiments_only2}) 
and it is most likely due to the missing velocity background 
which is not improved during the first iterations, and 
still missing.
In Figure~\ref{fig:salt_fwi_freq_dataplot}, we show the 
frequency-domain data at $2$ and $4$ \si{\Hz}: 
the observed data at $2$ \si{\Hz} are accurately obtained
with the model reconstructed with the decomposition in eigenvectors,
which confirms the pertinence of the method. Interestingly, at 
$4$ \si{\Hz}, while the frequency is not even used in the 
inversion scheme (we only use $2$ \si{\Hz} for 
Figure~\ref{fig:fwi_salt_experiments_only2_A}), 
we already have a good correspondence near the source and only 
the parts further away show a shift. 
In order to overcome the issue of recovering the 
background velocity, one would need lower frequency 
content, or one could employ alternative strategies, 
such as the MBTT method, based upon the decomposition 
of the background velocity model and the 
reflectivity \citep{Clement2001}. 
Here, the decomposition in eigenvectors appropriately 
recover the reflectivity part (better than traditional 
FWI), but the background model remains missing.

\setlength{\modelwidth} {5.25cm}
\setlength{\modelheight}{4.50cm}

\setlength{\plotwidth} {12cm}
\setlength{\plotheight}{2.25cm}

\subsubsection{On the choice of the number of eigenvectors}

\begin{figure}[ht!] \centering
  \pgfmathsetmacro{\ymin}    {0.17}
  \pgfmathsetmacro{\ymax}    {1.1}
  \pgfmathsetmacro{\xstart}  {0}
  \pgfmathsetmacro{\scale} {5240864} 
  \renewcommand{\datafile}{figures/salt_fwi_dataHnoise20db/misfit/geman_1e-9_multi-freq_single-N.txt}
  \renewcommand   {\dataA}  {10ev}
  \renewcommand   {\dataB}  {25ev}
  \renewcommand   {\dataC}  {50ev}
  \renewcommand   {\dataD} {100ev}
  \renewcommand   {\dataE} {250ev}
  \pgfmathsetmacro{\xend}   {29}
  \subfigure[Using fixed $N$.]{\includegraphics[scale=1]{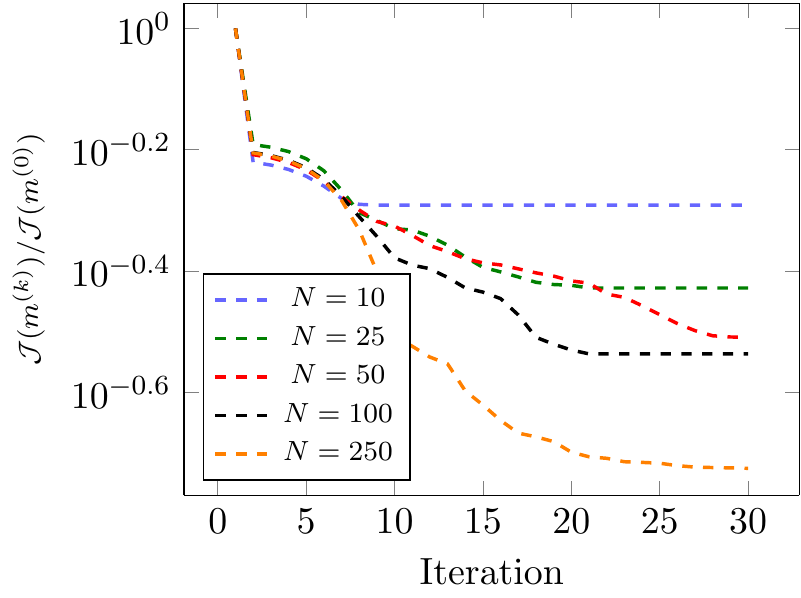}
                         \label{fig:fwi_salt_misfit_N_a}} \hfill
  \renewcommand{\datafile}{figures/salt_fwi_dataHnoise20db/misfit/geman_1e-9_single-freq_multi-N.txt}
  \renewcommand   {\dataA} {50ev}
  \pgfmathsetmacro{\xend}  {179}
  \subfigure[Using progression of $N$ 
             (from $50$ to $100$), 
             every $30$ iterations,
             following Experiment~2 of
             Table~\ref{table:fwi_experiments_2Dsalt}.]
             {\includegraphics[scale=1]{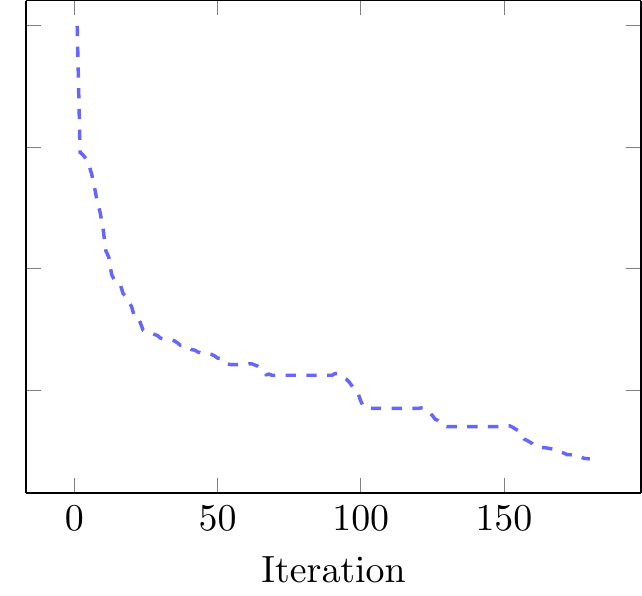}
                        \label{fig:fwi_salt_misfit_N_b}}
  \caption{Evolution of the misfit functional (scaled with
           the first iteration value) with iterations depending
           on the choice of $N$, using $2$ \si{\Hz} frequency.}
  \label{fig:fwi_salt_misfit_N}
\end{figure}

We have shown in Figures~\ref{fig:fwi_salt_50ev_2hz} 
and~\ref{fig:fwi_salt_100ev_2hz_only2} that one should
take an initial $N$ relatively low for the reconstruction
algorithm to succeed.
It remains to verify if the appropriate
$N$ can be selected `\emph{a priori}', or based upon
minimal experiments. In Figure~\ref{fig:fwi_salt_misfit_N_a},
we show the evolution of the misfit functional with 
thirty iterations, for different values of $N$, 
from $10$ to $250$. 
We compare, in Figure~\ref{fig:fwi_salt_misfit_N_b}, with
the progression of $N$, which follows Experiment~2 of
Table~\ref{table:fwi_experiments_2Dsalt}.

\setlength{\plotwidth} {6.25cm}
\setlength{\plotheight}{5.0cm}

From Figure~\ref{fig:fwi_salt_misfit_N_a}, we see
that all of the choices of $N$ have the same pattern: 
first the decrease of the functional and then its 
stagnation.
We notice that the good choice for $N$ 
\emph{is not} reflected by the misfit functional. 
Indeed, it shows lower error for larger $N$, 
while they are shown to result in erroneous 
reconstructions (Figure~\ref{fig:fwi_salt_100ev_2hz_only2}
compared to Figure~\ref{fig:fwi_salt_50ev_2hz}).
It is most likely that using larger $N$ leads to 
local minima and/or deteriorates the stability
(see, for the piecewise constant case, \cite{Beretta2016}).
It results in the false impression 
(from the misfit functional point of view) that it would improve
the efficiency of the method.
Using a progression of $N$, Figure~\ref{fig:fwi_salt_misfit_N_b},
eventually gives the same misfit functional value than
the large $N$, but it needs more iterations. 
This increase of iterations and `slow' convergence
is actually required, because it leads to an 
appropriate reconstruction, see Figure~\ref{fig:fwi_salt_experiments_only2}.

Therefore, we \emph{cannot} anticipate a good
choice for $N$ a priori (with a few evaluation
of the misfit functional). the guideline we 
propose, as a safe and robust approach, is the 
progression of increasing $N$, from low to high:
it costs more in terms of iterations, but it
converges properly.

\subsection{Reconstruction of the SEAM Phase I model}

\begin{figure}[ht!] \centering
  \graphicspath{{figures/seam/main/}}
  \renewcommand{\modelfile}{cp_true}
  \subfigure[Target model.]
            {\includegraphics[scale=1]{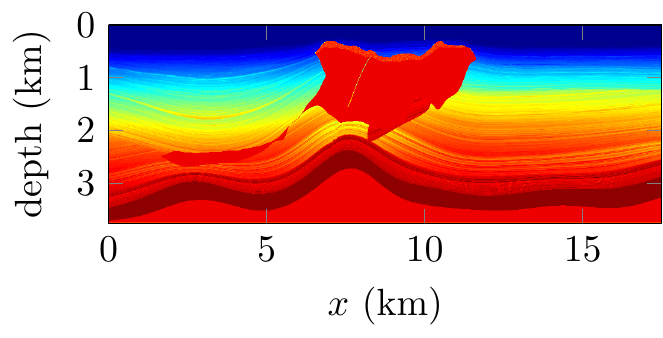}
             \label{fig:seam_fwi_main_true}}  \hfill
  \renewcommand{\modelfile}{cp_start}
  \subfigure[Initial guess for inversion.]
            {\includegraphics[scale=1]{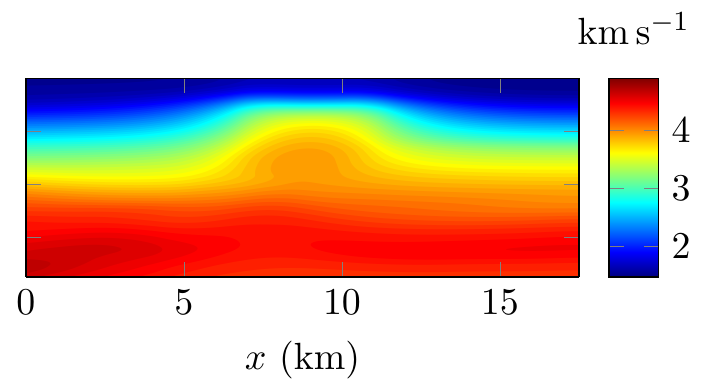}
             \label{fig:seam_fwi_main_start}}
  \caption{Target model and starting model for 
           FWI, of size $17.5 \times 3.75$ \si{\km}.}
  \label{fig:seam_fwi_main}
\end{figure}

We now consider the recovery of the SEAM Phase I model, 
which is expected to be more challenging as it contains 
both a salt-dome and sub-salt layers.
The starting model for the reconstruction is shown in 
Figure~\ref{fig:seam_fwi_main}, where, for the sake of 
clarity, we also picture the reference model which was
previously used for the decomposition. 
This medium is of size $17.5 \times 3.75$ \si{\km} and
the starting model we use is a smooth version of the 
reference one, where the contrasting objects and layers 
are missing.

\setlength{\modelwidth} {7.20cm}
\setlength{\modelheight}{3.60cm}

We follow the same configuration as in the previous experiment: 
the data are generated in the time-domain, and noise is 
incorporated to the seismograms using a signal-to-noise ratio
of 10 \si{\decibel}.
Next, we proceed with the discrete Fourier-transform to employ 
the iterative procedure with time-harmonic waves.
In this experiment, the smallest available frequency is $2$ \si{\Hz}.

\begin{figure}[ht!] \centering
  \graphicspath{{figures/seam/fwi/}}
  \renewcommand{\modelfile}{cp_start2hz_2hz_classic}
  \subfigure[Reconstruction after $180$ iterations 
             without using the eigenvector decomposition.]
            {\includegraphics[scale=1]{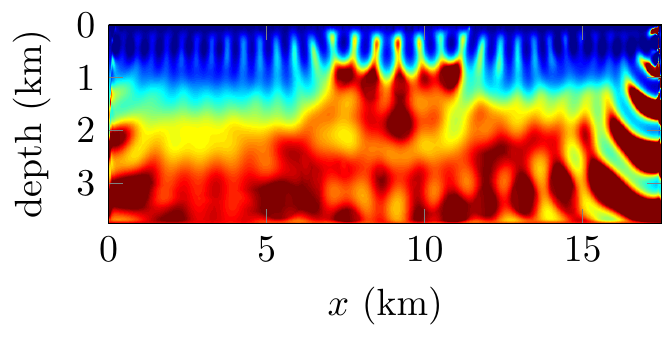}} \hfill
  \renewcommand{\modelfile}{cp_start2hz-25ev_2hz_geman_80ev_180iter}
  \subfigure[Reconstruction using eigenvector decomposition 
             with $\eta_3$ and a sequence $N=\{25, 35, 50, 60, 70, 80\}$
             with $30$ iterations for each, i.e. a total of 180 iterations.]
            {\includegraphics[scale=1]{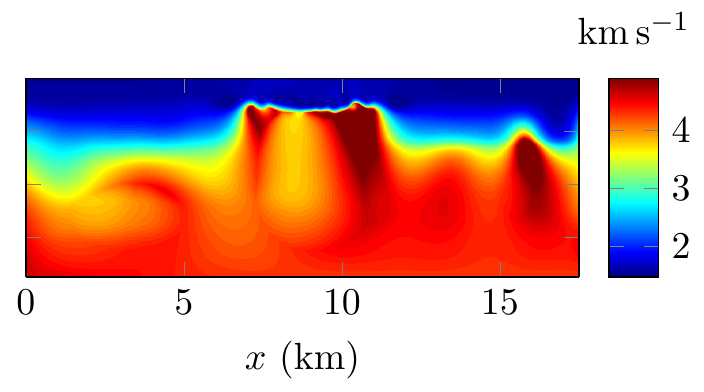}
             \label{fig:seam_fwi_2hz_ev}} 
  \caption{Reconstruction of the SEAM phase I model from
           an initially smooth model, 
           see Figure~\ref{fig:seam_fwi_main},
           using $2$ \si{\Hz} frequency data.
           }  \label{fig:seam_fwi_2hz}
\end{figure}

For the sake of conciseness, we only present the reconstruction
results with the the representation using $\eta_3$ 
(which was the most efficient).
We follow a slow increase of $N$ in a fixed basis 
(analogous to Experiment~2 of Table~\ref{table:fwi_experiments_2Dsalt})
which was the more stable approach.
The reconstruction using $2$ \si{\Hz} is shown in 
Figure~\ref{fig:seam_fwi_2hz}.

\setlength{\modelwidth} {7.20cm}
\setlength{\modelheight}{3.60cm}

\begin{figure}[ht!]\centering
  \renewcommand{\datafile}{figures/seam/data_frequency/data_2hz.txt}
  \renewcommand{\dataA}   {true_real}
  \renewcommand{\dataB}   {start_real}
  \renewcommand{\dataC}   {ev_real}
  \pgfmathsetmacro{\ymin}{-1.5e-3} \pgfmathsetmacro{\ymax}{2e-3}
  \subfigure[Comparison of the real parts of the pressure fields.]
            {\includegraphics[scale=1]{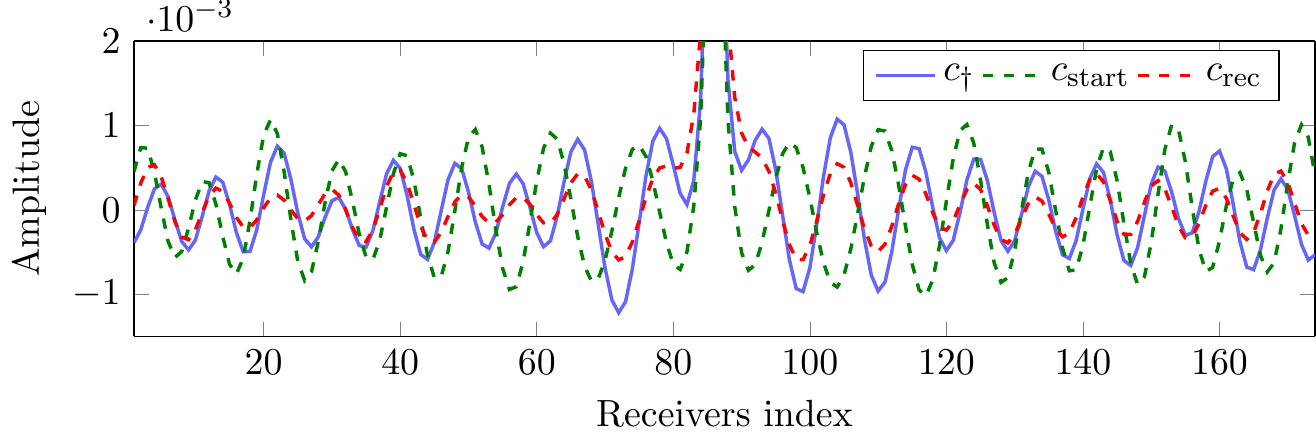}} \\
  \renewcommand{\dataA}   {true_imag}
  \renewcommand{\dataB}   {start_imag}
  \renewcommand{\dataC}   {ev_imag}
  \pgfmathsetmacro{\ymin}{-1.2e-3} \pgfmathsetmacro{\ymax}{1.5e-3}
  \subfigure[Comparison of the imaginary parts of the pressure fields.]
            {\includegraphics[scale=1]{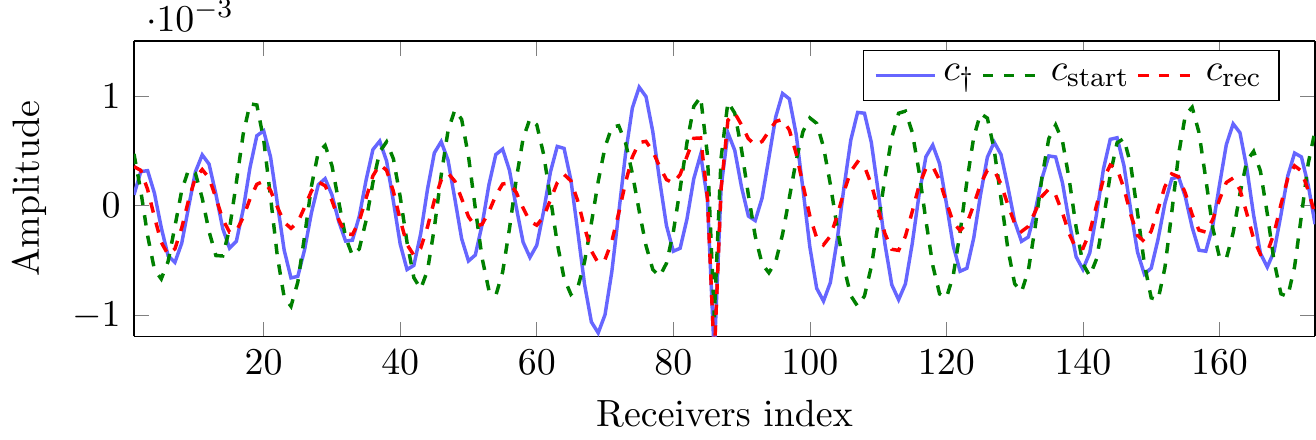}} 
  \caption{Comparisons of the frequency-domain data (pressure field) at
           the receivers location for a centrally located source at 
           $2$ \si{\Hz} frequency for the SEAM Phase I wave speed 
           of Figure~\ref{fig:seam_fwi_main_true}, 
           the starting wave speed of Figure~\ref{fig:seam_fwi_main_start} 
           and the reconstruction using the eigenvector decomposition, 
           Figure~\ref{fig:seam_fwi_2hz_ev}.
           Those three models are respectively denoted $c_{\dagger}$,
           $c_{\text{start}}$ and $c_{\text{rec}}$.}
  \label{fig:seam_fwi_2hz_dataplot}
\end{figure}

We observe that the standard FWI algorithm gives
artifacts over the medium, with oscillatory patterns,
in particular on the sides. 
On the other hand, the reconstruction using the representation 
based upon the eigenvector decomposition is stable, and is able 
to capture accurately the upper boundary of the salt dome.
Because of the frequency employed, only the long wavelengths 
are recovered at this stage. 
In Figure~\ref{fig:seam_fwi_2hz_dataplot}, we further 
illustrate the recovery by showing the frequency-domain 
data at $2$ \si{\Hz} frequency using the different wave 
speed models.
We can see that the data from the starting model are 
out of phase as soon as we move away from the source,
with cycle-skipping effects.
However, the reconstruction using the eigenvector decomposition 
is able to retrieve this information and accurately capture 
the oscillations of the signal, and only the amplitude is 
inaccurate.

\setlength{\plotwidth} {12cm}
\setlength{\plotheight}{3cm}

\begin{figure}[ht!] \centering
  \graphicspath{{figures/seam/fwi/}}
  \renewcommand{\modelfile}{cp_start2hz-after-geman_10hz_classic}
  \includegraphics[scale=1]{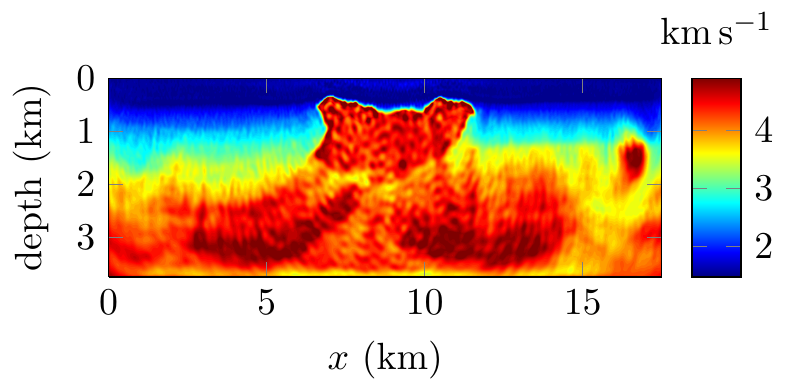}
  \caption{Reconstruction of the SEAM phase I model starting
           from the reconstruction obtained with eigenvector
           decomposition of Figure~\ref{fig:seam_fwi_2hz_ev},
           using data of frequency from $2$ to $10$ \si{\Hz},
           with 30 iterations per (sequential) frequency.}
  \label{fig:seam_fwi_final}
\end{figure}

We now continue the procedure using increasing frequencies. 
Because of the smoothing effect of the decomposition, 
we employ the algorithm without the eigenvector representation
(alternatively, one could use the decomposition but with large $N$).
Therefore, the velocity model obtained from the FWI with 
eigenvector decomposition of Figure~\ref{fig:seam_fwi_2hz_ev}
is used as an initial model for restarted FWI with multiple 
frequencies, from $2$ to $10$ \si{\Hz} 
(we use the sequential progression advocated by \cite{Faucher2019RR}).
Eventually, the reconstruction after $10$ \si{\Hz} iterations 
is pictured in Figure~\ref{fig:seam_fwi_final}.

\setlength{\modelwidth} {7.20cm}
\setlength{\modelheight}{3.60cm}

\begin{figure}[ht!] \centering
  \pgfmathsetmacro{\tlim}  {5}
  \graphicspath{{figures/seam/time-domain_data/}}
  \renewcommand{\modelfile}{trace_cp-true_noise_[-3_3]e5_10sec}
  \subfigure[Measured seismic trace corresponding to the SEAM Phase I model of Figure~\ref{fig:seam_fwi_main_true}.]
            {\includegraphics[scale=1]{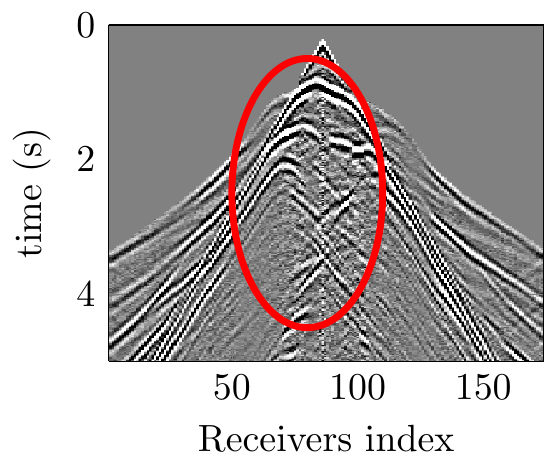}} \hspace*{.5cm}
  \renewcommand{\modelfile}{trace_cp-start_[-3_3]e5_10sec}
  \subfigure[Seismic trace corresponding to the starting model of Figure~\ref{fig:seam_fwi_main_start}.]
            {\includegraphics[scale=1]{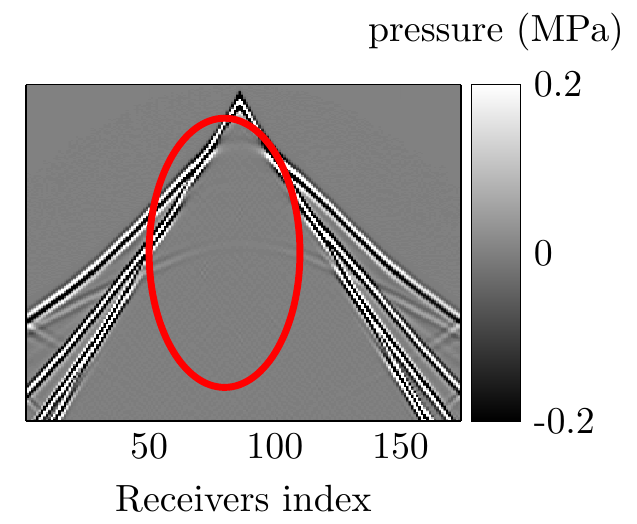}}
            
  \renewcommand{\modelfile}{trace_cp-ev-2hz_[-3_3]e5_10sec}
  \subfigure[Seismic trace corresponding to the partial reconstruction 
             using the eigenvector decomposition, Figure~\ref{fig:seam_fwi_2hz_ev}.]
            {\includegraphics[scale=1]{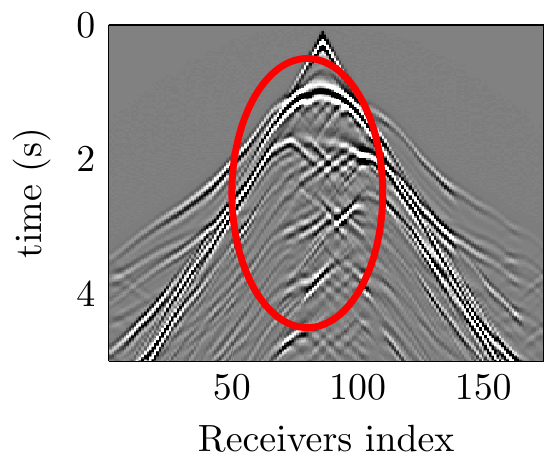}} \hspace*{.5cm}
  \renewcommand{\modelfile}{trace_cp-final-10hz_[-3_3]e5_10sec}
  \subfigure[Seismic trace corresponding to the final reconstruction 
             given Figure~\ref{fig:seam_fwi_final}.]
            {\includegraphics[scale=1]{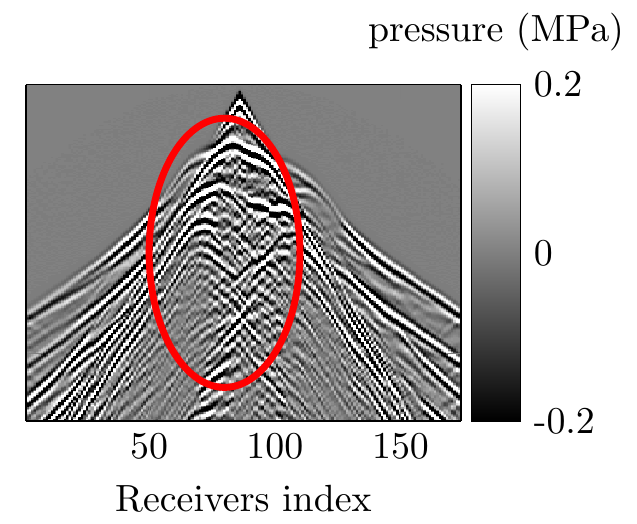}}
  \caption{Comparison of the time-domain seismic traces for a central shot
           using different velocity models.}
  \label{fig:seam_fwi_timetrace}
\end{figure}

The reconstruction is able to capture the finer 
details of the velocity model:
the salt dome is clearly defined, and the sub-salt
layer starts to appear.
Therefore, the eigenvector decomposition method 
can also serve to build initial models for the 
FWI algorithm, where it appears as an interesting 
alternative to overcome the lack of low frequency.
To illustrate the different steps of the reconstruction,
we show the time-domain seismograms at the different stages
of the reconstruction in Figure~\ref{fig:seam_fwi_timetrace}
(while the $2$ \si{\Hz} frequency-domain data 
are given in Figure~\ref{fig:seam_fwi_2hz_dataplot}).
We compare the traces resulting from the initial and true 
wave speed models, from the partial reconstruction obtained 
with the eigenvector decomposition (Figure~\ref{fig:seam_fwi_2hz_ev}), 
and from the final reconstruction (Figure~\ref{fig:seam_fwi_final}).

\setlength{\modelwidth} {6.00cm}
\setlength{\modelheight}{5.00cm}

The trace that uses the starting model mainly shows the 
first arrivals, with some minor reflections coming from 
the smoothing of the salt dome. 
The one using the partial reconstruction with eigenvenctor
decomposition and only the $2$ \si{\Hz} data accurately 
resolves the main multiple reflections between the 
salt upper boundary and the surface (the thick line in the 
center of the red ellipses in Figure~\ref{fig:seam_fwi_timetrace})
but it misses the events of smaller importance.
Eventually, the final reconstruction, that uses up to $10$ \si{\Hz}
data contents is able to reproduce some of the events of smaller 
amplitudes. We also note that the amplitude of the trace for
the final reconstruction is slightly high, which indicates
that the contrasts are even too strong.

\subsection{Implementation of the method}
\label{subsection:implementation}

The eigenvector decomposition for the model 
representation depends on different choices of 
parameters, and it is not trivial to 
efficiently implement the method in the
iterative reconstruction procedure.
From the experiments we have carried out, 
we propose the following strategy.
\begin{enumerate}
  \renewcommand{\theenumi}{(\arabic{enumi})}
  \item Regarding the choice of weighting coefficient, 
        $\eta_3$ from~\cite{Geman1992} is the most
        efficient in these applications, and supersedes the 
        standard Total Variation approach (i.e., $\eta_8$).
\end{enumerate}

Secondly, the difficulty resides in the number 
of eigenvectors to take for the decomposition: $N$. 
More important, it appears that the misfit functional
\emph{does not} provide us with a good indication 
(see Figure~\ref{fig:fwi_salt_misfit_N}).
\begin{enumerate} \setcounter{enumi}{1}
  \renewcommand{\theenumi}{(\arabic{enumi})}
  \item The number of eigenvector $N$ for the decomposition
        should initially takes a low value, and progressively 
        evolves to higher values with iterations.
        It may not be the fastest convergence, but it 
        is the most robust approach.
\end{enumerate}

Finally, the reconstruction can serve as an initial model for
multi-frequency data:
\begin{enumerate} \setcounter{enumi}{2}
  \renewcommand{\theenumi}{(\arabic{enumi})}
  \item the (partial) reconstruction with eigenvector
        decomposition is used as a starting model 
        for multi-frequency algorithm. 
        It allows the recovery of the finer details,
        which depend on the smaller wavelengths and 
        where the smoothing effect is misleading.        
\end{enumerate}
Namely, the decomposition is particularly efficient
to overcome the lack of low-frequency in the data.

\subsection{Three-dimensional experiment}
\label{subsection:3D_fwi}

\begin{figure}
\centering
\renewcommand{\modelfile}{cp_true}
\includegraphics[scale=1]{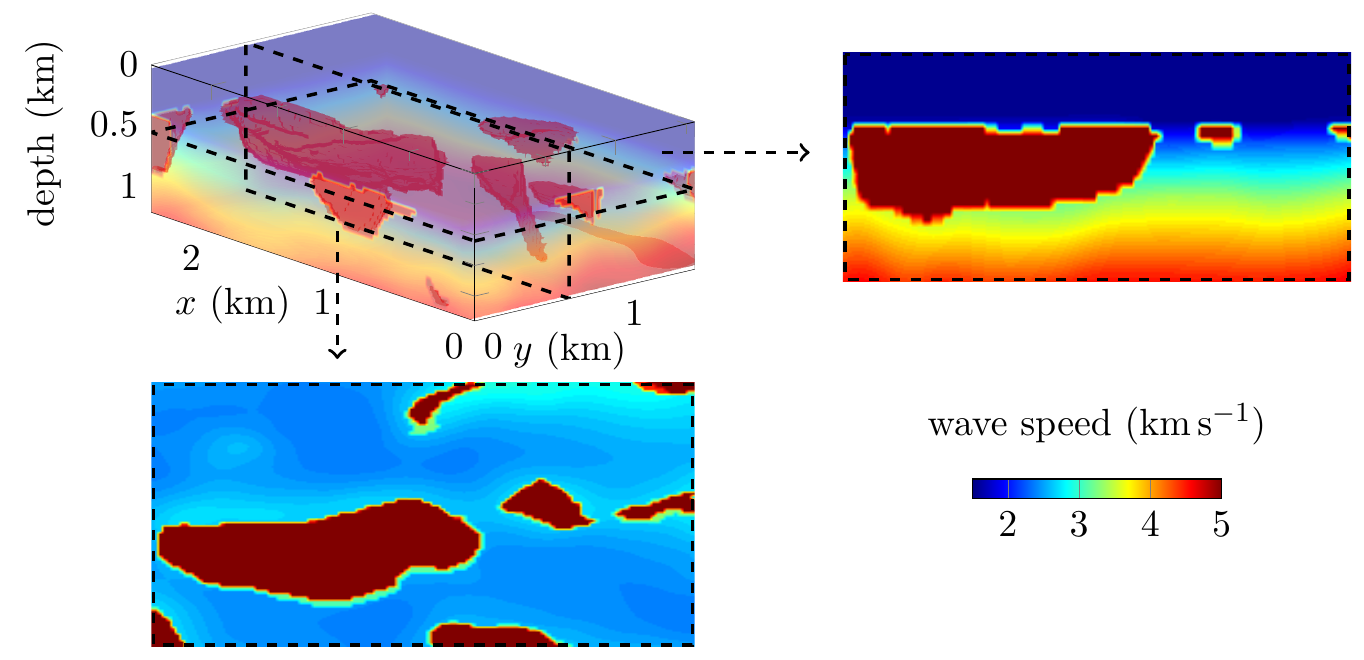}
\caption{Three-dimensional model incorporating 
         contrasting objects. The domain is of
         size $2.46 \times 1.56 \times 1.2$ \si{\km}.
         We highlight a horizontal section at $550$ \si{\meter} 
         depth and vertical section at $y=670$ \si{\meter}.}
\label{fig:3dsalt}
\end{figure}

\begin{figure}
\centering
\renewcommand{\modelfile}{cp_start}
\includegraphics[scale=1]{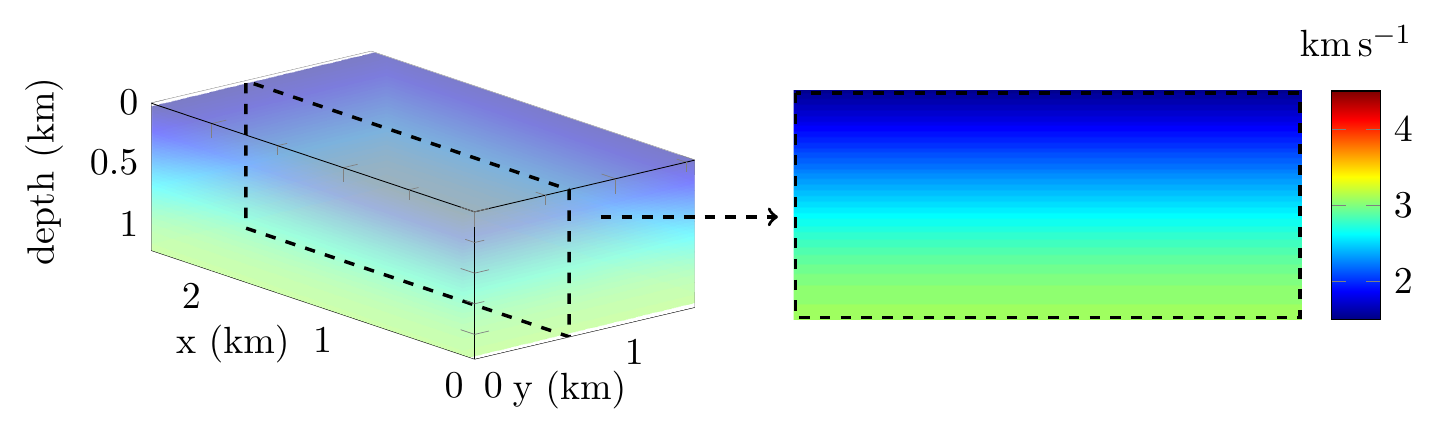}
\caption{Initial model taken for the reconstruction 
         of the three-dimensional medium with 
         vertical section at $y=670$ \si{\meter}.
         It consists in a one-dimensional variation 
         in depth.}
\label{fig:3dsalt_start}
\end{figure}

The method extends readily for three-dimensional model
reconstruction, simply incurring a larger computational 
cost (as larger matrices are involved for the eigenvector
decomposition and the forward problem discretization). 
We proceed with a three-dimensional experiment, 
where we consider a subsurface medium of size 
$2.46 \times 1.56 \times 1.2$ \si{\km}, encompassing 
several salt domes, illustrated in Figure~\ref{fig:3dsalt}.
The seismic acquisition consists in \num{96} sources, 
positioned on a two-dimensional plane at \num{10} \si{\meter}
depth; one thousand receivers are positioned at 
\num{100} \si{\meter} depth.
Similar to the previous experiments, the data are first
generated in the time-domain and we incorporate noise 
before we proceed to the Fourier transform. 
Figure~\ref{fig:3dsalt:data_traces} shows the time-domain 
data associated with a centrally located source, and the 
corresponding Fourier transform at $5$ \si{\Hz} frequency.
For the reconstruction, we start with a one-dimensional 
variation, in depth only, where none of the objects is
intuited, see Figure~\ref{fig:3dsalt_start}, and the velocity 
background is incorrect.


\graphicspath{{figures/salt3d_fwi/init/}}
\setlength{\modelwidth}{7.10cm}

\setlength{\modelwidth}{7.1cm}
\graphicspath{{figures/salt3d_fwi/init/}}

\begin{figure}
\centering
  \subfigure[Time-domain measurements with noise.]
            {\includegraphics[scale=1]{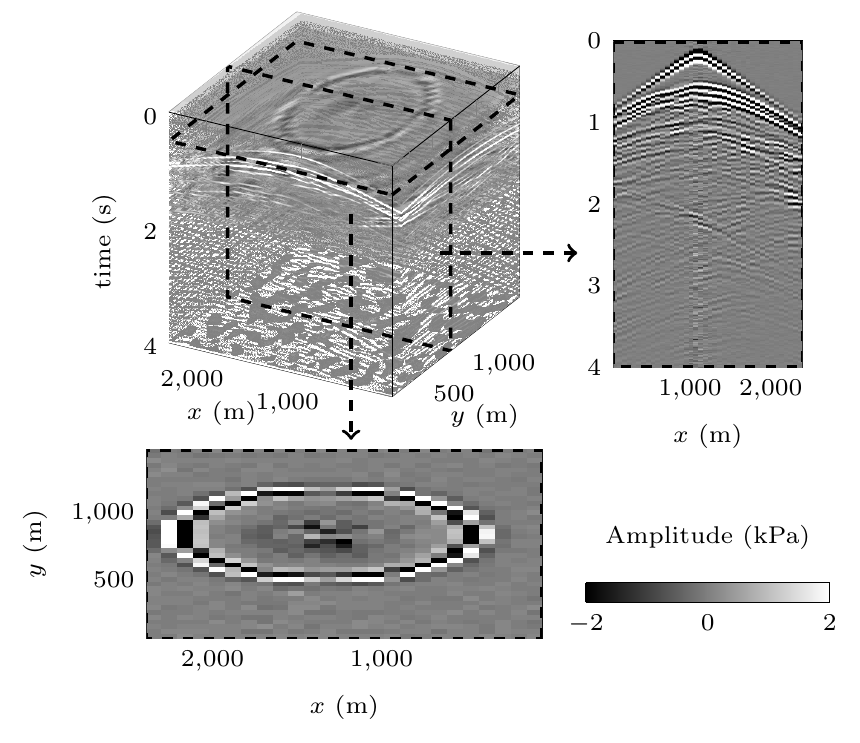}}
  \graphicspath{{figures/salt3d_fwi/data-trace/}}
  \renewcommand{\modelfile}{fourierT_5hz_scale2e2}
  \setlength{\modelwidth} {5cm}
  \setlength{\modelheight}{7cm}
  \subfigure[Fourier transform at $5$ \si{\Hz}.]
            {\raisebox{.5cm}{{\includegraphics[scale=1]{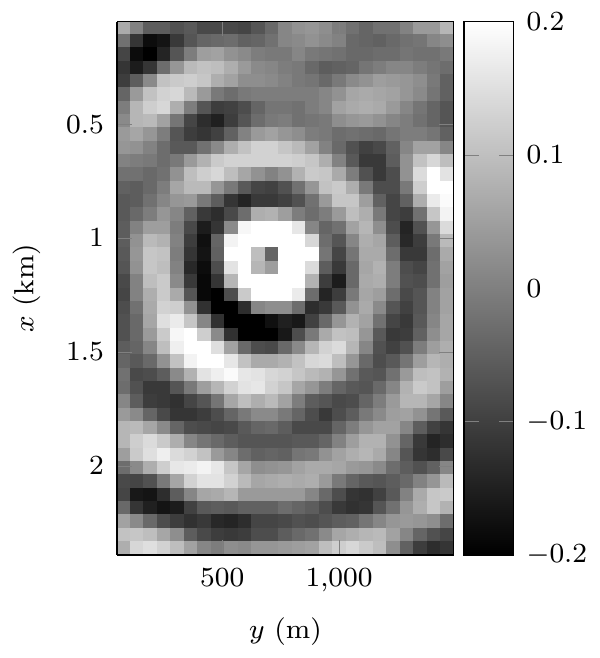}}}}
  \caption{Time-domain data and corresponding Fourier transform 
           at $5$ \si{\Hz} frequency. The three-dimensional trace
           corresponds with the evolution receivers recordings 
           (positioned on a $2$D map in the $x$--$y$ plane) with time.
           There are $1000$ data points per time step (i.e. 1000 receivers
           on the domain) and we highlight sections at fixed time 
           ($0.5$ \si{\second}) and for a line of receivers (positioned 
           at $y=710$ \si{\meter}).}
  \label{fig:3dsalt:data_traces}
\end{figure}

\begin{figure}
\centering
\renewcommand{\modelfile}{cp_5hz_classic}
\includegraphics[scale=1]{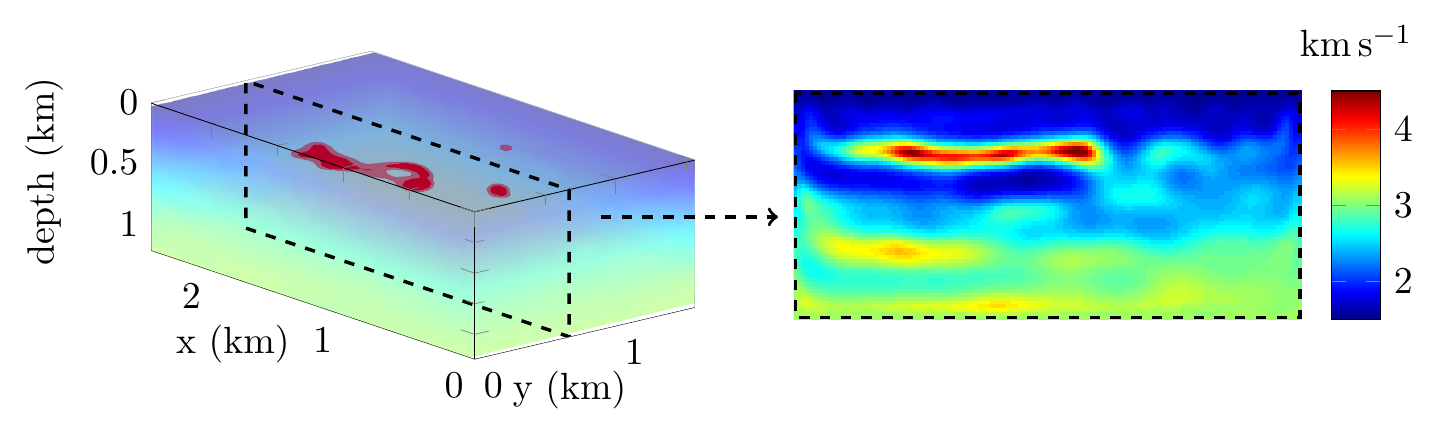}
\caption{Velocity model reconstruction  with 
         vertical section at $y=670$ \si{\meter},
         from the starting medium Figure~\ref{fig:3dsalt_start}
         using $5$ \si{\Hz} frequency data.
         The reconstruction \emph{does not} apply the eigenvector 
         decomposition. The model is parametrized following
         the domain discretization, using piecewise constant
         representation with one value per node on a 
         $\num{124}\times\num{79}\times\num{61}$ grid.}
\label{fig:3dsalt_fwi_classic}
\end{figure}

\begin{figure}
\centering
\renewcommand{\modelfile}{cp_set1-5hz_50ev_geman}
\includegraphics[scale=1]{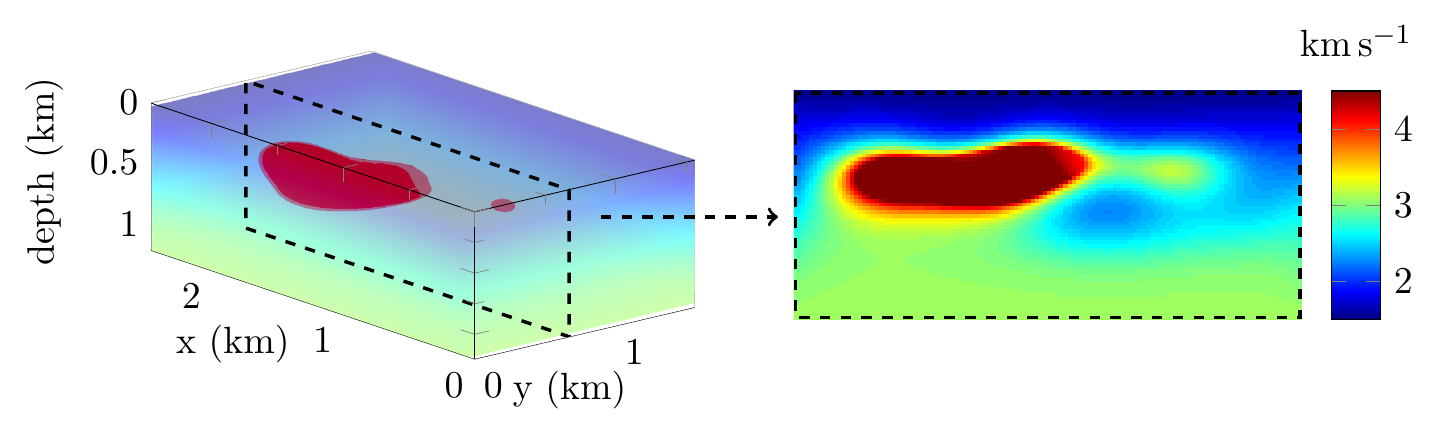}
\caption{Velocity model reconstruction  with 
         vertical section at $y=670$ \si{\meter},
         from the starting medium Figure~\ref{fig:3dsalt_start}
         using $5$ \si{\Hz} frequency data.
         The reconstruction applies the eigenvector 
         decomposition with $\eta_3$ and 
         $N=50$ with $30$ iterations, see 
         Algorithm~\ref{algo:fwi}.}
\label{fig:3dsalt_fwi_geman_1N}
\end{figure}

\begin{figure}
 \centering
 \renewcommand{\modelfile}{cp_5hz_start20ev-100ev_geman}
 \includegraphics[scale=1]{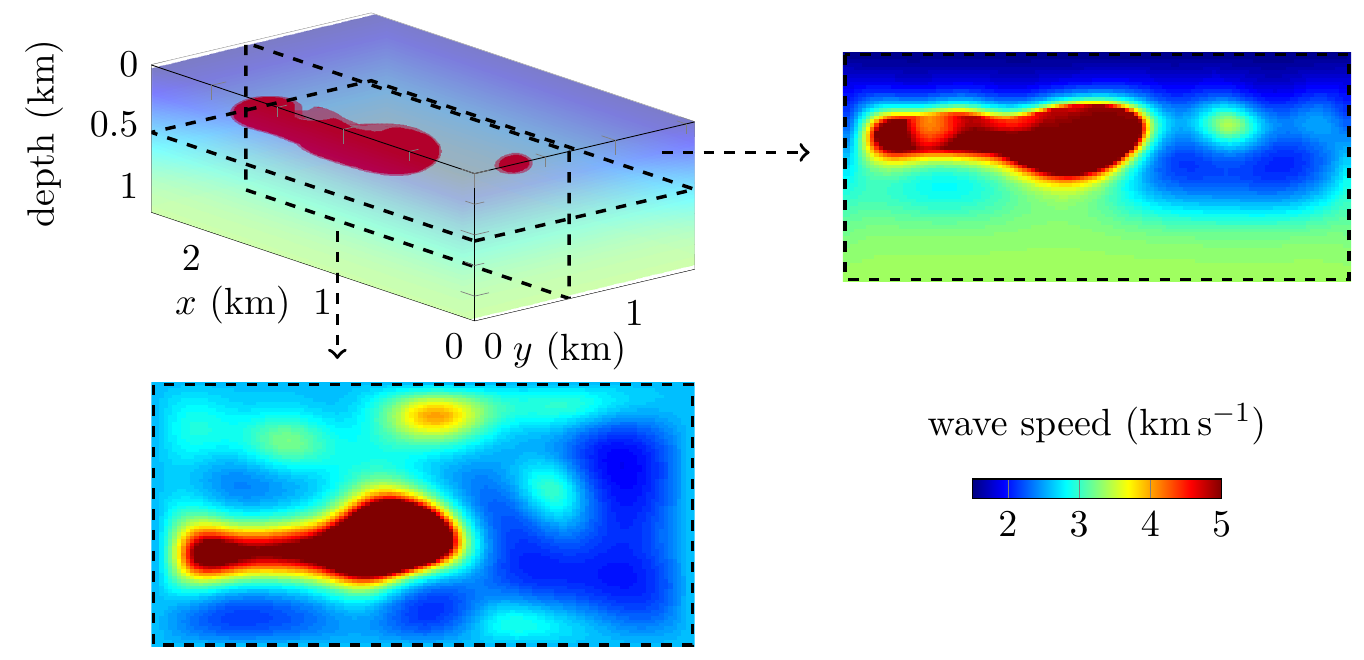}
\caption{Velocity model reconstruction  with 
         vertical section at $y=670$ \si{\meter}
         and horizontal section at $550$ \si{\meter} 
         depth, from the starting medium 
         Figure~\ref{fig:3dsalt_start} and
         using $5$ \si{\Hz} frequency data.
         The reconstruction applies the eigenvector 
         decomposition with $\eta_3$ and 
         progression of $N=\{20,30,50,75,100\}$,
         with \num{30} iterations per $N$, see 
         Algorithm~\ref{algo:fwi}.}
 \label{fig:3dsalt_fwi_geman_multiN}
\end{figure}

In Figure~\ref{fig:3dsalt_fwi_classic}, 
we show the reconstruction without 
employing the eigenvector decomposition, 
where the wave speed has a piecewise 
constant representation on a 
$124 \times 79 \times 61$ nodal grid. 
We only use $5$ \si{\Hz} frequency 
data, and $30$ iterations.
Next, we employ the eigenvector model 
representation with Algorithm~\ref{algo:fwi}.
Following the discussion in Subsection~\ref{subsection:implementation}, 
we select $\eta_3$, which is the most robust, and try two situations: 
\begin{itemize}
  \item \emph{single frequency} ($5$ \si{\Hz}),
        \emph{fixed $N$} reconstruction using 
        $\eta_3$, $N=50$ and $30$ iterations,
        the final reconstruction is Figure~\ref{fig:3dsalt_fwi_geman_1N}; 
  \item \emph{single frequency} ($5$ \si{\Hz}),
        \emph{multiple $N$} reconstruction using 
        $\eta_3$, $N=\{20,30,50,75,100\}$ and 
        $30$ iterations per $N$, i.e. $180$ 
        iterations in total: the final reconstruction is 
        Figure~\ref{fig:3dsalt_fwi_geman_multiN}.
\end{itemize}

\setlength{\modelwidth}{7.1cm}
\graphicspath{{figures/salt3d_fwi/1basis_1freq_multiN_5hz/}}

\setlength{\modelwidth}{7.1cm}
\graphicspath{{figures/salt3d_fwi/1basis_multi-freq_1N/}}

\graphicspath{{figures/salt3d_fwi/1basis_1freq_multiN_5hz/}}
\setlength{\modelwidth}{7.1cm}

For visualization, we focus on the two-dimensional 
vertical and horizontal sections which illustrate the 
positions and shapes of the objects. 
This experiment is consistent with the
two-dimensional results, and we observe 
the following.
\begin{itemize}
  \item The classical FWI reconstruction fails 
        to discover the subsurface objects, with
        only a narrow layer, misplaced, see 
        Figure~\ref{fig:3dsalt_fwi_classic}.
  \item The reconstruction with the representation 
        with eigenvector 
        decomposition is able to accurately 
        capture the largest subsurface salt 
        domes, see Figures~\ref{fig:3dsalt_fwi_geman_1N}
        and~\ref{fig:3dsalt_fwi_geman_multiN}.
  \item The best results are obtained when we 
        use the progression of $N$, with a single 
        frequency, see Figure~\ref{fig:3dsalt_fwi_geman_multiN}.
        The main salt dome is captured and smaller 
        ones start to appear, including near the 
        boundary. Using  a single $N$ is also a 
        good candidate, Figure~\ref{fig:3dsalt_fwi_geman_1N}, 
        as it necessitates much less iterations 
        ($30$ instead of $180$) hence less computational time.
\end{itemize}

\section{Conclusion}

In this paper, we have investigated the use of 
a representation based upon a basis of eigenvectors 
for image decomposition and quantitative reconstruction
in the seismic inverse problem, 
using two and three-dimensional experiments. 
We have implemented several diffusion coefficients,
and compared their performance, depending on the 
target medium.
\begin{enumerate}
  \renewcommand{\theenumi}{(\arabic{enumi})}
  \item In the context of image decomposition, 
        the case of contrasting objects (salt 
        domes) is clearly more appropriate 
        than the layered media (such as Marmousi). 
        All of the diffusion coefficients behave well and
        provide satisfactory results for salt domes
        image decomposition, 
        even in the presence of noise. 
  \item For the decomposition of images with layered patterns, 
        only a few formulations perform well ($\eta_1$, $\eta_3$, $\eta_6$).
        It would be interesting to investigate 
        further the performance of anisotropic or directional
        diffusion coefficients, mentioned by~\cite{Weickert1998,Grote2019}.
\end{enumerate}

Next, we have considered the quantitative 
reconstruction procedure in seismic, where 
only partial, backscattered (i.e. reflection)
data are available, from one side illumination.
We have probed the performance of the method 
by considering initial guesses with minimal a priori
information, and by avoiding the low frequency 
data, which are not accessible in seismic applications.
In this context, the FWI algorithm based upon 
the eigenvector model representation has shown 
promising results of subsurface $3$D 
salt dome media.
The method only requires the preliminary computation
of the basis of eigenvectors associated with the diffusion
operator, and a trivial modification of the gradient 
computation. Namely, the overall additional cost of the method
remains marginal compared to the cost of FWI.
Our findings are the following.
\begin{enumerate} \setcounter{enumi}{2}
  \renewcommand{\theenumi}{(\arabic{enumi})}
  \item For reconstruction, the result depends 
        on the choice of diffusion coefficient.
        We recommend $\eta_3$, from~\cite{Geman1992},
        which was the most robust in our applications, 
        even with a fixed $N$ and a few iterations.
  \item We have shown that the choice of $N$
        is not trivial, and one cannot rely on the misfit
        functional evaluation. Therefore, we have 
        proposed a progression of increasing $N$, which
        appears to stabilize the reconstruction.
  \item Because the method has a smoothing effect, 
        it focuses on the long wavelength structures.
        Thus, the reconstruction using the decomposition 
        can serve as an initial model to iterate with 
        higher frequency contents, in order to improve 
        the resolution.
\end{enumerate}
Following these analyses, some difficulties 
remain regarding the optimal choice of 
parameters. 
For instance, it is possible that $\eta$ has 
to be selected differently depending on the 
model (as illustrated with the Marmousi 
decomposition). 
Similarly, the scaling coefficient $\beta$ 
would affect the performance but the acceptable 
range appears, hopefully, quite large.

For future work, it seems that the lack of background 
velocity information would not be overcome by the 
decomposition, possibly resulting in artifacts.
Therefore, we envision the use of multiple basis to 
parametrize the velocity 
(e.g., using the background/reflectivity decomposition
idea of \cite{Clement2001,Faucher2019IP}, with a dedicated
smooth eigenvector basis to represent the background, and 
another to represent the reflectors).

Eventually, the method can readily extend to multi-parameter
inversion (e.g. for elastic medium), upon taking a separate basis
per parameter (e.g. one for each of the Lam\'e parameters
in linear elasticity). However, a more appropriate approach would 
be the use of joint-basis by, e.g., considering a system of PDE 
instead of the scalar diffusion operator. 
We have in mind strategies such that joint-sparsity and tensor decomposition.

\section*{Acknowledgments}
The authors would like to thank Prof. Marcus 
Grote for thoughtful comments and discussions.
The research of FF is supported by the 
Inria--TOTAL strategic action DIP.
FF acknowledges funding from the
Austrian Science Fund (FWF) under the 
Lise Meitner fellowship M 2791-N.
OS is supported by the FWF, 
with SFB F68, project F6807-N36 
(Tomography with Uncertainties).
OS also acknowledges support from the FWF via 
the project I3661-N27 (Novel Error Measures and Source 
Conditions of Regularization Methods for 
Inverse Problems).
HB acknowledges funding from the European 
Union's Horizon 2020 research and innovation 
program under the Marie Sklodowska-Curie grant
agreement Number 777778 (Rise action Mathrocks)
and with the E2S--UPPA CHICkPEA project.

\appendix
\renewcommand\thefigure{\thesection\arabic{figure}} 
\setcounter{figure}{0}
\section{Two-dimensional reconstruction with evolution of basis}
\label{appendix:multi_basis_experiments}

In this appendix, we experiment the re-computation 
of the eigenvectors basis along with the iterations. 
We consider two additional experiments, which 
are derived from Experiments~1 and 2 of 
Table~\ref{table:fwi_experiments_2Dsalt}, with one
major difference: 
\emph{the set of eigenvectors is recomputed from 
the current iteration model every $30$ iterations}. 
This is advocated in \cite{DeBuhan2013,Grote2019} where,
contrary to our experiments, the background velocity
is mostly known.

\begin{figure}[ht!] \centering
  \includegraphics[scale=1]{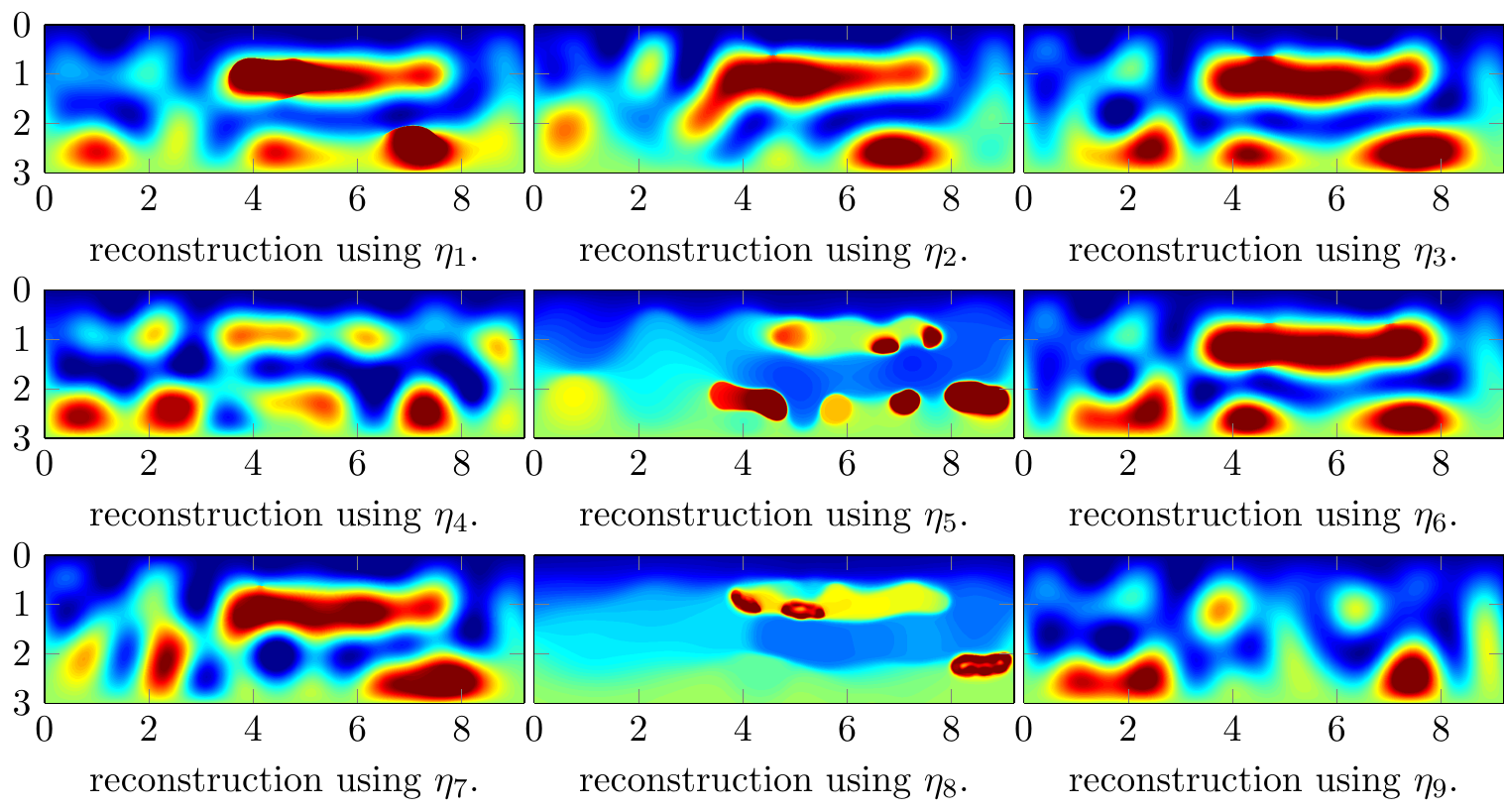}
  \caption{Reconstruction after $180$ iterations 
           from starting medium Figure~\ref{fig:fwi_start_salt}
           (same color scale) using $2$ \si{\Hz} data.
           The decomposition uses
           $N=50$ and the set of eigenvectors
           is re-computed every $30$ iterations,
           from the current iteration model.
           The formulations of $\eta$ are 
           Table~\ref{table:list_of_formula}.}
  \label{fig:fwi_salt_multibasis_single-freq_fixedN}
\end{figure}

\begin{figure}[ht!] \centering
  \includegraphics[scale=1]{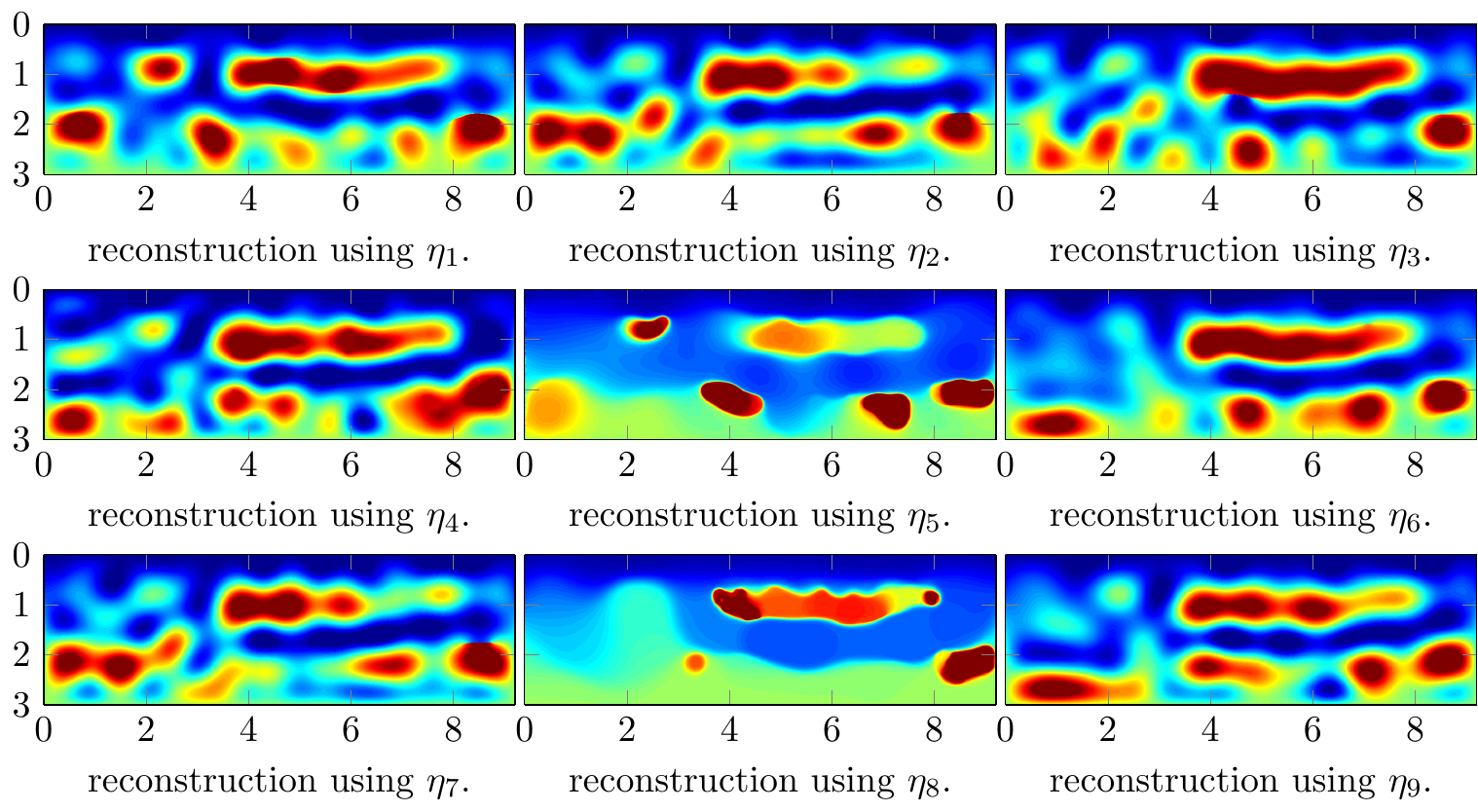}
  \caption{Reconstruction 
           from starting medium Figure~\ref{fig:fwi_start_salt}
           (same color scale)           
           using $2$ \si{\Hz} data.
           The decomposition employs the sequence
           $N=\{50,60,70,80,90,100\}$ with $30$
           iterations per $N$.
           The set of eigenvectors
           is re-computed every $30$ iterations,
           from the current iteration model
           (i.e. every time we update $N$).
           The formulations of $\eta$ are from 
           Table~\ref{table:list_of_formula}.}
  \label{fig:fwi_salt_multibasis_single-freq_multiN}
\end{figure}

Therefore, for the $180$ total iterations, 
the basis are 
\begin{itemize}
  \item computed from the initial model when starting 
        the very first iteration,
  \item re-computed from the current reconstruction at 
        the beginning of iterations $31$, $61$, $91$, $121$, $151$.
\end{itemize}
Compared to Algorithm~\ref{algo:fwi}, instead of using
a fixed $\basisloc$ from the initial model, it is 
recomputed from the current $\decomp^{(k)}$.
The reconstructions using the update of basis are shown
Figures~\ref{fig:fwi_salt_multibasis_single-freq_fixedN}
and~\ref{fig:fwi_salt_multibasis_single-freq_multiN},
for $N=50$ and $N=\{50,60,70,80,90,100\}$ respectively. 

\setlength{\modelwidth} {6.50cm}
\setlength{\modelheight}{3.10cm}
\graphicspath{{figures/salt_fwi_dataHnoise20db/21_basis-30iter_single-freq_fixed-N/}}
\renewcommand{\modelfileA} {cp_2hz_50ev_180iter}  
\renewcommand{\modelfileB} {cp_2hz_50ev_180iter}  
\renewcommand{\modelfileC} {cp_2hz_50ev_180iter}  
\renewcommand{\modelfileD} {cp_2hz_50ev_180iter}  
\renewcommand{\modelfileE} {cp_2hz_50ev_180iter}  
\renewcommand{\modelfileF} {cp_2hz_50ev_180iter}  
\renewcommand{\modelfileG} {cp_2hz_50ev_180iter}  
\renewcommand{\modelfileH} {cp_2hz_50ev_180iter}  
\renewcommand{\modelfileI} {cp_2hz_50ev_180iter}  

\setlength{\modelwidth} {6.50cm}
\setlength{\modelheight}{3.15cm}
\graphicspath{{figures/salt_fwi_dataHnoise20db/22_basis-30iter_single-freq_multi-N/}}
\renewcommand{\modelfileA} {cp_2hz-start50ev_180iter-100ev}  
\renewcommand{\modelfileB} {cp_2hz-start50ev_180iter-100ev}  
\renewcommand{\modelfileC} {cp_2hz-start50ev_180iter-100ev}  
\renewcommand{\modelfileD} {cp_2hz-start50ev_180iter-100ev}  
\renewcommand{\modelfileE} {cp_2hz-start50ev_180iter-100ev}  
\renewcommand{\modelfileF} {cp_2hz-start50ev_180iter-100ev}  
\renewcommand{\modelfileG} {cp_2hz-start50ev_180iter-100ev}  
\renewcommand{\modelfileH} {cp_2hz-start50ev_180iter-100ev}  
\renewcommand{\modelfileI} {cp_2hz-start50ev_180iter-100ev}  

We observe that changing the basis along with the 
iterations provides similar accuracy of the 
reconstruction compared to keeping the basis fixed.
For some formulations, the salt dome appears slightly
larger than with the fixed basis but the artifacts
at the bottom are also stronger, with patches of
high velocities.

We believe that the difficulty is due to the lack
of prior information regarding the velocity background.
Namely, the starting model has a background profile 
of low amplitude compared to the target model, 
see Figure~\ref{fig:fwi_start_salt}.
While our reconstruction captures the contrasting 
object, the background remains mostly erroneous
(i.e. $\mathfrak{m}_0$ in~\eqref{eq:methodEV_decomposition}),
impacting the decomposition.
Then, when we update the decomposition, the new set 
of eigenvectors will try to encompass the background 
variation, creating the bottom artifacts.

\setcounter{figure}{0}
\section{Extended set of figures}
\label{appendix:additional_figures}

We provide here the reconstructions obtained for
all $\eta$ for the two-dimensional salt dome FWI
test case of Subsection~\ref{subsection:fwi_2Dsalt_reconstruction}.

\graphicspath{{figures/salt_fwi_dataHnoise20db/01_fixed-basis_single-freq_fixed-N/start2hz/}}
\renewcommand{\modelfileA} {cp_2hz-30iter_100ev_perona1}  
\renewcommand{\modelfileB} {cp_2hz-30iter_100ev_perona2}  
\renewcommand{\modelfileC} {cp_2hz-30iter_100ev_geman}    
\renewcommand{\modelfileD} {cp_2hz-30iter_100ev_green}    
\renewcommand{\modelfileE} {cp_2hz-30iter_100ev_aubert}   
\renewcommand{\modelfileF} {cp_2hz-30iter_100ev_lorentz}  
\renewcommand{\modelfileG} {cp_2hz-30iter_100ev_gaussian} 
\renewcommand{\modelfileH} {cp_2hz-30iter_100ev_rudin}    
\renewcommand{\modelfileI} {cp_2hz-30iter_100ev_tikhonov} 
\begin{figure}[ht!] \centering
  \includegraphics[scale=1]{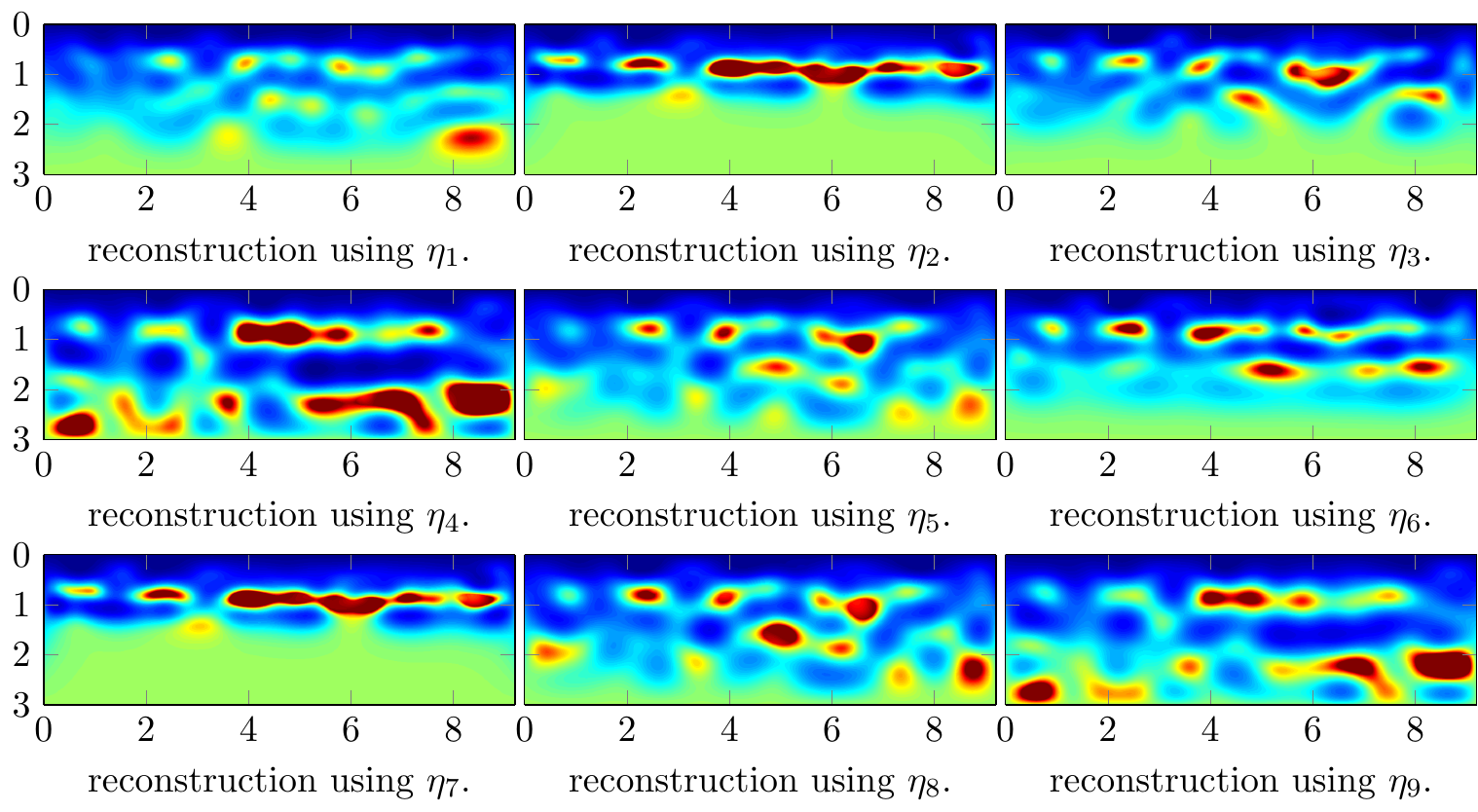}
  \caption{Reconstruction of the salt velocity model 
           from starting medium Figure~\ref{fig:fwi_start_salt}
          (with the same color scale)           
           using $2$ \si{\Hz} frequency data.          
           The eigenvector decomposition employs $N=100$ 
           and the formulations of $\eta$ from 
           Table~\ref{table:list_of_formula}.}
  \label{fig:fwi_salt_100ev_2hz}
\end{figure}

\graphicspath{{figures/salt_fwi_dataHnoise20db/02_fixed-basis_single-freq_multi-N/start2hz/}}
\renewcommand{\modelfileA} {cp_2hz_50ev_perona1_180iter-100ev}  
\renewcommand{\modelfileB} {cp_2hz_50ev_perona2_180iter-100ev}  
\renewcommand{\modelfileC} {cp_2hz_50ev_geman_180iter-100ev}    
\renewcommand{\modelfileD} {cp_2hz_50ev_green_180iter-100ev}    
\renewcommand{\modelfileE} {cp_2hz_50ev_aubert_180iter-100ev}   
\renewcommand{\modelfileF} {cp_2hz_50ev_lorentz_180iter-100ev}  
\renewcommand{\modelfileG} {cp_2hz_50ev_gaussian_180iter-100ev} 
\renewcommand{\modelfileH} {cp_2hz_50ev_rudin_180iter-100ev}    
\renewcommand{\modelfileI} {cp_2hz_50ev_tikhonov_180iter-100ev} 
\begin{figure}[ht!] \centering
  \includegraphics[scale=1]{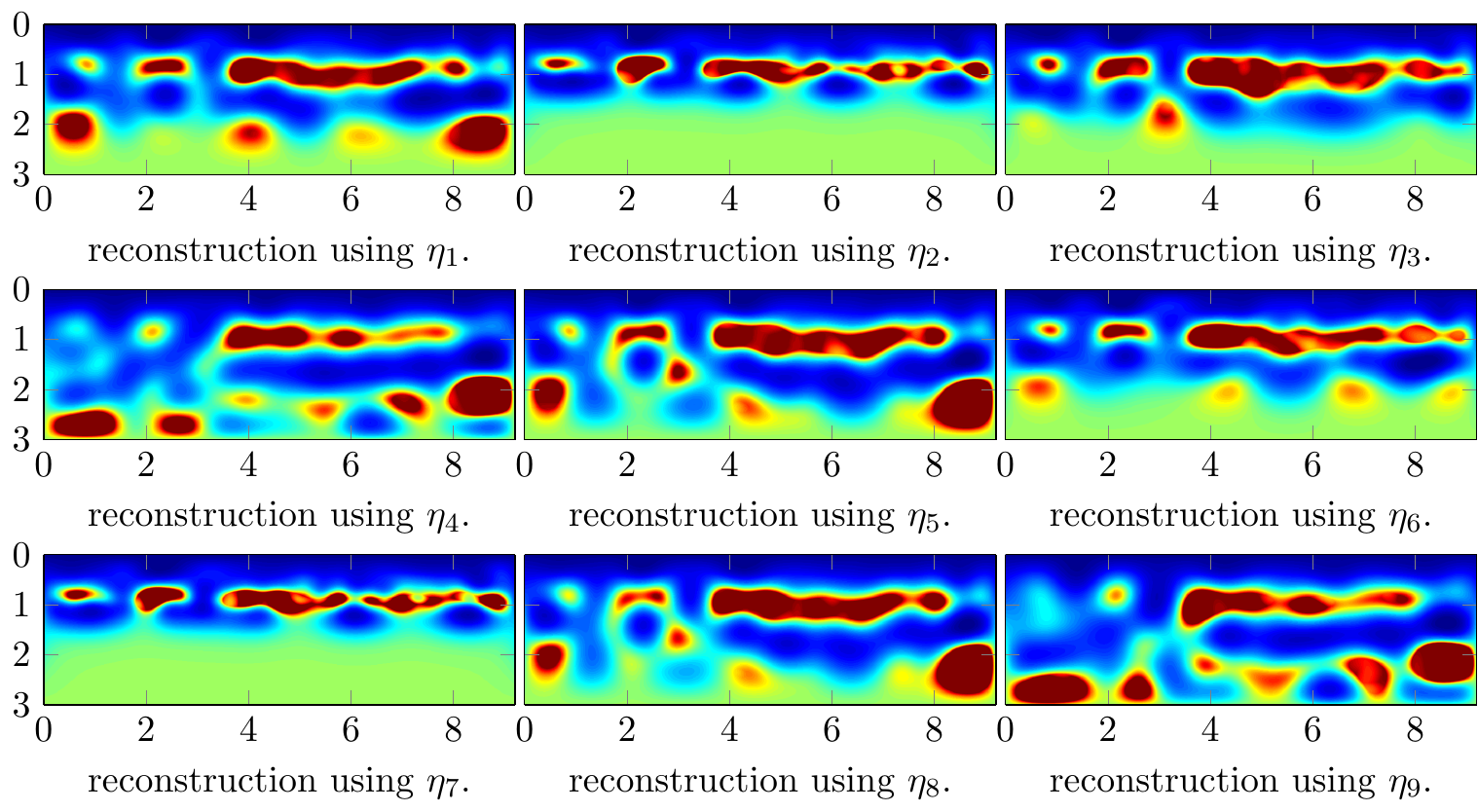}
  \caption{Experiment~2 of Table~\ref{table:fwi_experiments_2Dsalt}
           (multiple $N$, single frequency data)
           for the reconstruction of the salt velocity model 
           from starting medium Figure~\ref{fig:fwi_start_salt}
          (with the same color scale),
           the formulations of $\eta$ follow 
           Table~\ref{table:list_of_formula}.}
  \label{fig:fwi_salt_multiN_single-freq_2hz}
\end{figure}

\graphicspath{{figures/salt_fwi_dataHnoise20db/03_fixed-basis_multi-freq_fixed-N/start2hz/}}  
\renewcommand{\modelfileA} {cp_5hz_50ev}
\renewcommand{\modelfileB} {cp_5hz_50ev}
\renewcommand{\modelfileC} {cp_5hz_50ev}
\renewcommand{\modelfileD} {cp_5hz_50ev}
\renewcommand{\modelfileE} {cp_5hz_50ev}
\renewcommand{\modelfileF} {cp_5hz_50ev}
\renewcommand{\modelfileG} {cp_5hz_50ev}
\renewcommand{\modelfileH} {cp_5hz_50ev}
\renewcommand{\modelfileI} {cp_5hz_50ev}
\begin{figure}[ht!] \centering
  \includegraphics[scale=1]{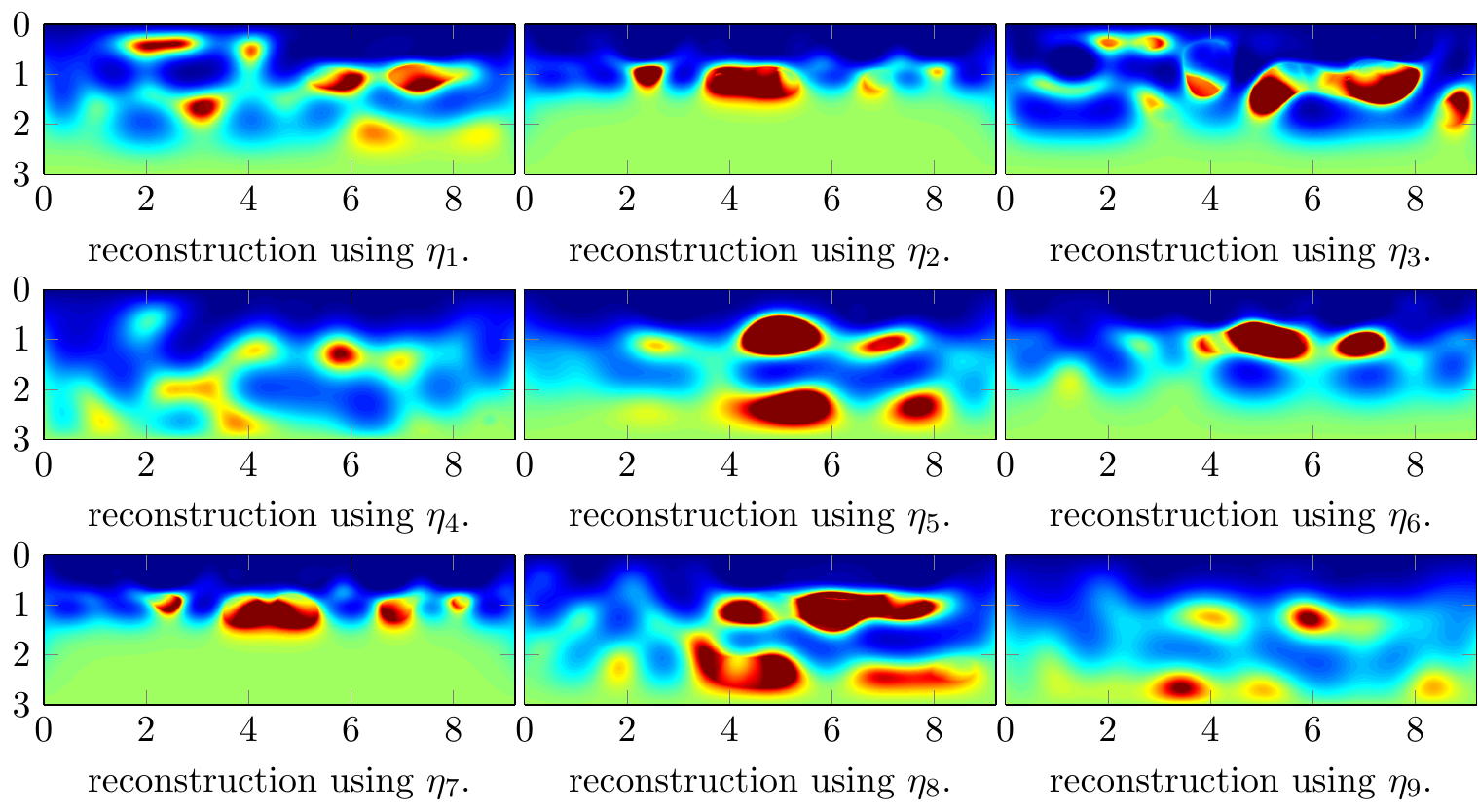}
  \caption{Experiment~3 of Table~\ref{table:fwi_experiments_2Dsalt}
           (single $N$, multiple frequency data)
           for the reconstruction of the salt velocity model 
           from starting medium Figure~\ref{fig:fwi_start_salt}
          (with the same color scale),
           the formulations of $\eta$ follow 
           Table~\ref{table:list_of_formula}.}
  \label{fig:fwi_salt_50ev_multi-freq_5hz}
\end{figure}

\graphicspath{{figures/salt_fwi_dataHnoise20db/04_fixed-basis_multi-freq_multi-N/start2hz/}}
\renewcommand{\modelfileA} {cp_start50ev_5hz-80ev}  
\renewcommand{\modelfileB} {cp_start50ev_5hz-80ev}  
\renewcommand{\modelfileC} {cp_start50ev_5hz-80ev}  
\renewcommand{\modelfileD} {cp_start50ev_5hz-80ev}  
\renewcommand{\modelfileE} {cp_start50ev_5hz-80ev}  
\renewcommand{\modelfileF} {cp_start50ev_5hz-80ev}  
\renewcommand{\modelfileG} {cp_start50ev_5hz-80ev}  
\renewcommand{\modelfileH} {cp_start50ev_5hz-80ev}  
\renewcommand{\modelfileI} {cp_start50ev_5hz-80ev}  
\begin{figure}[ht!] \centering
  \includegraphics[scale=1]{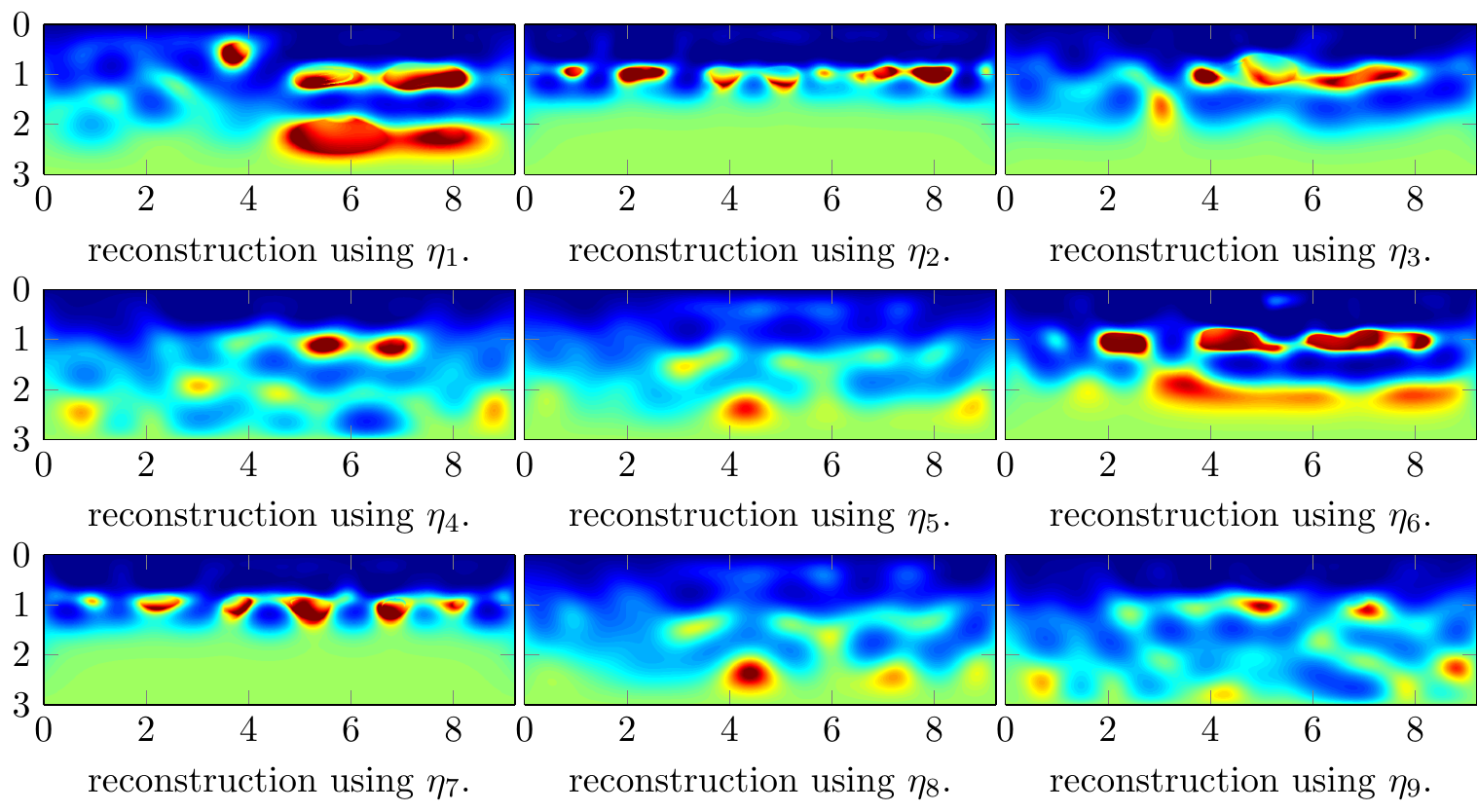}
  \caption{Experiment~4 of Table~\ref{table:fwi_experiments_2Dsalt}
          (multiple $N$, multiple frequency data)
           for the reconstruction of the salt velocity model 
           from starting medium Figure~\ref{fig:fwi_start_salt}
          (with the same color scale),
           the formulations of $\eta$ follow 
           Table~\ref{table:list_of_formula}.}
  \label{fig:fwi_salt_multiN_multi-freq_5hz}
\end{figure}


\bibliographystyle{gji}
\bibliography{bibliography}

\end{document}